\documentclass[nohyperref]{article}

\usepackage{color}

\usepackage{ulem} % for \erase

\usepackage{microtype}
\usepackage{graphicx}
\usepackage{subfigure}
\usepackage{booktabs} % for professional tables

\usepackage{hyperref}

\newcommand{\E}{\mathbb{E}}

\usepackage{arxiv}

\usepackage{amsmath}
\usepackage{amssymb}
\usepackage{mathtools}
\usepackage{amsthm}
\usepackage[capitalize,noabbrev]{cleveref}

\theoremstyle{plain}
\newtheorem{theorem}{Theorem}[section]

\newtheorem{lemma}[theorem]{Lemma}

\theoremstyle{definition}
\newtheorem{definition}{Definition}[section]
\theoremstyle{remark}

\theoremstyle{definition}
\newtheorem{custom}{Assumption}
 \newenvironment{assumption}[1]
  {\renewcommand\thecustom{#1}\custom}
  {\endcustom}

\usepackage[textsize=tiny]{todonotes}

\RequirePackage{etoolbox}
\RequirePackage{xcolor}
\definecolor{ptblue}{RGB}{68, 119, 170}
\definecolor{ptcyan}{RGB}{102, 204, 238}
\definecolor{ptgreen}{RGB}{34, 136, 51}
\definecolor{ptyellow}{RGB}{204, 187, 68}
\definecolor{ptred}{RGB}{238, 102, 119}
\definecolor{ptpurple}{RGB}{170, 51, 119}
\definecolor{ptgrey}{RGB}{187, 187, 187}

\usepackage{xcolor}
\hypersetup{
    colorlinks=true,
    citecolor=ptgreen,
    linkcolor=ptpurple,
    urlcolor=ptblue,
}
\newcommand*{\cit}[1]{{[\cite{#1}]}}

\begin{document}

%\twocolumn[
%\icmltitle{Single Loop Gaussian Homotopy Method for Non-convex Optimization}
\title{Single Loop Gaussian Homotopy Method for\\ Non-convex Optimization}

\author{%
  Hidenori Iwakiri\thanks{The first two authors contributed equally.} \\
  The University of Tokyo, RIKEN AIP\\
  \texttt{iwakiri-hidenori2020@g.ecc.u-tokyo.ac.jp}\\
  % examples of more authors
  \And
  Yuhang Wang\textcolor{ptpurple}{\footnotemark[1]} \\
  The University of Tokyo\\
  \texttt{utyuuhikou@gmail.com}\\
  \AND
  Shinji Ito\\
  NEC Corporation, RIKEN AIP\\
  \texttt{i-shinji@nec.com}\\
  \And
  Akiko Takeda\\
  The University of Tokyo, RIKEN AIP\\
  \texttt{takeda@mist.i.u-tokyo.ac.jp}
}

\maketitle

\begin{abstract}
  The Gaussian homotopy (GH) method is a popular approach to finding better stationary points for non-convex optimization problems by gradually reducing a parameter value $t$, which changes the problem to be solved from an almost convex one to the original target one. Existing GH-based methods repeatedly call an iterative optimization solver to find a stationary point every time $t$ is updated, which incurs high computational costs. We propose a novel single loop framework for GH methods (SLGH) that updates the parameter $t$ and the optimization decision variables at the same. Computational complexity analysis is performed on the SLGH algorithm under various situations: either a gradient or gradient-free oracle of a GH function can be obtained for both deterministic and stochastic settings. The convergence rate of SLGH with a tuned hyperparameter becomes consistent with the convergence rate of gradient descent, even though the problem to be solved is gradually changed due to $t$. In numerical experiments, our SLGH algorithms show faster convergence than an existing double loop GH method while outperforming gradient descent-based methods in terms of finding a better solution.
\end{abstract}


\section{Introduction}
Let us consider the following non-convex optimization problem:
\begin{align}
  {\mathop{\rm minimize}\limits_{x\in {\mathbb{R}}^d}}\quad f(x), % := \mathbb{E}_\xi [F(x;\xi)]
  \label{nonconvprob}
\end{align}
where $f:{\mathbb{R}}^d\rightarrow \mathbb{R}$ is a non-convex function.
Let us also consider the following stochastic setting:
\begin{align}
  f(x) := \mathbb{E}_\xi [\bar{f}(x;\xi)],
  \label{stocprob}
\end{align}
where $\xi$ %\memo{the size $\xi \in \mathbb{R}^p$ or something?}
is the random variable following a probability distribution $P$ from which
i.i.d.~samples can be generated. % and the function $\bar{f}:{\mathbb{R}}^d\rightarrow \mathbb{R}$ is not convex.
Such optimization problems attract significant attention in machine learning,
and at the same time, the need for optimization algorithms that can find a stationary point with smaller objective value is growing. For example, though %the loss function in deep learning is generally considered to have a highly non-convex structure,
it is often said that simple gradient methods can find global minimizers for deep learning
(parameter configurations with zero or near-zero training loss),
such beneficial behavior is not universal, as noted in \cit{li2017visualizing};
the trainability of neural nets is highly dependent on network architecture design choices,
variable initialization, etc.
There are also various other highly non-convex optimization problems in machine learning
(see e.g., \cit{jain2017nonconv}).

The Gaussian homotopy (GH) method is designed to avoid
poor stationary points by building a sequence of
successively smoother approximations of the original objective function $f$,
and it is expected to find a good stationary point with a small objective value
for a non-convex problem.
More precisely, using the GH function $F(x,t)$ with a parameter $t \geq 0$ that satisfies
$F(x,0)=f(x)$, the method starts from solving an almost convex smoothed function $F(x,t_1)$ with
some sufficiently large $t_1 \geq 0$ and gradually changes the optimization problem $F(x,t)$
to the original one $f(x)$ while decreasing the parameter $t$.
The homotopy method developed so far, then, consists of a double loop structure;
the outer loop reduces $t$, and the inner loop solves $\min_x F (x, t)$ for the fixed $t$.

\paragraph{Related research on the GH method}
  The GH method is popular owing to its ease of implementation
  and the quality of its obtained stationary points, i.e., their function values. %often superior empirical performance.
  The nature of this method was first proposed in \cit{blake1987visual},
  and it was then successfully applied in various fields,
  including computer vision \cit{Nielsen1993,brox2010large,Zach2018},
  physical sciences \cit{hemeda2012homotopy}
  and computation chemistry \cit{wu1996effective}.
  \cit{hazan2016graduated} introduces machine learning applications for the GH method,
  and an application to tuning hyperparameters of kernel ridge regression \cit{shao2019graduated}
  has recently been introduced.
  %and machine learning (e.g., see \cit{hazan2016graduated} and references therein).
  %The nature of this method is firstly proposed in \cit{blake1987visual} and applied to solve optimization problems  appearing in surface reconstruction \cit{Nielsen1993},\cit{wu1996effective}.
Although there have been recent studies on the GH function $F(x,t)$ \cit{mobahi2015link, mobahi2015theoretical, hazan2016graduated}, all existing GH methods use the double loop approach noted above. Moreover, to the best of our knowledge, there are no existing works that give theoretical guarantee for the convergence rate except for \cit{hazan2016graduated}. It characterizes a family of non-convex functions for which a GH algorithm converges to a global optimum and derives the convergence rate to an $\epsilon$-optimal solution. However, the family covers only a small part of non-convex functions, and it is difficult to check whether the required conditions are satisfied for each function. See Appendix~\ref{sec:related_work} for more discussion on related work.

% It derives the convergence rate to $\epsilon$-optimal solution, but the function class of analysis in their work called ``$\sigma$-nice'' covers only a very small part of non-convex functions. % there has been almost no significant change in the double loop framework of GH algorithms from the earlier one referred to as “Graduated Non-Convexity” in \cit{blake1987visual}. To the best of our knowledge, all existing GH methods use the double loop approach noted above.

%, although the GH function is used for different purposes; for example, for finding better local minima in \cit{mobahi2015link,hazan2016graduated} and making the nonsmooth function $f$ smoothed in \cit{chen2012smoothing,jin2018local}.
%  and are restricted for applications due to the weakness of traditional Gaussian Homotopy method.

\begin{table}[tb] \label{table:summary}
\caption{Each theorem shows the iteration complexity of SLGH with respect to $\epsilon$
and the dimension of input space $d$  to reach an $\epsilon$-stationary point in
the corresponding problem setting. ``const.~$\gamma$'' shows the complexity
when we treat the decreasing parameter $\gamma$ as a constant.
``tuned $\gamma$'' shows the lowest complexity of SLGH attained by
updating $t$ appropriately, which matches the complexity of
the standard first- or zeroth-order methods (see e.g., Theorem \ref{iter_determin}).
We also consider two cases of a zeroth-order setting:
``exact $f$'', in which we can query the exact or stochastic function value,
and ``err.~$f$'', in which we can only access the function value with bounded error.}
\centering
\begin{tabular}{c|ccc}
    &  1) first-order  &  \multicolumn{2}{c}{zeroth-order} \\\cline{3-4}
    &    & 2) exact $f$ & 3) err. $f$ \\ \hline
  a) deterministic~& Thm.~\ref{iter_determin}  & Thm.~\ref{thm:zo_determinstic} &  Thm.~\ref{thm:error_deterministic}\\
  const. $\gamma$ &  $O\left(\frac{d^{3/2}}{\epsilon^2}\right)$ & $O\left(\frac{d^{2}}{\epsilon^2}\right)$ & $O\left(\frac{d^{3}}{\epsilon^2}\right)$ \\ 
  tuned $\gamma$ & $O\left(\frac{1}{\epsilon^2}\right)$ & $O\left(\frac{d}{\epsilon^2}\right)$ & $O\left(\frac{d}{\epsilon^2}\right)$\\ \hline 
  b) stochastic~& Thm.~\ref{thm:iter_stochastic} & Thm.~\ref{thm:zo_stochastic} & Thm.~\ref{thm:error_stochastic} \\
  const.~$\gamma$ &   $O \left(
    \frac{d}{\epsilon^4} + \frac{d^{3/2}}{\epsilon^2}
\right)$ & $O\left(\frac{d^2}{\epsilon^4}\right)$ & $O\left(\frac{d^2}{\epsilon^4}+\frac{d^3}{\epsilon^2}\right)$\\
  tuned $\gamma$ & $O\left(\frac{1}{\epsilon^4}\right)$ &  $O\left(\frac{d}{\epsilon^4}\right)$ & $O\left(\frac{d}{\epsilon^4}\right)$\\ \hline
\end{tabular}
\end{table}

\paragraph{Motivation for this work}
This paper proposes novel deterministic and stochastic GH methods employing
a single loop structure in which the decision variables $x$
and the smoothing parameter $t$ are updated at the same time
using individual gradient/derivative information.
Using a well-known fact in statistical physics on the relationship between
the {\it heat equation} and Gaussian convolution of $f$,
together with {\it the maximum principle} (e.g., \cit{evans2010partial})
for the {\it heat equation},
we can see that a solution $(x^\ast, t^\ast)$
minimizing the GH function $F(x,t)$ satisfies $t^\ast=0$;
thus, $x^\ast$ is also a solution for (\ref{nonconvprob}).
This observation leads us to a single loop GH method (SLGH, in short),
which updates the current point $(x_k,t_k)$ simultaneously for $\min_{x \in \mathbb{R}^d, t \geq 0} F(x,t)$.
The resulting SLGH method can be regarded as an application of
the steepest descent method to the optimization problem,
with $(x, t)$ as a variable. We are then able to investigate the convergence rate of
our SLGH method so as to achieve an $\epsilon$-stationary point of
\eqref{nonconvprob} and \eqref{stocprob} by following existing theoretical complexity analyses.

%Table~\ref{table:summary} summarizes  the convergence rate of our SLGH method under several problems settings shown below.

We propose two variants of the SLGH method: $\text{SLGH}_\text{d}$ and $\text{SLGH}_\text{r}$, which have different update rules for $t$. $\text{SLGH}_\text{d}$ updates $t$ using the derivative of $F(x, t)$ in terms of $t$, based on the idea of viewing $F(x, t)$ as the objective function with respect to the variable $(x, t)$. Though this approach is effective in finding good solutions (as demonstrated in Appendix \ref{subsec:toy}), it requires additional computational cost due to the calculation of $\frac{\partial F}{\partial t}$. To avoid this additional computational cost, we also consider $\text{SLGH}_\text{r}$ that uses fixed-rate update rule for $t$. We also show that both $\text{SLGH}_\text{d}$ and $\text{SLGH}_\text{r}$ have the same theoretical guarantee.

Table 1 summarizes the convergence rate of our SLGH method to reach
an $\epsilon$-stationary point under a number of problem settings.
Since the convergence rate depends on the decreasing speed of $t$,
we list two kinds of complexity in the table; details are described in the caption.

We consider the three settings in which available oracles differ. In Case 1),
the full (or stochastic) gradient of $F(x,t)$ in terms of $x$ is available for
the deterministic problem \eqref{nonconvprob} (or stochastic problem \eqref{stocprob},
respectively). However, in this setting, we have to calculate Gaussian convolution
for deriving GH functions and their gradient vectors, which becomes expensive,
especially for high-dimensional applications, unless closed-form expression of Gaussian convolution is
possible. While \cit{mobahi2016closed} provides closed-form expression for some specific functions $f$, such as polynomials, Gaussian RBFs, and trigonometric functions,
such problem examples are limited.
% a closed-form smoothing convolution, which is expensive especially in high-dimension applications.
As Case 2), we extend our deterministic and stochastic GH methods
to the zeroth-order setting, for which the convolution computation is approximated
using only the function values. %% Zeroth-order optimization has become increasingly popular due to its wide application to situations in which the explicit gradient cannot be calculated, and in which often only the exact function values can be queried.
%Such a class of applications appear in black-box adversarial attacks on deep neural networks \cit{chen2019zo}, structured prediction \cit{sokolov2016stochastic} and reinforcement learning \cit{xu2020zeroth}. Many zeroth-order methods (ZOSGD \cit{ghadimi2013stochastic}, ZOAdaMM \cit{chen2019zo}, ZOSVRG \cit{liu2018zeroth} have been proposed for such black-box situation. All of them are developed from ZOGD in \cit{nesterov2017random}, which introduces random gradient-free oracles based on Gaussian smoothing.
%% This trend also applies to research on the GH method; \cit{hazan2016graduated} first developed a zeroth order variant of the GH method, and \cit{shao2019graduated} used a modified one for hyperparameter tuning. Their methods are zeroth-order double loop GH methods, while ours adopts a single loop structure. See Appendix~\ref{sec:related_work} for more discussion on machine learning applications and popular zeroth-order methods.
%although their smoothed gradient oracle is not very efficient \memo{AT}{ I do not check this but our method is batter than Hazan's work? Should we say ``slightly different setting''? the paper assumes that the objective is only accessible through a noisy value oracle.}.
%Therefore, we adopt their proposed GradOpt with the random gradient-free oracles in \cit{nesterov2017random} to ZOGradOpt and choose it as a baseline of zeroth-order GH algorithm in our experiment.
Another zeroth-order setting, Case 3), is also considered in this paper:
the inexact function values (more precisely, the function value with bounded error) can be
queried similarly as in the setting in \cit{jin2018local}. See Appendix \ref{sec:error_appendix} for more details.
%In order to better compare our proposed algorithms with existing methods, we conducted two experiments. The first toy example is constructed by adding noise to a simple convex function. And the non-convex landscapes of such artificial non-convex example can be adjusted with different noises. The second experiment is for generating black-box adversarial attack examples.

Although no existing studies have analyzed the complexity of a double loop GH method to find an $\epsilon$-stationary point, we can see that its inner loop requires the same complexity as GD (gradient descent) method up to constants. Furthermore, as noted above, the complexity of the SLGH method with a tuned hyperparameter matches that of GD method. Thus, the SLGH method becomes faster than a double loop GH method by around the number of outer loops. The SLGH method is also superior to double loop GH methods from practical perspective, because in order to ensure convergence of their inner loops, we have to set the stepsize conservatively, and furthermore a sufficiently tuned terminate condition must be required.

\textbf{Contributions} \quad We can summarize our contribution as follows:

(1) We propose novel deterministic and stochastic single loop GH (SLGH) algorithms
and analyze their convergence rates to an $\epsilon$-stationary point.
As far as we know, this is the first analysis of convergence rates of GH methods for general non-convex problems \eqref{nonconvprob} and \eqref{stocprob}. For non-convex optimization, the convergence rate of SLGH with a tuned hyperparameter becomes consistent with the convergence rate of gradient descent, even though the problem to be solved is gradually changed due to $t$.
% For smooth non-convex optimization, the complexity with respect to $\epsilon$ matches that of the standard gradient descent (GD) method. The complexity with respect to $d$ depends on the decreasing speed of the smoothing parameter $t$, and the complexity matches that of the standard GD method if the decreasing speed is appropriate.
At this time, the SLGH algorithms become faster than a double loop one by around its number of outer loops.

(2) We propose zeroth-order SLGH (ZOSLGH) algorithms based on zeroth-order estimators of
gradient and Hessian values, which are useful when Gaussian smoothing convolution is
difficult. We also consider the possibly non-smooth case in which the accessible function
contains error, and we derive the upper bound of the error level for convergence guarantee.

(3) We empirically compare our proposed algorithm and other algorithms in experiments,
including artificial highly non-convex examples and black-box adversarial attacks.
Results show that the proposed algorithm converges much faster than an existing
double loop GH method, while it is yet able to find better solutions than are
GD-based methods.

\section{Standard Gaussian homotopy methods}
\paragraph{Notation:}
For an integer $N$, let $[N]:=\{1,...,N \}$. We express $\chi_{[N]}:=\{\chi_1, \ldots, \chi_N \}$ for a set of some vectors. We also express the range of the smoothing parameter $t$ as $\mathcal{T}:= [0, t_1]$, where $t_1$ is an initial value of the smoothing parameter. Let $\|\cdot\|$ denote the Euclidean norm and $\mathcal{N}(0,\mathrm{I}_d)$ denote the $d$-dimensional standard normal distribution.

Let us first define Gaussian smoothed function.
\begin{definition}
%  Given a function $f(x)$,
  Gaussian smoothed function $F(x,t)$ of $f(x)$ is defined as follows:
    \begin{align}
    F(x,t) &:=\mathbb{E}_{u\sim\mathcal{N}(0,\mathrm{I}_d)}[f(x+tu)]
    = \int f(x+ty)k(y)dy,
    \label{00}
    \end{align}
    where $k(y) = (2\pi)^{-d/2}\exp{(-{\|y\|^2}/2)}$ is referred to as the Gaussian kernel.
\end{definition}
The idea of Gaussian smoothing is to take an expectation over the function value with
a Gaussian distributed random vector $u$. For any $t>0$,
the smoothed function $F(x,t)$ is a $C^\infty$ function,
and $t$ plays the role of a smoothing parameter that controls the level of smoothing.

Here, let us show the link between Gaussian smoothing and
the {\it heat equation} \cit{widder1976heat}. The Gaussian smoothing convolution is
basically the solution of the {\it heat equation} \cit{widder1976heat}.
\begin{align}
    \frac{\partial}{\partial t}\hat{u} = \Delta_x \hat{u}, \quad \hat{u}(\cdot,0) = f(\cdot),
    \label{heatequation}
\end{align}
where $\Delta_x$ denotes the Laplacian. The solution of the {\it heat equation} is
$\hat{u}(x,t) = (\frac{1}{4\pi t})^{\frac{d}{2}}\int f(y) e^{-\frac{\|x-y\|^2}{4t}}dy$. %\memo{About $\hat{u}(x;t)$: ``;'' is standard?}
This can be made the same as the Gaussian smoothing function $F(x,t)$
by scaling its coefficient, which only changes the speed of progression.
%This property plays an important role in the proof of Theorem~\ref{thm: main}, shown in Appendix \ref{subsec:proof_optimality}. %the following theorem.

Corollary 9 in  \cit{mobahi2012gaussian} shows a sufficient condition for ensuring that
$f$ has the asymptotic strict convexity in which the smoothed function $F(x,t)$ becomes
convex if a sufficiently large smoothing parameter $t$ is chosen.
%The following theorem shows a sufficient condition to ensure that $f$ has the asymptotic strict convexity
%where the smoothed function $F(x,t)$ becomes convex if a sufficient large smoothing parameter is chosen.
%\begin{theorem}(Corollary 9 in \cit{MobahiMa})\label{thm:asymptoticConvex}
%  $f\in L^1(\mathbb{R}^n)$, $\int f(x)\, dx<0$  and $f(x)=O(\|x\|^{-(n+3)})$ as $\|x\|\to \infty$. Then, $f$ is asymptotically strictly convex.
%\end{theorem}
On this basis, the standard GH method, Algorithm~\ref{algo:Continuation}, starts with
a (almost) convex optimization problem $F(x,t)$ with
large parameter value $t \in \mathbb{R}$ and gradually changes the problem toward
the target non-convex $f(\cdot)=F(\cdot,0)$ by decreasing $t$ gradually.
\cit{hazan2016graduated} reduces $t$ by multiplying by a factor of $1/2$ for
each iteration $k$.  \cit{mobahi2015theoretical} focuses more on
theoretical work w.r.t.~the general setting and do not discuss the update rule for $t$.

%  can achieve a better solution by optimizing the above smoothed function $F(x;t)$ for the $k$th iteration ($k\in[T]$) with different smoothing parameter $t$. However, this double-loop structure method becomes much slower because the original optimization problem (1) is turned into $T$ sub-problems.

\begin{algorithm}[H]
    \caption{Standard GH method (\cit{mobahi2015theoretical,hazan2016graduated})} \label{algo:Continuation}
    \begin{algorithmic}
      \REQUIRE Objective function $f$, iteration number $T$, sequence %$(t_i)_{i \in [1,T]}$
      $\{t_1,\ldots,t_T\}$ satisfying $t_1 > \cdots > t_T$.
      Find a solution $x_1$ for minimizing $F(x, t_1)$.
      \FOR {$k = 1$ to $T$}
      \STATE Find a stationary point $x_{k+1}$ of $F(x, t_{k+1})$  with the initial solution $x_{k}$.
      \ENDFOR
      \RETURN $x_T$
    \end{algorithmic}
\end{algorithm}

\section{Single loop Gaussian homotopy algorithm}
\label{sec:first-order}
A function $h(x)$ is $L_{0}$-$Lipschitz$ with a constant $L_0$
if for any $x,y \in \mathbb{R}^d$, $|h(x)-h(y)|\leq L_0 \|x-y\|$ holds.
In addition, $h(x)$ is $L_{1}$-$smooth$ with
a constant $L_1$ if for any $x,y \in \mathbb{R}^d$, $\|\nabla h(x)-\nabla h(y)\|\leq L_1 \|x-y\|$ holds.
%We sometimes denote $L_0$ and $L_1$ of a function $h$ by $L_0(h)$ and $L_1(h)$ as necessary.
Let us here list assumptions for developing algorithms with convergence guarantee.

\begin{assumption}{A1}$\ $\label{A1}
\renewcommand{\labelenumi}{(\roman{enumi})}
\begin{enumerate}
    \item Objective function $f$ satisfies $\sup_{x\in\mathbb{R}^d}\mathbb{E}_u[|f(x+tu)|] < \infty$ (In the stochastic setting,  $f$ satisfies $\sup_{x\in\mathbb{R}^d,\xi}\mathbb{E}_u[|\bar{f}(x+tu;\xi)|] < \infty$). % and $\sup_{x\in\mathbb{R}^d}\mathbb{E}_\xi[|\bar{f}(x;\xi)|] < \infty$).
    \item The optimization problem (\ref{nonconvprob}) has an optimal value $f^\ast$.
    \item Objective function $f(x)$ is $L_{0}$-$Lipschitz$ and $L_{1}$-$smooth$ on $\mathbb{R}^d$ (In the stochastic setting, $\bar{f}(x;\xi)$ is $L_{0}$-$Lipschitz$ and $L_{1}$-$smooth$ on $\mathbb{R}^d$ in terms of $x$ for any $\xi$).
\end{enumerate}
\end{assumption}

Assumption (i) for making $F(x,t)$ well-defined and enabling to exchange the order of differentiation and integration, as well as Assumption (ii),
is mandatory for theoretical analysis with the GH method.
Assumption (iii) is often imposed for gradient-based methods.
This is a regular boundedness and smoothness assumption
in recent non-convex optimization analyses
(see e.g., \cit{NEURIPS2019_50a074e6, NEURIPS2020_0cb5ebb1, NEURIPS2019_b8002139}).

In the remainder of this section, we consider the nature of
the GH method and propose a more efficient algorithm, a SLGH algorithm.
We then provide theoretical analyses for our proposed SLGH algorithm.

\subsection{Motivation}
The standard GH algorithm needs to solve an optimization problem for
a given smoothing factor $t$ in each iteration and manually reduce $t$,
e.g., by multiplying some decreasing factor. To simplify this process,
we consider an alternative problem as follows:
\begin{align}
  %  {\mathop{\rm minimize}\limits_{x\in \mathbb{R}^d, t\in\mathcal{T}}}\quad F(x,t)
  {\mathop{\rm minimize}\limits_{x\in \mathbb{R}^d, t \in \mathcal{T}}} \quad F(x,t),
  \label{xtvar_prob}
\end{align}
where %$\mathcal{T} \subset {\{x\in\mathbb{R}|x\geq0\}}$ and
$F(x,t)$ is the Gaussian smoothed function of $f(x)$. This single loop structure can
reduce the number of iterations by optimizing $x$ and $t$ at the same time.

The following theorem is a (almost) special case of Theorem 6 in \cit{evans2010partial},\footnote{Although the assumptions in Theorem 3.1 are stronger than those in the theorem proved by Evans, the statement of ours is also stronger than that of his theorem, in a sense that our theorem guarantees that all optimal solutions satisfy $t=0$.}
which is studied in statistical physics but may not be well-known in machine learning and
optimization communities. This theorem shows that
the optimal solution of (\ref{xtvar_prob}) $(x^\ast, t^\ast)$ satisfies $t^\ast=0$, and thus $x^\ast$ is also a solution for (\ref{nonconvprob}).
%This enables us to provide convergence analyses for the SLGH method.
Therefore, we can regard $F(x,t)$ as an objective function in the SLGH method.

\begin{theorem}\label{thm: main}
  %$F(x,t) $ is the Gaussian smoothing function of $f:\mathbb{R}^n \rightarrow \mathbb{R}$.
  Suppose that Assumptions \ref{A1} (i) and (ii) are satisfied. Unless $f$ is constant a.e., the minimum of
  the GH function $F(x,t)$ will be always found at $t=0$, and the corresponding $x$ will be
  an optimal solution for \eqref{nonconvprob}.
\end{theorem}
%\cit{evans2010partial} introduces this theorem
%in a general way as {\it the maximum principle for the Cauchy problem}.
% {\it The maximum principle} for the {\it heat equation}, together with
%a well-known fact in statistical physics that solutions of the {\it heat equation} are
%given by Gaussian convolution with the original function $f$ leads to
%the statement in Theorem~\ref{thm: main}.
We present a proof of this theorem in Appendix \ref{subsec:proof_optimality}. 
The proof becomes much easier than
that in \cit{evans2010partial} due to its considering a specific case.

Let us next introduce an update rule for $t$ utilizing the derivative information.
When we solve the problem \eqref{xtvar_prob} using a gradient descent method,
the update rule for $t$ becomes $t_{k+1}=t_k-\eta\frac{\partial F}{\partial t}$,
where $\eta$ is a step size. The formula \eqref{heatequation} in the {\it heat equation}
implies that the derivative $\frac{\partial F}{\partial t}$ is equal to the Laplacian $\Delta_x F$,
i.e., $\frac{\partial F}{\partial t} = \mathrm{tr}(\mathrm{H}_F(x))$,
where $\mathrm{H}_F(x)$ is the Hessian of $F$ in terms of $x$.
Since $\mathrm{tr}(\mathrm{H}_F(x))$ represents the sharpness of minima \cit{Dinh2017sharp}, %; sharpness of a minimizer can be characterized by the magnitude of the eigenvalues of the Hessian matrix of $f$.
this update rule can sometimes decrease $t$ quickly around a minimum and
find a better solution. See Appendix \ref{subsec:toy} for an example of such a problem.
%\memo{Can you explain the same pace?? If difficult maybe removed? : If $\frac{\partial F}{\partial t}$ is used, the descent speed of $t$ will keep the same pace with the change of $x$.}
\subsection{SLGH algorithm}
\label{sec:GH}

Let us next introduce our proposed SLGH algorithm, %We choose gradient descent based algorithms to solve the above optimization problem \eqref{xtvar_prob}.
%Here, we propose two algorithms,
which has two variants with different update rules for $t$:
SLGH with a fixed-ratio update rule ($\text{SLGH}_\text{r}$) and
SLGH with a derivative update rule ($\text{SLGH}_\text{d}$).
$\text{SLGH}_\text{r}$ updates $t$ by multiplying a decreasing factor $\gamma$ (e.g., 0.999)
at each iteration. In contrast to this, $\text{SLGH}_\text{d}$ updates $t$
while using derivative information. Details are described in Algorithm~\ref{alg:GH}.
Algorithm~\ref{alg:GH} transforms a double loop Algorithm~\ref{algo:Continuation} into
a single loop algorithm. This single loop structure can significantly reduce
the number of iterations while ensuring the advantages of the GH method.

\begin{figure}[H]
\begin{algorithm}[H]
\caption{Deterministic/Stochastic Single Loop GH algorithm (SLGH)} \label{alg:GH}
\begin{algorithmic}
  \REQUIRE Iteration number $T$, initial solution $x_1$,
  initial smoothing parameter $t_1$, step size $\beta$ for $x$,
  step size $\eta$ for $t$, decreasing factor $\gamma \in (0,1)$,
  sufficient small positive value $\epsilon$\\
\FOR{$k=1$ to $T$}
    \STATE
            $
            x_{k+1} = x_k - \beta \widehat{G}_x,
            \ \widehat{G}_x = \left\{\begin{array}{l}
                \nabla_x F(x_k,t_k)\  (\text{determ.}) \\
                \nabla_x \bar{F}(x_k,t_k;\xi_k),\ \xi_k\sim P \  (\text{stoc.})
            \end{array}\right.
            $
            \begin{align*}
            t_{k+1} &= \left\{\begin{array}{l}
                \gamma t_{k}\  (\text{SLGH}_{\text{r}}) \\
                \text{max}\{\text{min} \{t_{k} - \eta \widehat{G}_t,
                \  \gamma t_{k}\}
                ,\ \epsilon'\}\ (\text{SLGH}_{\text{d}})
            \end{array},\right.
            \widehat{G}_t = \left\{\begin{array}{l}
                \frac{\partial F(x_{k},t_{k})}{\partial t}\   (\text{determ.}) \\
                \frac{\partial \bar{F}(x_{k},t_{k};\xi_{k})}{\partial t},\ \xi_k\sim P\ (\text{stoc.})
            \end{array}\right.
        \end{align*}
\ENDFOR
\end{algorithmic}
\end{algorithm}
\end{figure}

In the stochastic setting of \eqref{stocprob}, the gradient of $F(x,t)$ in terms of $x$ is approximated by $\nabla_x \bar{F}(x,t;\xi)$ with randomly chosen $\xi$, where $\bar{F}(x,t;\xi)$ is the GH function of $\bar{f}(x;\xi)$. Likewise, the derivative of $F(x,t)$ in terms of $t$ is approximated by $\frac{\partial \bar{F}(x,t;\xi)}{\partial t}$.  The stochastic algorithm in Algorithm~\ref{alg:GH} uses one sample $\xi_k$.
% Note that stochastic gradient descent is chosen %as the first-order method
% to update $x$ and $t$ for simplicity and theoretical analysis.
We can extend the stochastic approach to a minibatch one by approximating $\nabla_x F(x,t)$ by
$\frac{1}{M}\sum_{i=1}^M\nabla_x \bar{F}(x,t;\xi_i)$ with samples $\{\xi_1,\ldots,\xi_M\}$ of some batch size $M$, but for the sake of simplicity, we here assume one sample in each iteration. In this setting, the gradient complexity matches the iteration complexity; thus, we also use the term ``iteration complexity'' in the stochastic setting.
Other methods, such as momentum-accelerated method \cit{sutskever2013importance} and Adam \cit{kingma2014adam} can also be applied here. According to Theorem \ref{thm: main}, the final smoothing parameter needs to be zero. Thus, we multiply $\gamma$ by $t$ even in $\text{SLGH}_{\text{d}}$ when the decrease of $t$ is insufficient. We also assure that $t$ is larger than a sufficiently small positive value $\epsilon'>0$ during an update to prevent $t$ from becoming negative.
%\memo{You can move this to numerical experiment part. This is the thing you compare with Algorithm~\ref{alg:GH} to emphasize that the update rule of $t$ is good. : And we proposed a simple version of SLGH algorithm named SLGH-constant, which only use $\gamma$ without the derivative $\frac{\partial F}{\partial t}$. We will show that although SLGH is a bit complicated compared to SLGH-constant, it can achieve fewer iterations and better solutions.}

\subsection{Convergence analysis for SLGH} \label{sec:convSLGH}
Let us next analyze the worst-case iteration complexity for both deterministic and stochastic SLGHs, but, before that, let us first show some properties for
Gaussian smoothed function $F(x,t)$ under Assumption \ref{A1} for the original function $f(x)$.
In the complexity analyses in this paper,
we always assume that $\gamma$ is bounded from above by a universal constant $\bar{\gamma} < 1$,
which implies $1/(1 - \gamma) = O(1)$.
%is utilized because it often has better properties than $f(x)$. At least, the assumptions for original function also hold for these smoothed functions.
\begin{lemma} 
  Let $f(x)$ be a $L_0$-$Lipschitz$ function. Then, for any $t>0$, its Gaussian smoothed function $F(x,t)$ will then also be $L_0$-$Lipschitz$ in terms of $x$.
  %and $L_0(F)\leq L_0(f)$. 
  Let $f(x)$ be a $L_1$-$smooth$ function. Then, for any $t>0$, $F(x,t)$ will also be $L_1$-$smooth$ in terms of $x$.
  % and $L_1(F)\leq L_1(f)$.
  \label{lem:Lip}
\end{lemma}
Lemma~\ref{lem:Lip} indicates that  Assumption \ref{A1} given to the function $f(x)$ also guarantees the same properties for $F(x,t)$.
%, though the Lipschitz constants are different.
Below, we give some bounds between the smoothed function $F(x,t)$ and the original function $f(x)$.
\begin{lemma}
    Let $f$ be a $L_0$-$Lipschitz$ function. Then, for any $x\in\mathbb{R}^d$, $F(x, t)$ is also $L_0\sqrt{d}$-$Lipschitz$ in terms of $t$, i.e., for any $x$, smoothing parameter values $t_1, t_2>0$,  we have $|F(x,t_1) - F(x,t_2)|\leq L_0\sqrt{d}|t_1-t_2|.$
    \label{lem:Lip_t}
\end{lemma}
On the basis of Lemmas~\ref{lem:Lip} and \ref{lem:Lip_t}, the convergence results of our deterministic and stochastic SLGH algorithms can be given as in
Theorems~\ref{iter_determin} and \ref{thm:iter_stochastic}, respectively. Proofs of the following theorems are given in Appendix \ref{subsec:proof_first_order}.
% , and the definitions of $\hat{x}$ are provided in the proofs.
Let us first deal with the deterministic setting.
\begin{theorem}[\textbf{Convergence of SLGH, Deterministic setting}]
    Suppose Assumption \ref{A1} holds
    % The output sequence of the SLGH algorithm $\{x_k\}_{k=1}^T$ will then satisfy
    % $
    %     \min_{k \in [T]} \| \nabla f( x_k ) \|^2
    %     =
    %     O \left(
    %         \frac{1}{T} \left( 1 + d^{3/2} \sum_{k=1}^T t_k \right)
    %     \right)
    % $.
    % Consequently,
    , and let $\hat{x}:=x_{k'},\ k' = \mathop{\mathrm{argmin}}_{k\in[T]} \| \nabla f(x_k)\|$. Set the stepsize for $x$ as $\beta=1/L_1$. Then,
    % if we update $t_k$ as in Algorithm~\ref{alg:GH} with an arbitrary constant $\gamma \in (0, 1)$,
    for any setting of the parameter $\gamma$, $\hat{x}$ satisfies
    $\|\nabla f(\hat{x})\|\leq\epsilon$ with the iteration complexity of
    $T = O\left(d^{3/2}/\epsilon^2\right)$.
    Further,
    if we choose
    $\gamma \leq d^{-\Omega(\epsilon^2)}$,
    %$\gamma = d^{-3 \epsilon^2/2}$,
    the iteration complexity can be bounded as
    $T = O({1}/{\epsilon^2})$.
    \label{iter_determin}
\end{theorem}
% \begin{remark}
% \label{remark}
% This theorem indicates that there exists a trade-off between the convergence speed and the possibility of finding better solutions. Let us consider rewriting this complexity into a more concrete form. If $t_{k+1}\leq\gamma t_k\ (k\in[T])$ are satisfied, we get $\sum_{k=1}^Tt_k=\frac{t_1(1-\gamma^T)}{1-\gamma}$. Thus, we can see that the SLGH algorithm has the iteration complexity $O\left(\frac{d^{3/2}}{\epsilon^2}\right)$ for achieving an $\epsilon$-stationary point when we regard $t_1$ as a constant. Since the iteration complexity for gradient descent methods to reach an $\epsilon$-stationary point is $O(\frac{1}{\epsilon^2})$ \cit{nesterov2004intro}, this result seems to imply that the proposed algorithm takes $O(d^{3/2})$ times longer than gradient descent methods. However, we can eliminate the dependence in terms of $d$
% when $\gamma$ is chosen as $\gamma=O(d^{-3\epsilon^2/2})$, which implies $\gamma^{1/\epsilon^2}=O(d^{-3/2})$;
% %if $\gamma^{1/\epsilon^2}=O(d^{-3/2}) \Rightarrow\gamma=O(d^{-3\epsilon^2/2})$ holds.
% in such a case, the value of $t$ becomes $O(d^{-3/2})$ in $1/\epsilon^2$ iterations, and we can see that the algorithm finds an $\epsilon$-stationary point after another $O(1/\epsilon^2)$ iterations by substituting $\sum_{k=1}^Tt_k$ for $O(d^{-3/2})$.
% \end{remark}
% This theorem indicates that there exists a trade-off between the convergence speed and the possibility of finding better solutions.
%If we focus on exploring better solutions,
%i.e.,
This theorem indicates that
if we choose $\gamma$ close to $1$,
then the iteration complexity can be $O\left( d^{3/2}/\epsilon^2 \right)$,
which is $O(d^{3/2})$ times larger than the $O(1/\epsilon^2)$-iteration complexity by the standard gradient descent methods \cit{nesterov2004intro}.
% standard gradient descent methods.
% Let us consider rewriting this complexity into a more concrete form. If $t_{k+1}\leq\gamma t_k\ (k\in[T])$ are satisfied, we get $\sum_{k=1}^Tt_k=\frac{t_1(1-\gamma^T)}{1-\gamma}$. Thus, we can see that the SLGH algorithm has the iteration complexity $O\left(\frac{d^{3/2}}{\epsilon^2}\right)$ for achieving an $\epsilon$-stationary point when we regard $t_1$ as a constant.
% Since the iteration complexity for gradient descent methods to reach an $\epsilon$-stationary point is $O(\frac{1}{\epsilon^2})$ \cit{nesterov2004intro},
% the iteration complexity of $O(\frac{d^{3/2}}{\epsilon^2})$ in Theorem~\ref{iter_determin} implies that the proposed algorithm may take $O(d^{3/2})$ times longer than standard gradient descent methods.
However,
we can remove this dependency on $d$
to obtain an iteration complexity matching that of the standard gradient descent,
by choosing $\gamma \leq d^{-\Omega(\epsilon^2)}$,
as shown in Theorem~\ref{iter_determin}.
Empirically,
settings of $\gamma$ close to $1$,
e.g., $\gamma = 0.999$,
seem to work well enough,
as demonstrated in Section~\ref{sec:exp}.

An inner loop of the double loop GH method using the standard GD requires the same complexity as the standard GD method up to constants since the objective smoothed function of inner optimization problem is $L_1$-smooth function. By considering the above results, we can see that the SLGH algorithm becomes faster than the double loop one by around the number of outer loops.

To provide theoretical analyses in the stochastic setting, we need additional standard assumptions.

\begin{assumption}{A2}\label{A2}$\ $
\renewcommand{\labelenumi}{(\roman{enumi})}
\begin{enumerate}
    \item The stochastic function $\bar{f}(x;\xi)$ becomes an unbiased estimator of $f(x)$. That is, for any $x\in\mathbb{R}^d$, $f(x) = \mathbb{E}_\xi [\bar{f}(x;\xi)]$ holds.
    \item For any $x\in\mathbb{R}^d$, the variance of the stochastic gradient oracle is bounded as $\mathbb{E}_{\xi}[\|\nabla_x\bar{f}(x;\xi)-\nabla f(x)\|^2] \leq \sigma^2$. Here, the expectation is taken w.r.t.~random vectors $\{\xi_k\}$.
\end{enumerate}

\end{assumption}
The following theorem shows the convergence rate in the stochastic setting.

\begin{theorem}[\textbf{Convergence of SLGH, Stochastic setting}]
    Suppose Assumptions \ref{A1} and \ref{A2} hold. 
    % The output sequence of the SLGH algorithm $\{x_k\}_{k=1}^T$ will then satisfy
    % $
    %     \frac{1}{T}\sum_{t = 1}^T \mathbb{E} [ \| \nabla f( x_k ) \|^2]
    %     = O \left(\frac{1+\sqrt{d}\sum_{k=1}^T\mathbb{E}_{\xi}[|t_{k+1}-t_k|]}{\sqrt{T}}+\frac{d^{3/2}}{T}\sum_{k=1}^{T}\mathbb{E}_{\xi}[t_k]\right)
    % $,
    % where the expectation is taken w.r.t.~random vectors $\{\xi_k\}$.
    % Consequently,
    % if we update $t_k$ as in Algorithm~\ref{alg:GH} with an arbitrary constant $\gamma \in (0, 1)$,
    Take $k_1 := \Theta(1/\epsilon^4)$ and $k_2 := O\left(\log_{\gamma} \min\{ d^{-1/2},  d^{-3/2} \epsilon^{-2} \}\right)$ and define $k_0 = \min \{ k_1, k_2 \}$. Let $\hat{x}:=x_{k'}$, where $k'$ is chosen from a uniform distribution over $\{ k_0+1, k_0+2, \ldots, T \}$. 
    Set the stepsize for $x$ as $\beta=\min\left\{1/L_1, 1/\sqrt{T-k_0}\right\}$.
    Then,
    for any setting of the parameter $\gamma$,
    $\hat{x}$ satisfies
    $\mathbb{E}[ \|\nabla f(\hat{x})\| ] \leq\epsilon$ with the iteration complexity of
    $T = O \left(
        d/\epsilon^4 + d^{3/2}/\epsilon^2
    \right)$
    where the expectation is taken w.r.t.~random vectors $\{\xi_k\}$.
    Further,
    if we choose $\gamma \leq (\max\{ d^{1/2},  d^{3/2} \epsilon^2 \})^{- \Omega(\epsilon^4)}$,
    the iteration complexity can be bounded as
    $T = O({1}/{\epsilon^4})$.
\label{thm:iter_stochastic}
\end{theorem}
We note that the iteration complexity of $T = O({1}/{\epsilon^4})$ for sufficiently small $\gamma$
matches that for the standard stochastic gradient descent (SGD) shown,
e.g.,
by \cit{ghadimi2013stochastic}.
% Again, this result yields the fact that if $t_{k+1}\leq\gamma t_k,\ k\in[T]$ are satisfied, the SLGH algorithm has the iteration complexity $O\left(\frac{d^2}{\epsilon^4}\right)$ in reaching an $\epsilon$-stationary point. Further, in a way similar to that in the deterministic setting, if $\gamma=O(d^{-3\epsilon^2/2})$ is satisfied, the dependence in terms of $d$ vanishes; thus, the iteration complexity matches that of stochastic GD (SGD) $O\left(\frac{1}{\epsilon^4}\right)$ \cit{ghadimi2013stochastic}.

% \memo{In Algorithm~\ref{alg:GH}, $\beta$ or $\{\beta_k\}$? Theorems assume $\beta$}
%================================================================================================%

\section{Zeroth-order single loop Gaussian homotopy algorithm}
\label{sec:zeroth-order}
In this section, we introduce a zeroth-order version of the SLGH algorithms. This ZOSLGH algorithm is proposed for those optimization problems in which Gaussian smoothing convolution is difficult to compute, or in which only function values can be queried.

\subsection{ZOSLGH algorithm}
For cases in which only function values are accessible, approximations for the gradient in terms of $x$ and derivative in terms of $t$ are needed.
\cit{nesterov2017random} has shown that the gradient of the smoothed function $F(x,t)$ can be represented as
\begin{align}
    \nabla_x F(x,t) &= \frac{1}{t}\mathbb{E}_u([f(x+tu)-f(x)]u),\ u \sim \mathcal{N}(0,\mathrm{I}_d).
\end{align}
Thus, the gradient $\nabla_x F(x,t)$ can be approximated by an unbiased estimator $\Tilde{g}_{x}(x,t;u)$ as %\memo{Why did you use the notation $\Tilde{g}_{x,u}(x,t)$ instead of $\Tilde{g}_{u}(x,t)$?}
\begin{align}
     \Tilde{g}_{x}(x,t;u) := \frac{1}{t}(f(x+tu)-f(x))u,\ u \sim \mathcal{N}(0,\mathrm{I}_d).
     \label{zogradient}
\end{align}

The derivative $\frac{\partial F}{\partial t}$ is equal to the trace of the Hessian of $F(x,t)$ because the Gaussian smoothed function is the solution of the {\it heat equation} $\frac{\partial F}{\partial t} = \mathrm{tr}(\mathrm{H}_F(x))$. We can estimate $\mathrm{tr}(\mathrm{H}_F(x))$ on the basis of the second order Stein's identity \cit{stein1972bound} as follows:
\begin{align}
    \mathrm{H}_F(x) \approx \frac{(vv^\top-\mathrm{I}_d)}{t^2}(f(x+tv)-f(x)),\ v \sim \mathcal{N}(0,\mathrm{I}_d).
\end{align}
Thus, the estimator for derivative can be written as:
\begin{align}
     \Tilde{g}_{t}(x,t;v) := \frac{(v^\top v-d)(f(x+tv)-f(x))}{{t}^2},\ v \sim \mathcal{N}(0,\mathrm{I}_d).
     \label{zoderivative}
\end{align}
%The following is our propose ZOGH algorithm, a zeroth-order algorithm using zeroth-order estimators for gradient and Hessian.%
As for the stochastic setting, $f(x)$ in \eqref{zogradient} and \eqref{zoderivative} is replaced by the stochastic function $\bar{f}(x;\xi)$ with some randomly chosen sample $\xi$. The gradient
$\nabla_x \bar{F}(x,t;\xi)$ of its GH function $\bar{F}(x,t;\xi)$ can then be approximated by %$\Tilde{G}_{x,u}(x,t;\xi) := \frac{F(x+tu;\xi)-F(x;\xi)}{t}u$
$\Tilde{G}_{x}(x,t;\xi,u) := \frac{\bar f(x+tu;\xi)-\bar f(x;\xi)}{t}u$,
and the derivative $\frac{\partial \bar{F}}{\partial t}$ can be approximated by $\Tilde{G}_{t}(x,t;\xi,v) := \frac{(v^\top v-d)(\bar f(x+tv;\xi)-\bar f(x;\xi))}{{t}^2}$ %\memo{$\bar{F}(x,t;\xi)$ is the GH function  for $\bar{f}(x;\xi)$? Or better notation? Anyway, notation need to be unified}
(see Algorithm \ref{alg:ZOGH} for more details).
\begin{figure}[H]
\begin{algorithm}[H]
\caption{Deterministic/Stochastic Zeroth-Order Single Loop GH algorithm (ZOSLGH)}
\label{alg:ZOGH}
\begin{algorithmic}
\REQUIRE Iteration number $T$, initial solution $x_1$, initial smoothing parameter $t_1$, step size $\beta$ for $x$, step size $\eta$ for $t$, decreasing factor $\gamma \in (0,1)$, sufficient small positive value $\epsilon$\\
\FOR{$k=1$ to $T$}
    \STATE Sample $u_k$ from $\mathcal{N}(0,\mathrm{I}_d)$
    \STATE
        $
            x_{k+1} = x_k - \beta \bar{G}_{x,u},
            \ \bar{G}_{x,u} = \left\{ \begin{array}{l}
                \Tilde{g}_{x}(x_k,t_k;u_k)\  (\text{determ.}) \\
                \Tilde{G}_{x}(x_k,t_k;\xi_k,u_k),\ \xi_k\sim P\  (\text{stoc.})
            \end{array}\right.
        $
    \STATE Sample $v_k$ from $\mathcal{N}(0,\mathrm{I}_d)$
    \begin{align*}
            t_{k+1} = \left\{\begin{array}{l}
                \gamma t_{k}\  (\text{SLGH}_{\text{r}}) \\
                \text{max}\{\text{min} \{t_{k} - \eta \bar{G}_{t,v},
                \gamma t_{k}\}
                ,\epsilon'\}\ (\text{SLGH}_{\text{d}})
            \end{array},\right.
            \bar{G}_{t,v} = \left\{\begin{array}{l}
                \Tilde{g}_{t}(x_k,t_k;v_k)\ (\text{determ.}) \\
                \Tilde{G}_{t}(x_k,t_k;\xi_k,v_k),\ \xi_k\sim P\  (\text{stoc.})
            \end{array}\right.
    \end{align*}
\ENDFOR
%\ENSURE ~~ $x_T$
\end{algorithmic}
\end{algorithm}
\end{figure}

\subsection{Convergence analysis for ZOSLGH}
We can analyze the convergence results using concepts similar to those used with the first-order SLGH algorithm.
Below are the convergence results for ZOSLGH in both the deterministic and stochastic settings.
Proofs of the following theorems are given in Appendix \ref{subsec:proof_zeroth_order}, and the definitions of $\hat{x}$ are provided in the proofs. 
We start from the deterministic setting, which is aimed at the deterministic problem \eqref{nonconvprob}.
\begin{theorem}[\textbf{Convergence of ZOSLGH, Deterministic setting}]\label{thm:zo_determinstic}
    Suppose Assumption \ref{A1} holds. 
    % The output sequence of the ZOSLGH algorithm $\{x_k\}_{k=1}^T$ will then satisfy $
    %     \frac{1}{T}\sum_{t = 1}^T \mathbb{E} [ \| \nabla f( x_k ) \|^2]
    %     = O\left(\frac{d}{T}\left(1+d\sum_{k=1}^T \mathbb{E}_u[t_k^2]+\sqrt{d}\sum_{k=1}^T \mathbb{E}_u[t_k]\right)\right)
    % $,
    % where the expectation is taken w.r.t.~random vectors $\{u_k\}$.
    % Consequently,
    % if we update $t_k$ as in Algorithm~\ref{alg:ZOGH} with an arbitrary constant $\gamma \in (0, 1)$,
Take $k_1 := \Theta(d/\epsilon^2)$ and $k_2 := O\left(\log_{\gamma} d^{-1/2} \right)$, and define $k_0 = \min \{ k_1, k_2 \}$. Let $\hat{x}:=x_{k'}$, where $k'$ is chosen from a uniform distribution over $\{ k_0+1, k_0+2, \ldots, T \}$. Set the stepsize for $x$ as $\beta=1/(2(d+4)L_1)$.
    Then,
    for any setting of the parameter $\gamma$,
    $\hat{x}$ satisfies
    $\mathbb{E}[ \|\nabla f(\hat{x})\| ] \leq\epsilon$ with the iteration complexity of
    $T = O(d^2 / \epsilon^2 )$,
    where the expectation is taken w.r.t.~random vectors $\{u_k\}$ and $\{v_k\}$.
    Further,
    if we choose $\gamma \leq d^{-\Omega(\epsilon^2 / d)}$,
    the iteration complexity can be bounded as
    $T = O(d / {\epsilon^2})$.
\end{theorem}
This complexity of $O({d}/{\epsilon^2})$ for $\gamma \leq d^{-\Omega(\epsilon^2 / d)}$ matches that of zeroth-order GD (ZOGD) \cit{nesterov2017random}.
%  if $\gamma=O(d^{-\epsilon^2/2})$ is satisfied.

Let us next introduce the convergence result for the stochastic setting. As shown in \cit{ghadimi2013stochastic}, if we take the expectation for our stochastic zeroth-order gradient oracle with respect to both $\xi$ and $u$, under Assumption \ref{A2} (i), we will have
\begin{align*}
    \mathbb{E}_{\xi,u} [\Tilde{G}_{x}(x,t;\xi,u)] =  \mathbb{E}_u [\mathbb{E}_\xi [\Tilde{G}_{x}(x,t;\xi,u)|u]] = \nabla_x F(x,t).
\end{align*}
Therefore, $\zeta_k := (\xi_k,u_k)$ behaves similarly to $u_k$ in the deterministic setting.
\begin{theorem}[\textbf{Convergence of ZOSLGH, Stochastic setting}]
    \label{thm:zo_stochastic}
    Suppose Assumptions \ref{A1} and \ref{A2} hold.
    Take $k_1 := \Theta(d/\epsilon^4)$ and $k_2 := O\left(\log_{\gamma} d^{-1/2} \right)$, and define $k_0 = \min \{ k_1, k_2 \}$. Let $\hat{x}:=x_{k'}$, where $k'$ is chosen from a uniform distribution over $\{ k_0+1, k_0+2, \ldots, T \}$. Set the stepsize for $x$ as $\beta=\mathop{\rm min}\{\frac{1}{2(d+4)L_1}, \frac{1}{\sqrt{(T - k_0)(d+4)}}\}$. 
    % The output sequence of the ZOSLGH algorithm $\{x_k\}_{k=1}^T$ will then satisfy 
    % \begin{align*}
    %     &\frac{1}{T}\sum_{t = 1}^T \mathbb{E} [ \| \nabla f( x_k ) \|^2]
    %     \\
    %     &= O\left(\frac{\sqrt{d}\left(1+\sqrt{d}\sum_{k=1}^T\mathbb{E}_\zeta[|t_{k+1}-t_k|]\right)}{\sqrt{T}} \right.
    %     \\
    %     &\quad \quad + \left. \frac{d\left(d\sum_{k=1}^T \mathbb{E}_\zeta[t_k^2]+\sqrt{d}\sum_{k=1}^T\mathbb{E}_\zeta[t_k]+1\right)}{T}\right),
    % \end{align*}
    % where the expectation is taken w.r.t.~random vectors $\{u_k\}$ and $\{\xi_k\}$.
    % Consequently,
    % if we update $t_k$ as in Algorithm~\ref{alg:ZOGH} with an arbitrary constant $\gamma \in (0, 1)$,
    Then,
    for any setting of the parameter $\gamma$,
    $\hat{x}$ satisfies
    $\mathbb{E}[ \|\nabla f(\hat{x})\| ] \leq\epsilon$ with the iteration complexity of
    $T = O({d^{2}}/{\epsilon^4})$,
    where the expectation is taken w.r.t.~random vectors $\{u_k\}$, $\{v_k\}$, and $\{\xi_k\}$.
    Further,
    if we choose $\gamma \leq d^{- \Omega(\epsilon^4/d)}$,
    the iteration complexity can be bounded as
    $T = O({d}/{\epsilon^4})$.
\end{theorem}

This complexity of $O(d/\epsilon^4)$ for $\gamma \leq d^{- \Omega(\epsilon^4/d)}$ also matches that of ZOSGD \cit{ghadimi2013stochastic}.

\section{Experiments}
\label{sec:exp}

In this section, we present our experimental results. We conducted two experiments. The first was to compare the performance of several algorithms including the proposed ones, using test functions for optimization. We were able to confirm the effectiveness and versatility of our SLGH methods for highly non-convex functions. We also created a toy problem in which $\text{ZOSLGH}_{\text{d}}$, which utilizes the derivative information $\frac{\partial F}{\partial t}$ for the update of $t$, can decrease $t$ quickly around a minimum and find a better solution than that with $\text{ZOSLGH}_{\text{r}}$. The second experiment was to generate examples for a black-box adversarial attack with different zeroth-order algorithms. The target models were well-trained DNNS for CIFAR-10 and MNIST, respectively. All experiments were conducted using Python and Tensorflow on Intel Xeon CPU and NVIDIA Tesla P100 GPU. We show the results of only the adversarial attacks due to the space limitations; other results are given in Appendix \ref{sec:test_funcs}.

\textbf{Generation of per-image black-box adversarial attack example.} Let us consider the unconstrained black-box attack optimization problem in \cit{chen2019zo}, which is given by
\begin{align}
    {\mathop{\rm minimize}\limits_{x\in \mathbb{R}^d}}\  f(x) :=& \lambda \ell(0.5\text{tanh}(\text{tanh}^{-1}(2a)+x)) + \|0.5\text{tanh}(\text{tanh}^{-1}(2a)+x)-a\|^2,\nonumber
\end{align}
where %$x\in \mathbb{R}^d$ is the optimization variable,
$\lambda$ is a regularization parameter, $a$ is the input image data, and  $tanh$ is the element-wise operator which helps eliminate the constraint representing the range of adversarial examples. The first term $\ell(\cdot)$ of $f(x)$ is the loss function for the untargeted attack in \cit{carlini2017towards}, and the second term $L_2$ distortion is the adversarial perturbation (the lower the better). The goal of this problem is to find the perturbation that makes the loss $\ell(\cdot)$ reach its minimum while keeping $L_2$ distortion as small as possible. The initial adversarial perturbation $x_0$ was set to $0$. We say a successful attack example has been generated when the loss $\ell(\cdot)$ is lower than the attack confidence (e.g., $1e-10$).

Let us here compare our algorithms, $\text{ZOSLGH}_{\text{r}}$ and $\text{ZOSLGH}_{\text{d}}$, to three zeroth-order algorithms: ZOSGD \cit{ghadimi2013stochastic}, ZOAdaMM \cit{chen2019zo}, and
ZOGradOpt \cit{hazan2016graduated}. ZOGradOpt is a homotopy method with a double loop structure. In contrast to this, ZOSGD and ZOAdaMM are SGD-based zeroth-order methods and thus do not change the smoothing parameter during iterations.
%\memo{Maybe this can be moved to supplimentary. It is required to write but sounds too detailed}

%The most important parameter here is the smoothing parameter $t$. For better comparison, we chose $10$ for GH methods and $0.005$ for the  others. Other common parameters, including step size and attack confidence, were selected with the same value for all algorithms. The initial adversarial perturbation $x_0$ was set to $0$.
%For CIFAR-10 task, the smoothing parameter settings are: $0.005$ and $10$ for ZOSGD, $0.005$ for ZOAdaMM, $10$ for ZOGradOpt, ZOSLGH and ZOSLGH-constant. For MNIST task, the smoothing parameter settings are: $0.005$ and $100$ for ZOSGD, $0.005$ for ZOAdaMM, $100$ for ZOGradOpt, ZOSLGH and ZOSLGH-constant. Note that other parameters including learning rate, attack confidence are chosen with the same value for all algorithms. See appendix B for more details.

\begin{figure}[H]
\centering
\begin{tabular}{ccc}
    \subfigure[CIFAR-10, Image ID = 56]{
    \label{fig2a}
    \begin{minipage}[t]{0.32\linewidth}
        \centering
        \includegraphics[width=2in]{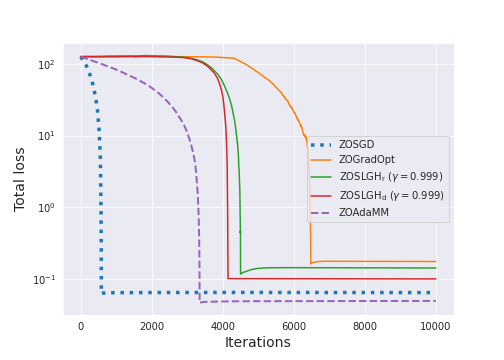}
    \end{minipage}%
    }
    \subfigure[CIFAR-10, Image ID = 34]{
    \label{fig2b}
    \begin{minipage}[t]{0.32\linewidth}
        \centering
        \includegraphics[width=2in]{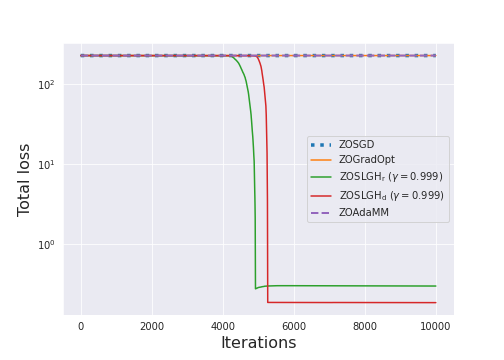}
     \end{minipage}
    }
    \subfigure[CIFAR-10, Image ID = 102]{
    \label{fig2c}
    \begin{minipage}[t]{0.32\linewidth}
        \centering
        \includegraphics[width=2in]{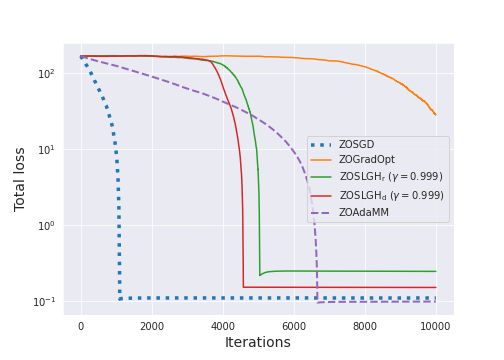}
     \end{minipage}
    }
    \\\vspace{-1mm}
    \subfigure[MNIST, Image ID = 32]{
    \label{fig2d}
    \begin{minipage}[t]{0.32\linewidth}
        \centering
        \includegraphics[width=2in]{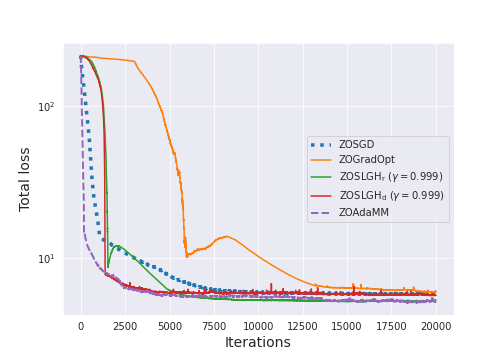}
    \end{minipage}%
    }
    \subfigure[MNIST, Image ID = 45]{
    \label{fig2e}
    \begin{minipage}[t]{0.32\linewidth}
        \centering
        \includegraphics[width=2in]{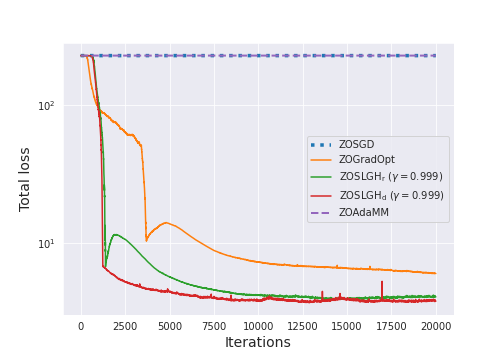}
     \end{minipage}
    }
    \subfigure[MNIST, Image ID = 95]{
    \label{fig2f}
    \begin{minipage}[t]{0.32\linewidth}
        \centering
        \includegraphics[width=2in]{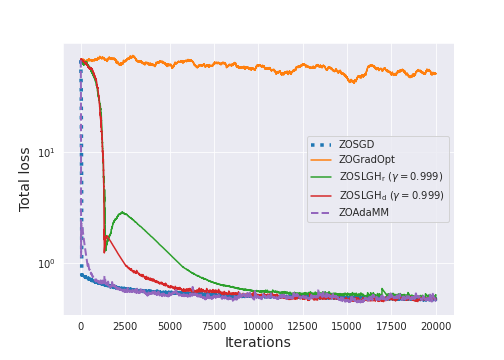}
     \end{minipage}
    }
\end{tabular}
\centering
\caption{\begin{tabular}{l}Total loss for generating per-image black-box adversarial examples for different images \\ of CIFAR-10 and MNIST (log scale).
\end{tabular}}
\label{fig:plots}
\end{figure}

Table \ref{table:results} and Figure \ref{fig:plots} show results for our experiment. We can see that SGD-based algorithms are able to succeed in the first attack with far fewer iterations than our GH algorithms (e.g., Figure \ref{fig2a}, Figure \ref{fig2d}). Accordingly, the value of $L_2$ distortion decreases slightly more than GH methods. However, SGD-based algorithms have lower success rates than do our SLGH algorithms. This is because SGD-based algorithms remain around a local minimum $x=0$ when it is difficult to attack, while GH methods can escape the local minima due to sufficient smoothing (e.g., Figure \ref{fig2b}, Figure \ref{fig2e}). Thus, the SLGH algorithms are, on average, able to decrease total loss over that with SGD-based algorithms. In a comparison within GH methods, ZOGradOpt requires more than 6500 iterations to succeed in the first attack due to its double loop structure (e.g., Figure \ref{fig2c}, Figure \ref{fig2f}). In contrast to this, our SLGH algorithms achieve a high success rate with far fewer iterations. Please note that $\text{SLGH}_{\text{d}}$ takes approximately twice the  computational time per iteration than the other algorithms because it needs additional queries for the computation of the derivative in terms of $t$. See Appendix \ref{sec:black_box} for a more detailed presentation of the experimental setup and results.

\begin{table}[H]
\centering
  \caption[]{Performance of a per-image attack over $100$ images of CIFAR-10 under $T = 10000$ iterations. ``Succ. rate'' indicates the ratio of success attack, ``Avg. iters to 1st succ.''  is the average number of iterations to reach the first successful attack , ``Avg. $L_2$ (succ.)'' is the average of $L_2$ distortion taken among successful attacks, and ``Avg. total loss'' is the average of total loss $f(x)$ over 100 samples. Please note that the standard deviations are large since the attack difficulty varies considerably from sample to sample.}

\begin{tabular}{cc|c|c|c|c}
\toprule
&Methods & \begin{tabular}{c}Succ. rate\end{tabular} & \begin{tabular}{c} Avg. iters \\ to 1st succ.\end{tabular}  & \begin{tabular}{c} Avg. $L_2$\ \\ (succ.)\end{tabular} & \begin{tabular}{c} Avg. total loss\end{tabular}\\ \hline
SGD algo. &ZOSGD & $88\%$ & $\textbf{835} \pm 1238$ & $0.076 \pm 0.085$ & $27.70 \pm 74.80$\\
&ZOAdaMM& $85\%$ & $3335 \pm 2634$ & $\textbf{0.050} \pm 0.055$ & $20.24 \pm 62.48$  \\\hline
GH algo. &ZOGradOpt & $65\%$ & $6789 \pm 1901$ & $0.249 \pm 0.159$ & $41.45 \pm 76.04$ \\
&$\text{ZOSLGH}_{\text{r}}\ (\gamma=0.999)$ & $\textbf{93\%}$ & $4979 \pm 756$ & $0.246 \pm 0.178$ & $\textbf{14.26} \pm 54.61$ \\
&$\text{ZOSLGH}_{\text{d}}\ (\gamma=0.999)$ & $\textbf{92\%}$ & $4436 \pm 805$ & $0.150 \pm 0.084$ & $\textbf{16.49} \pm 58.69$\\
\bottomrule
\end{tabular}%
\label{table:results}%
\end{table}%

% \footnotetext{The standard deviations become large since the attack difficulty varies considerably from sample to sample.}
\vspace{-3mm}

%=================================================================================================%

\section{Summary and future work}
\label{sec:conclu}
We have presented here the deterministic/stochastic SLGH and ZOSLGH algorithms as well as their convergence results. They have been designed for the purpose of finding better solutions with fewer iterations by simplifying the homotopy process into a single loop. We consider this work to be a first attempt to improve the standard GH method.

Although this study has considered the case in which the accessible function contains some error and is possibly non-smooth, we assume the underlying objective function to be smooth. Further work should be carried out to investigate the case in which the objective function itself is non-smooth.

\paragraph{Acknowledgements} This work was supported by JSPS KAKENHI Grant Number 19H04069, JST ACT-I Grant Number JPMJPR18U5, and JST ERATO Grant Number JPMJER1903.

%Our work is restricted to a differentiable objective function because of Assumption A1, though the smoothed function becomes differentiable and our algorithms will work for a non-smooth function. In addition, the fixed step size is chosen for all proposed algorithms for simplicity. We want to replace Assumption A1 with a weaker one and accelerate our algorithms by using adaptive step size as our future research topics.

\newpage
\bibliographystyle{icml2022}
\bibliography{ref}

%%%%%%%%%%%%%%%%%%%%%%%%%%%%%%%%%%%%%%%%%%%%%%%%%%%%%%%%%%%%%%%%%%%%%%%%%%%%%%%
%%%%%%%%%%%%%%%%%%%%%%%%%%%%%%%%%%%%%%%%%%%%%%%%%%%%%%%%%%%%%%%%%%%%%%%%%%%%%%%
% APPENDIX
%%%%%%%%%%%%%%%%%%%%%%%%%%%%%%%%%%%%%%%%%%%%%%%%%%%%%%%%%%%%%%%%%%%%%%%%%%%%%%%
%%%%%%%%%%%%%%%%%%%%%%%%%%%%%%%%%%%%%%%%%%%%%%%%%%%%%%%%%%%%%%%%%%%%%%%%%%%%%%%
\newpage
\appendix
\onecolumn

\section{Related work}
\label{sec:related_work}
%Homotopy method was firstly proposed in 1987 by \cit{blake1987visual} and one recent application is tuning hyperparameters of  kernel ridge regression \cit{shao2019graduated}. During this development,
\paragraph{Iteration complexity analysis for GH methods} To the best of our knowledge, there are no existing works that give theoretical guarantee for the convergence rate except for \cit{hazan2016graduated}.\footnote{Their method is not exactly a GH method because it smooths the objective function using random variables sampled from the unit ball (or the unit sphere in a zeroth-order setting) rather than Gaussian random variables. However, for the sake of simplicity, we treat it as a GH method in this paper.} It characterized a parameterized family of non-convex functions referred to as ``$\sigma$-nice'', for which a GH algorithm converges to a global optimum. Moreover, it derived the convergence rate to an $\epsilon$-optimal solution for the $\sigma$-nice function. The framework of $\sigma$-nice imposes the two conditions: (i) the solution obtained in each inner loop is located sufficiently close to an optimal solution of the optimization problem in the next inner loop; (ii) the optimization problem in each inner loop is strongly convex around its optimal solutions. Unfortunately, it is not obvious whether we can efficiently judge a function is ``$\sigma$-nice'', and we cannot apply the analysis results to general non-convex functions. On the other hand, this work tackles a problem of different nature from \cit{hazan2016graduated} since it analyzes the convergence rate to an $\epsilon$-stationary point for general non-convex functions.

\paragraph{Guarantee for the value of the objective function}
\cit{mobahi2015theoretical} provided an upper bound on the objective value attained by a homotopy method.
The bound was characterized by a quantity that they referred to as ``optimization complexity'', which can be analytically computed when the objective function is expressed in some suitable basis functions such as Gaussian RBFs.

\paragraph{Other smoothing methods} Smoothing methods other than Gaussian smoothing include \cit{chen1993non, chen2012smoothing}. The smoothing kernel in those works is simpler but restricted to specific problem settings. For example, \cit{chen2012smoothing} constructs smoothing approximations for optimization problems that can be reformulated by using the plus function $(t)_{+}:=\max\{0, t\}$. 
% The use of Gaussian smoothing does not limit application examples over those of other smoothing methods.

\paragraph{Zeroth-order techniques} In problem settings in which the explicit gradient of the objective function cannot be calculated but the exact function values can be queried,
zeroth-order optimization has become increasingly popular due to its potential for wide application. Such a class of applications appears in black-box adversarial attacks on deep neural networks \cit{chen2019zo}, structured prediction \cit{sokolov2016stochastic}, and reinforcement learning \cit{xu2020zeroth}. Various zeroth-order methods (ZOSGD \cit{ghadimi2013stochastic}, ZOAdaMM \cit{chen2019zo}, ZOSVRG \cit{liu2018zeroth}) have been proposed for such black-box situations. All of them have been developed from ZOGD in \cit{nesterov2017random}, which introduces random gradient-free oracles based on Gaussian smoothing with fixed $t$. This trend also applies to research on the GH method.
\cit{hazan2016graduated} developed a GH method in the zeroth-order setting for which the objective is only accessible through a noisy value oracle. \cit{shao2019graduated} proposed a GH method for hyperparameter tuning  based on %Hazan's work
\cit{hazan2016graduated} using two-point zeroth-order estimators \cit{nesterov2017random}.

\section{Proofs for theorems and lemmas in Sections \ref{sec:first-order} and \ref{sec:zeroth-order}}

\paragraph{Notation:} We sometimes denote the expectation with respect to random variables $\chi_{S+1}, \ldots, \chi_T\ (S, T\in\mathbb{N}, T>S)$ as $\mathbb{E}_\chi[\cdot]$ for the sake of simplicity.

\subsection{Theorem \ref{thm: main}}
\label{subsec:proof_optimality}
\textbf{Proof for Theorem \ref{thm: main}:}
Since the optimization problem \eqref{nonconvprob} has an optimal value $f^\ast$ by Assumption \ref{A1} (ii), for any $t\in\mathcal{T}$ and for any $x\in\mathbb{R}^d$, we have
\begin{align*}
    F(x,t)-f^\ast = \mathbb{E}_u[f(x+tu)-f^\ast]\geq 0.
\end{align*}

Together with the relationship $F(x,0)=f(x)$, for any $x\in\mathbb{R}^d$, for any $t\in\mathcal{T}$ and for any optimal solution $x^\ast\in\mathbb{R}^d$ of the optimization problem \eqref{nonconvprob}, we have $F(x,t)-F(x^\ast, 0)\geq 0$. Furthermore, if we exclude cases where $f(x)$ is constant (a.e.), for any $(x,t)\in\mathbb{R}^d\times \mathcal{T}\setminus\{(x,0)\mid f(x)=f(x^\ast)\}$, we obtain
%the inequality in \eqref{31_ineq} becomes strict:
\begin{align*}
    F(x,t)-f^\ast = \mathbb{E}_u[f(x+tu)-f^\ast]> 0.
\end{align*}
Therefore, a minimum of the optimization problem of the GH function ${\mathop{\rm minimize}\limits_{x\in \mathbb{R}^d, t \in \mathcal{T}}}\ F(x,t)$ holds only at $t=0$ and the corresponding $x$ becomes an optimal solution of the original optimization problem ${\mathop{\rm minimize}\limits_{x\in {\mathbb{R}}^d}}\ f(x)$.
%we can see that the optimal solution of $\mathop{\rm min}\limits_{x\in \mathbb{R}^d} F(x,0)$ becomes also optimal for our optimization problem  ${\mathop{\rm minimize}\limits_{x\in {\mathbb{R}}^d}}\ f(x)$.
\hfill $\Box$

\subsection{First-order SLGH algorithm}
\label{subsec:proof_first_order}
At the beginning of the subsection, we introduce a lemma that gives upper bounds for moments of Gaussian random variables, and then prove the two lemmas which appeared in the main paper.

\begin{lemma}[\textbf{Lemma 1 in \cit{nesterov2017random}}]\label{lem:gauss_norm}
Let $u\in\mathbb{R}^d$ be a standard normal random variable. For $p\in[0,2]$, we have $\mathbb{E}_u[\|u\|^p]\leq d^{p/2}$. If $p\geq 2$, $\mathbb{E}_u[\|u\|^p]\leq (d+p)^{p/2}$ holds.
\end{lemma}

\textbf{Proof for Lemma \ref{lem:Lip}:}
According to the definition of Gaussian smoothing in the main paper, we have
\begin{align}
|F(x,t)-F(y,t)| &=  \left|\int (f(x+tz)k(z)-f(y+tz)k(z)) dz\right| \nonumber\\
&\leq \int \left|f(x+tz)-f(y+tz)\right|k(z) dz \nonumber\\
&\leq \int L_0\|x-y\|k(z)dz\nonumber\\
&\leq L_0\|x-y\|\nonumber.
\end{align}
The proof of $L_1$-$smooth$ is similar to that of $L_0$-$Lipschitz$:
\begin{align}
|\nabla_xF(x,t)-\nabla_xF(y,t)| &\leq \int |\nabla f(x+tz)-\nabla f(y+tz)|k(z) dz\nonumber\\
&\leq \int L_1\|x-y\|k(z)dz\nonumber\\
&\leq L_1\|x-y\|\nonumber.
\end{align}
\hfill $\Box$

The lemma has proved that the Lipschitz constants of $F(x,t)$ and $\nabla_x F(x,t)$ in terms of $x$ are smaller than
those of $f(x)$ and $\nabla f(x)$, respectively. Therefore we can use the Lipschitz constants $L_0$ and $L_1$ of $f(x)$ and $\nabla f(x)$
for $F(x,t)$ and $\nabla_x F(x,t)$. 

\textbf{Proof for Lemma \ref{lem:Lip_t}:}
\begin{align}
|F(x,t_1)-F(x,t_2)| &=  |\mathbb{E}_{u}[f(x+t_1u)-f(x+t_2u)]| \nonumber\\
&\leq \mathbb{E}_u[|f(x+t_1u)-f(x+t_2u)|]\nonumber\\
&\leq \mathbb{E}_u[L_0|t_1-t_2| \|u\|]\nonumber\\
&\leq L_0|t_1-t_2|\sqrt{d}\nonumber,
\end{align}
where the last inequality holds due to Lemma \ref{lem:gauss_norm}.
\hfill $\Box$

Before going to the convergence theorems, we introduce an additional useful lemma to estimate the gap between the gradient of the smoothed function and the true gradient.\\
\begin{lemma}\label{lem:grad_square_diff}
    Let $f$ be a $L_1$-$smooth$ function.\\
\textbf{(i) (\textbf{Lemma 4 in \cit{nesterov2017random}})}
    For any $x\in \mathbb{R}^d$ and $t>0$, we have
    \begin{align*}
        \|\nabla f(x)\|^2 \leq 2\|\nabla_x F(x,t)\|^2 + \frac{t^2}{2}L_1^2(d+6)^3.
    \end{align*}
\textbf{(ii)} Further, if $f$ is $L_0$-Lipschitz, for any $x\in \mathbb{R}^d$ and $t>0$, we have
\begin{align*}
        \|\nabla f(x)\|^2 \leq \|\nabla_x F(x,t)\|^2 + tL_0L_1(d+3)^{3/2}.
    \end{align*}
\end{lemma}

\textbf{Proof for (ii):} 
We have
\begin{align*}
    \|\nabla f(x)\|^2-\|\nabla_x F(x,t)\|^2 &= (\|\nabla f(x)\|+\|\nabla_x F(x,t)\|)(\|\nabla f(x)\|-\|\nabla_x F(x,t)\|)\\
    &\leq 2L_0(\|\nabla f(x)\|-\|\nabla_x F(x,t)\|)\\
    &\leq 2L_0\|\nabla_x F(x,t) - \nabla f(x)\|.
\end{align*}
The term $\|\nabla_x F(x,t) - \nabla f(x)\|$ can be upper bounded as follows:
\begin{align*}
    \|\nabla_x F(x,t) - \nabla f(x)\| &\leq \left\|\mathbb{E}_u\left[\left(\frac{f(x+tu)-f(x)}{t}-\langle \nabla f(x), u\rangle\right) u\right]\right\|\\
    &\leq \mathbb{E}_u\left[\left|\frac{1}{t}\left(f(x+tu)-f(x)-t\langle \nabla f(x), u\rangle\right)\right|\|u\|\right]\\
    &\leq \mathbb{E}_u\left[\frac{tL_1}{2}\|u\|^3\right]\\
    &\leq \frac{tL_1}{2}(d+3)^{3/2},
\end{align*}
where the last second inequality follows from a property of $L_1$-smooth function ($\forall x,y\in\mathbb{R}^d,\ |f(y)-f(x)-\langle \nabla f(x), y-x\rangle|\leq \frac{L_1}{2}\|y-x\|^2$), and the last inequality holds due to Lemma \ref{lem:gauss_norm}. Therefore, we obtain 
\begin{align*}
        \|\nabla f(x)\|^2 \leq \|\nabla_x F(x,t)\|^2 + tL_0L_1(d+3)^{3/2}.
\end{align*}

Now, we are ready to prove Theorem \ref{iter_determin}.

\textbf{Proof for Theorem \ref{iter_determin}:}
We follow the convergence analysis of gradient descent. According to Assumption \ref{A1} and Lemma \ref{lem:Lip}, $F(x,t)$ is  $L_0$-$Lipschitz$ and $L_1$-$smooth$ in terms of $x$. Therefore, we have
\begin{align}
F(x_{k+1},t_k) &\leq F(x_k,t_k) + \left<\nabla_x F(x_k,t_k),(x_{k+1}-x_k)\right> + \frac{L_1}{2}\|x_{k+1}-x_k\|^2 \nonumber\\ 
&= F(x_k,t_k) - \left(\beta-\frac{L_1}{2}\beta^2\right)\|\nabla_x F(x_k,t_k)\|^2\nonumber,
\end{align}
where the last equation holds due to the updating rule of the gradient descent: $x_{k+1} - x_k = -\beta\nabla_x F(x_k,t_k)$. Then, we can get the upper bound for $\|\nabla_x F(x,t)\|^2$:
\begin{align}
\left(\beta-\frac{L_1}{2}\beta^2\right)\|\nabla_x F(x_k,t_k)\|^2 &\leq F(x_k,t_k) - F(x_{k+1},t_k)\nonumber\\
&= F(x_k,t_k) - F(x_{k+1},t_{k+1}) + F(x_{k+1},t_{k+1}) - F(x_{k+1},t_k)\nonumber\\
&\leq F(x_k,t_k) - F(x_{k+1},t_{k+1}) + L_0|t_{k+1}-t_k|\sqrt{d},\nonumber
\end{align}
where the last inequality follows from Lemma \ref{lem:Lip_t}.\\
Now, sum up the above inequality for all iterations  $k_0+1\leq k\leq T\ (T>k_0\in\mathbb{N})$, and denote the minimum of $f$ as $f^*$, then we have
\begin{align}
\left(\beta-\frac{L_1}{2}\beta^2\right)\sum_{k=k_0+1}^{T}\|\nabla_x F(x_k,t_k)\|^2 &\leq F(x_{k_0+1},t_{k_0+1}) - F(x_{T+1},t_{T+1}) + L_0\sqrt{d}\sum_{k=k_0+1}^T|t_{k+1}-t_k|\nonumber\\
&\leq F(x_{k_0+1},t_{k_0+1}) - f^* + L_0\sqrt{d}\sum_{k={k_0+1}}^T|t_{k+1}-t_k|\nonumber\\
&\leq f(x_{k_0+1}) - f^* + L_0\sqrt{d}\left(t_{k_0+1} + \sum_{k=k_0+1}^T|t_{k+1}-t_k|\right)
\label{normGradF},
\end{align}
where the last inequality holds due to Lemma \ref{lem:Lip_t}. 
Then,
we can get the upper bound for $\|\nabla f(\hat{x})\|^2$ as
\begin{align*}
    & \| \nabla f(\hat{x}) \|^2
    =
    {\mathop{\rm min}\limits_{k\in [T]}}\|\nabla f(x_k)\|^2 
    \\
    &
    \leq 
    {\mathop{\rm min}\limits_{k = k_0+1, \ldots, T}}\|\nabla f(x_k)\|^2 
    \\
    &\leq \frac{1}{T-k_0}\sum_{k=k_0+1}^{T}\|\nabla f(x_k)\|^2\nonumber
    \\
    &\leq \frac{1}{T-k_0}\sum_{k=k_0+1}^{T}\|\nabla_x F(x_k,t_k)\|^2 + \frac{1}{T-k_0}L_0L_1(d+3)^{3/2}\sum_{k=k_0+1}^Tt_k\nonumber\\
    &\leq \frac{2\left(f(x_{k_0+1})-f^*+L_0\sqrt{d}\left(t_{k_0+1} + \sum_{k=k_0+1}^T|t_{k+1}-t_k|\right)\right)}{(T-k_0)(2\beta-L_1\beta^2)} + \frac{1}{T-k_0}L_0L_1(d+3)^{3/2}\sum_{k=k_0+1}^Tt_k,
\end{align*}
where the third inequality holds due to Lemma \ref{lem:grad_square_diff} (ii) and the last inequality follows from \eqref{normGradF}.

If we choose the step size $\beta$ as $\frac{1}{L_1}$, we have
\begin{align}
    \nonumber
    & \|\nabla f(\hat{x})\|^2\\
    \nonumber &\leq
\frac{2L_1\left(f(x_{k_0+1})-f^*+L_0\sqrt{d}\left(t_{k_0+1} + \sum_{k=k_0+1}^T|t_{k+1}-t_k|\right)\right)}{T-k_0} + \frac{1}{T-k_0}L_0L_1(d+3)^{3/2}\sum_{k=k_0+1}^Tt_k
\\
    &= O\left(\frac{1}{T - k_0}\left(1+ d^{3/2}\sum_{k= k_0 + 1}^Tt_k\right)\right),
    \label{eq:boundsquarednorm_iter_determ}
\end{align}
where the last equality holds since $\sum_{k= k_0 + 1}^T|t_{k+1}-t_k|=O\left(\sum_{k= k_0 + 1}^Tt_k\right)$ is satisfied.
If we update $t_k$ as in Algorithm~\ref{alg:GH},
we have
% \begin{align}
$
    \sum_{k= k_0 + 1}^T t_{k}
    \leq
    \sum_{k= k_0 + 1}^T \max\{ t_1 \gamma^{k-1}, \epsilon'\}
    \leq
    \sum_{k= k_0 + 1}^T \left( t_1 \gamma^{k-1} + \epsilon' \right)
    \leq
    \frac{t_{1}\gamma^{k_0}}{1 - \gamma} + \epsilon' (T-k_0).
$
By taking $\epsilon'$ sufficiently close to $0$, together with the assumption of $1/(1-\gamma) = O(1)$, we have $\sum_{k= k_0 + 1}^T t_k =O(\gamma^{k_0})$. This implies that
$\| \nabla f(\hat{x}) \|^2 \leq O( \frac{1 + \gamma^{k_0}d^{3/2}}{T - k_0} )$.
Hence,
we can obtain
$\| \nabla f(\hat{x}) \| \leq \epsilon$ in $T = k_0 + O\left(\frac{1 + \gamma^{k_0} d^{3/2}}{\epsilon^2}\right)$ iterations.

Now, set $k_0$ as $k_0=O\left(\frac{1}{\epsilon^2}\right)$, then, the iteration complexity can be bounded as $T=O\left(\frac{d^{3/2}}{\epsilon^2}\right)$. Furthermore, when $\gamma$ is chosen as $\gamma\leq d^{-3\epsilon^2/2}$, we can obtain
$
\gamma^{k_0} = O\left(d^{-3/2}\right)
$
for some $k_0=O\left( \frac{1}{\epsilon^2} \right)$. This yields the iteration complexity of $T=O\left(\frac{1}{\epsilon^2}\right)$.

\hfill $\Box$

Before going to the proof of Theorem \ref{thm:iter_stochastic} in the stochastic setting, we prove that the gradient of the smoothed stochastic function $\nabla F(x,t;\xi)$ is unbiased, and it has a finite variance.

\begin{lemma}
\label{lem:grad_smoothed_and_stochastic}
Suppose that $f$ satisfies Assumption \ref{A1} (i) and Assumption \ref{A2}.\\
\textbf{(i)} 
The stochastic gradient of the smoothed function $\nabla_x \bar{F}(x,t;\xi)$ becomes an unbiased estimator of $\nabla_x F(x,t)$. That is, for any $x\in\mathbb{R}^d$ and $t>0$, $\mathbb{E}_\xi[\nabla_x\bar{F}(x,t;\xi)]=\nabla_x F(x,t)$ holds.\\
\textbf{(ii)} 
For any $x\in\mathbb{R}^d$ and $t>0$, the variance of $\nabla_x \bar{F}(x,t;\xi)$ is bounded as $\mathbb{E}_{\xi}[\|\nabla_x\bar{F}(x,t;\xi)-\nabla_x F(x,t)\|^2] \leq \sigma^2$.
\end{lemma}
\textbf{Proof for (i):}
From Assumption~\ref{A1}~(i), we can exchange the order of integration in terms of $\xi$ and $u$, which yields that 
\begin{align*}
    \mathbb{E}_\xi[\nabla_x \bar{F}(x,t;\xi)]
    &= \mathbb{E}_\xi \left[ \mathbb{E}_u \left[ \frac{\bar{f}(x+tu;\xi) - \bar{f}(x;\xi)}{t}u \right]\right] \\
    &= \mathbb{E}_u \left[ \mathbb{E}_\xi \left[ \frac{\bar{f}(x+tu;\xi)- \bar{f}(x;\xi)}{t}u \right]\right] \\
    &= \mathbb{E}_u \left[ \frac{f(x+tu) - f(x)}{t}u \right] \\
    &= \nabla_x F(x,t).
\end{align*}

\textbf{Proof for (ii):}
We have
\begin{align*}
    \mathbb{E}_\xi[\|\nabla_x \bar{F}(x,t;\xi)-\nabla_xF(x,t)\|^2] &= \mathbb{E}_\xi[\|\nabla_x \mathbb{E}_u[\bar{f}(x+tu;\xi)]-\nabla_x\mathbb{E}_u[f(x+tu)]\|^2]\\
    &= \mathbb{E}_\xi[\|\mathbb{E}_u[\nabla_x\bar{f}(x+tu;\xi) - \nabla f(x+tu)]\|^2]\\
    &\leq \mathbb{E}_\xi[\mathbb{E}_u[\|\nabla_x\bar{f}(x+tu;\xi) - \nabla f(x+tu)\|^2]]\\
    &= \mathbb{E}_u[\mathbb{E}_\xi[\|\nabla_x\bar{f}(x+tu;\xi) - \nabla f(x+tu)\|^2]]\\
    &\leq \sigma^2,
\end{align*}
where the second and third equalities hold due to Assumption~\ref{A1}~(i), and the last inequality follows from Assumption~\ref{A2}~(ii).

\textbf{Proof for Theorem \ref{thm:iter_stochastic}:}
Denote $\delta_k := \nabla_x \bar{F}(x_k,t_k;\xi_k)-\nabla_x F(x_k,t_k)$.
We follow the convergence analysis of stochastic gradient descent. According to Lemma \ref{lem:Lip}, since $f(x)$ is $L_0$-$Lipschitz$ and $L_1$-$smooth$, $F(x,t)$ is also $L_0$-$Lipschitz$ and $L_1$-$smooth$ in terms of $x$.
Thus, we have
\begin{align}
F(x_{k+1},t_k) &\leq F(x_k,t_k) + \left<\nabla_x F(x_k,t_k),(x_{k+1}-x_k)\right> + \frac{L_1}{2}\|x_{k+1}-x_k\|^2\nonumber\\ 
&=  F(x_k,t_k) - \beta \left<\nabla_x F(x_k,t_k), \nabla_x \bar{F}(x_k,t_k;\xi_k)\right> + \frac{L_1}{2}\beta^2\|\nabla_x \bar{F}(x_k,t_k;\xi_k)\|^2\nonumber\\
&=  F(x_k,t_k) - \left(\beta-\frac{L_1}{2}\beta^2\right)\|\nabla_x F(x_k,t_k)\|^2 - (\beta-L_1\beta^2) \left<\nabla_x F(x_k,t_k),\delta_k\right> + \frac{L_1}{2}\beta^2\|\delta_k\|^2,
\label{4}
\end{align}
where the first equation holds due to the updating rule $x_{k+1} - x_k = -\beta \nabla_x \bar{F}(x_k,t_k;\xi_k)$, and the last equation holds due to the definition of $\delta_k$. Denote
\[
A_k := - (\beta-L_1\beta^2) \left<\nabla_x F(x_k,t_k),\delta_k\right> + \frac{L_1}{2}\beta^2\|\delta_k\|^2
\]
for simplicity. From \eqref{4}, we obtain the upper bound for $\|\nabla_x F(x,t)\|^2$ as follows:
\begin{align}
\left(\beta-\frac{L_1}{2}\beta^2\right)\|\nabla_x F(x_k,t_k)\|^2 &\leq F(x_k,t_k) - F(x_{k+1},t_k) + A_k\nonumber\\
&= F(x_k,t_k) - F(x_{k+1},t_{k+1}) + F(x_{k+1},t_{k+1}) - F(x_{k+1},t_k) + A_k\nonumber\\
&\leq F(x_k,t_k) - F(x_{k+1},t_{k+1}) + L_0|t_{k+1}-t_k|\sqrt{d} + A_k,\nonumber
\end{align}
where the last inequality follows from Lemma \ref{lem:Lip_t}.\\
Now, sum up the above inequality for all iterations $k_0 + 1\leq k\leq T\ (k_0<T)$. Then we have
\begin{align}
&\left(\beta-\frac{L_1}{2}\beta^2\right)\sum_{k=k_0 + 1}^{T}\|\nabla_x F(x_k,t_k)\|^2 \nonumber \\
&\leq F(x_{k_0 + 1},t_{k_0 + 1}) - F(x_{T+1},t_{T+1}) + L_0\sqrt{d}\sum_{k=k_0 + 1}^T |t_{k+1}-t_k| + \sum_{k=k_0 + 1}^{T} A_k\nonumber\\
&\leq F(x_{k_0 + 1},t_{k_0 + 1}) - f^* + L_0\sqrt{d}\sum_{k=k_0 + 1}^T |t_{k+1}-t_k| + \sum_{k=k_0 + 1}^{T} A_k\nonumber.\\
&\leq f(x_{k_0 + 1}) - f^* + L_0\sqrt{d}\left(t_{k_0 + 1}+\sum_{k=k_0 + 1}^T |t_{k+1}-t_k|\right) + \sum_{k=k_0 + 1}^{T} A_k\nonumber.
\end{align}
Take the expectation with respect to the random vectors $\{ \xi_{k_0+1}, \ldots, \xi_T \}$, then we have
\begin{align}
&\left(\beta-\frac{L_1}{2}\beta^2\right)\sum_{k=k_0 + 1}^{T}\mathbb{E}_{\xi}[\|\nabla_x F(x_k,t_k)\|^2] \nonumber
\\
&\leq f(x_{k_0 + 1}) - f^* + L_0\sqrt{d}\left(t_{k_0 + 1}+\sum_{k=k_0 + 1}^T \mathbb{E}_{\xi}[|t_{k+1}-t_k|]\right) + \sum_{k=k_0 + 1}^{T} \mathbb{E}_{\xi}[A_k].
\label{5}
\end{align}
The expectation of $A_k$ is evaluated as 
\begin{align}
\sum_{k=k_0 + 1}^{T} \mathbb{E}_{\xi}[A_k] &= - \sum_{k=k_0 + 1}^{T}(\beta-L_1\beta^2)\mathbb{E}_{\xi}[ \left<\nabla_x F(x_k,t_k),\delta_k\right>] + \sum_{k=k_0 + 1}^{T}\frac{L_1}{2}\beta^2\mathbb{E}_{\xi}[\|\delta_k\|^2]\nonumber\\
&\leq (T - k_0) \frac{L_1}{2}\beta^2\sigma^2,
\label{6}
\end{align}
where the last equality holds due to Lemma~\ref{lem:grad_smoothed_and_stochastic}~(ii) ($\mathbb{E}_{\xi}[\|\delta_k\|^2]\leq \sigma^2$) and the fact that each point $x_k$ is a function of the history $\xi_{[k-1]}$ in the random process, thus $\mathbb{E}_{\xi_k}[ \left<\nabla_x F(x_k,t_k),\delta_k\right>\mid\xi_{[k-1]}] = 0$.

%Next, we estimate the upper bound for $\frac{1}{T}\sum_{k=1}^{T}t_k^2$, which is useful for the following proof.
%%\begin{align}
%%\frac{1}{T}\sum_{k=1}^{T}t_k^2 &\leq \frac{1}{T}{(t_1^2+\gamma^2 t_1^2 + \gamma^4 t_1^2 + ... + \gamma^{2(T-1)}t_1^2)}\nonumber\\
%%&= \frac{t_1^2(1-\gamma^{2T})}{T(1-\gamma^2)}
%%\label{25}
%%\end{align}

Then,
we can estimate the upper bound for $\mathbb{E}_{\xi,k'}[\|\nabla f(\hat{x})\|^2]$ as 
\begin{align*}
&\mathbb{E}_{\xi, k'}[\|\nabla f(\hat{x})\|^2]=\frac{1}{T - k_0}\sum_{k=k_0 + 1}^{T}\mathbb{E}_{\xi}[\|\nabla f(x_k)\|^2] 
\\
&
\leq \frac{1}{T - k_0}\sum_{k=k_0 + 1}^{T}\mathbb{E}_{\xi}[\|\nabla_x F(x_k,t_k)\|^2] +\frac{1}{T - k_0}L_0L_1(d+3)^{3/2}\sum_{k=k_0 + 1}^{T}\mathbb{E}_{\xi}[t_k]\\
&\leq \frac{2\left(f(x_{k_0 + 1})-f^*+L_0\sqrt{d}\left(t_{k_0 + 1}+\sum_{k=k_0 + 1}^T \mathbb{E}_{\xi}[|t_{k+1}-t_k|]\right)\right)}{(T - k_0)(2\beta-L_1\beta^2)} \\ 
&+
\frac{1}{T - k_0}L_0L_1(d+3)^{3/2}\sum_{k=k_0 + 1}^{T}\mathbb{E}_{\xi}[t_k]+  \frac{L_1\beta^2\sigma^2}{2\beta-L_1\beta^2},
\end{align*}
where the first inequality holds due to Lemma \ref{lem:grad_square_diff} (ii) and the last inequality follows from  \eqref{5} and \eqref{6}.

If the step size $\beta$ is chosen as $\beta = \mathop{\rm min}\ \{\frac{1}{L_1}, \frac{1}{\sqrt{T - k_0}}\}$, then we have
\begin{align*}
    \frac{1}{2\beta-L_1\beta^2} \leq \frac{1}{\beta},
\end{align*}
\begin{align*}
    \frac{1}{\beta} \leq L_1 + \sqrt{T - k_0}.
\end{align*}
Hence, we can obtain
\begin{align*}
    & \frac{2\left(f(x_{k_0 + 1})-f^*+L_0\sqrt{d}\left(t_{k_0 + 1}+\sum_{k=k_0 + 1}^T \mathbb{E}_{\xi}[|t_{k+1}-t_k|]\right)\right)}{(T - k_0)(2\beta-L_1\beta^2)}\\
    & + \frac{1}{T - k_0}L_0L_1(d+3)^{3/2}\sum_{k=k_0 + 1}^{T}\mathbb{E}_{\xi}[t_k]  +  \frac{L_1\beta^2\sigma^2}{2\beta-L_1\beta^2}\\
    & = O \left(\frac{1+\sqrt{d}\mathbb{E}_{\xi}\left[\sum_{k=k_0 + 1}^T|t_{k+1}-t_k|\right]}{\sqrt{T - k_0}}+\frac{d^{3/2}}{T - k_0}\mathbb{E}_{\xi}\left[\sum_{k=k_0 + 1}^{T}t_k\right]\right).
\end{align*}
If $t_k$ is updated as in Algorithm~\ref{alg:GH},
we have
$
    \sum_{k=k_0 + 1}^T |t_{k+1} - t_k| \leq t_1\gamma_{k_0} = O(\gamma^{k_0})
$
and 
$
    \sum_{k=k_0 + 1}^T t_{k} \leq \frac{t_1\gamma^{k_0}}{1-\gamma} + \epsilon'T = O(\gamma^{k_0})
$ in the same argument that showed Theorem~\ref{iter_determin}.
Combining the above inequalities,
we obtain
\begin{align}
    \mathbb{E}_{\xi, k'}[\|\nabla f(\hat{x})\|^2]
    =
    \frac{1}{T - k_0}\sum_{k=k_0 + 1}^{T}\mathbb{E}_{\xi}[\|\nabla f(x_k)\|^2] 
    =
    O \left( \frac{1 + \sqrt{d}\gamma^{k_0}}{\sqrt{T - k_0}} + \frac{d^{3/2} \gamma^{k_0}}{T - k_0} \right).
    \label{eq:boundsquarednorm_iter_stoc}
\end{align}

Here, we have 
$
    k_0 = O \left(
        \frac{1}{\epsilon^4}
    \right)
$
by the definition of $k_0$.
Thus, by setting
$
    T
    = k_0 + O \left(
        \frac{d}{\epsilon^4} + \frac{d^{3/2}}{\epsilon^2}
    \right)
    = 
    O \left(
        \frac{d}{\epsilon^4} + \frac{d^{3/2}}{\epsilon^2}
    \right)
$,
we can obtain
$
    \mathbb{E}_{\xi, k'}[\|\nabla f(\hat{x})\|^2] \leq \epsilon^2
$.
% As we have
This implies 
$
    \mathbb{E}_{\xi, k'}[\|\nabla f(\hat{x})\|] \leq \epsilon
$
as
$
    \mathbb{E}_{\xi, k'}[\|\nabla f(\hat{x})\|]^2 
    \leq
    \mathbb{E}_{\xi, k'}[\|\nabla f(\hat{x})\|^2]
$
follows from Jensen's inequality. Furthermore, when $\gamma$ is chosen as
$\gamma \leq (\max\{ d^{1/2},  d^{3/2} \epsilon^{2} \})^{- \epsilon^4}$
, we have 
$
\log_\gamma \min\{ d^{-1/2}, d^{-3/2} \epsilon^{-2} \} = O\left(\frac{1}{\epsilon^4}\right)
$, which implies $k_0 = \Omega(\log_\gamma \min\{ d^{-1/2}, d^{-3/2} \epsilon^{-2} \}).$
Therefore, we can obtain
$
\gamma^{k_0} = O\left(\min\{ d^{-1/2}, d^{-3/2} \epsilon^{-2} \})\right)
$
, which yields the iteration complexity of $T=O\left(\frac{1}{\epsilon^4}\right)$.

\hfill $\Box$

\subsection{Zeroth-order SLGH algorithm}
\label{subsec:proof_zeroth_order}
In the zeroth-order setting, we can evaluate the gap between the zeroth-order gradient estimator and the true gradient using the following lemma.\\
\begin{lemma}
    [\textbf{Theorem 4 in \cit{nesterov2017random}}]\label{lem:zo_grad_square_diff}
    Let $f$ be a $L_1$-$smooth$ function, then for any $x\in \mathbb{R}^d$ and for any $t>0$, we have
    \begin{align*}
        \mathbb{E}_u\left[\frac{1}{t^2}(f(x+tu)-f(x))^2\|u\|^2\right] \leq \frac{t^2}{2}L_1^2(d+6)^3 + 2(d+4)\|\nabla f(x)\|^2.
    \end{align*}
\end{lemma}

\textbf{Proof for Theorem \ref{thm:zo_determinstic}:}
Let $w_k:=(u_k, v_k),\ k\in[T]$, and denote $\delta_k := \Tilde{g}_x(x_k,t_k;u_k)-\nabla_x F(x_k,t_k)$, where $\Tilde{g}_x(x_k,t_k;u_k)$ is the zeroth-order estimator of gradient defined in the main paper. Utilize the updating rule of $x$ and $L_1$-smoothness of $F(x,t)$ in terms of $x$. Then we have
\begin{align}
F(x_{k+1},t_k) &\leq F(x_k,t_k) + \left<\nabla_x F(x_k,t_k),(x_{k+1}-x_k)\right> + \frac{L_1}{2}\|x_{k+1}-x_k\|^2\nonumber\\
&= F(x_k,t_k) - \beta \left<\nabla_x F(x_k,t_k),\Tilde{g}_x(x_k,t_k;u_k)\right> + \frac{L_1}{2}\beta^2\|\Tilde{g}_x(x_k,t_k;u_k)\|^2\nonumber\\
&= F(x_k,t_k) - \beta\|\nabla_x F(x_k,t_k)\|^2 - \beta \left<\nabla_x F(x_k,t_k),\delta_k\right> + \frac{L_1}{2}\beta^2\|\Tilde{g}_x(x_k,t_k;u_k))\|^2,
\label{7}
\end{align}
where the first equation holds due to the updating rule $x_{k+1} - x_k = -\beta \Tilde{g}_x(x_k,t_k;u_k)$.\\
Denote
\[B_k := - \beta \left<\nabla_x F(x_k,t_k),\delta_k\right> + \frac{L_1}{2}\beta^2\|\Tilde{g}_x(x_k,t_k;u_k)\|^2
  \]
  for simplicity. From Lemma \ref{lem:Lip_t} and \eqref{7}, we get the upper bound for $\|\nabla_x F(x,t)\|^2$ as 
\begin{align}
\beta \|\nabla_x F(x_k,t_k)\|^2 &\leq F(x_k,t_k) - F(x_{k+1},t_k) + B_k\nonumber\\
&= F(x_k,t_k) - F(x_{k+1},t_{k+1}) + F(x_{k+1},t_{k+1}) - F(x_{k+1},t_k) + B_k\nonumber\\
&\leq F(x_k,t_k) - F(x_{k+1},t_{k+1}) + L_0|t_{k+1}-t_k|\sqrt{d} + B_k\nonumber.
\end{align}
Now, sum up the above inequality for all iterations $k_0+1\leq k\leq T\ (k_0<T)$. Then we have
\begin{align}
\sum_{k=k_0+1}^{T}\beta\|\nabla_x F(x_k,t_k)\|^2 &\leq F(x_{k_0 + 1},t_{k_0 + 1}) - F(x_{T+1},t_{T+1}) + L_0\sum_{k=k_0 + 1}^T |t_{k+1}-t_k|\sqrt{d} + \sum_{k=k_0 + 1}^{T}B_k\nonumber\\
&\leq F(x_{k_0 + 1},t_{k_0 + 1}) - f^* + L_0\sqrt{d}\sum_{k=k_0 + 1}^T |t_{k+1}-t_k| + \sum_{k=k_0 + 1}^{T}B_k\nonumber\\
&\leq f(x_{k_0 + 1}) - f^* + L_0\sqrt{d}\left(t_{k_0 + 1}+\sum_{k=k_0 + 1}^T |t_{k+1}-t_k|\right) + \sum_{k=k_0 + 1}^{T}B_k\nonumber.
\end{align}
Next, take the expectations with respect to random vectors $\{w_{k_0+1}, \ldots, w_T \}$ on both sides. Then we can get
\begin{align}
\sum_{k=k_0 + 1}^{T}\beta \mathbb{E}_{w}[\|\nabla_x F(x_k,t_k)\|^2] &\leq f(x_{k_0 + 1}) - f^* + L_0\sqrt{d}\left(t_{k_0 + 1}+\sum_{k=k_0 + 1}^T \mathbb{E}_w[|t_{k+1}-t_k|]\right)\nonumber \\
&+ \sum_{k=k_0 + 1}^{T} \mathbb{E}_{w} [B_k].
\label{33}
\end{align}
Observe by the definition of $\Tilde{g}_x(x_k,t_k;u_k)$ in the main paper that $\mathbb{E}_{u_k}[\Tilde{g}_x(x_k,t_k;u_k)\mid u_{[k-1]}] = \nabla_x F(x_k,t_k)$, thus $\mathbb{E}_{w_k}[\left<\nabla_x F(x_k,t_k),\delta_k\right>\mid w_{[k-1]}] = 0$ holds. 
Then we have
\begin{align}\label{B_k}
    \mathbb{E}_{w_k} [B_k\mid w_{[k-1]}] 
    &= -\beta \mathbb{E}_{w_k}[ \left<\nabla_x F(x_k,t_k),\delta_k\right>\mid w_{[k-1]}] + \frac{L_1}{2}\beta^2\mathbb{E}_{w_k}[\|\Tilde{g}_x(x_k,t_k;u_k)\|^2\mid w_{[k-1]}]\nonumber\\
    &\leq \frac{L_1}{2}\beta^2\left(\frac{\mathbb{E}_{w_k} [t_k^2\mid w_{[k-1]}]}{2}L_1^2(d+6)^3+2(d+4) \mathbb{E}_{w_k}[\|\nabla f(x_k)\|^2\mid w_{[k-1]}]\right)\nonumber\\
    &= \frac{\mathbb{E}_{w_k} [t_k^2\mid w_{[k-1]}]}{4}L_1^3\beta^2(d+6)^3+L_1\beta^2(d+4) \mathbb{E}_{w_k}[\|\nabla f(x_k)\|^2\mid w_{[k-1]}],
\end{align}
where the inequality holds due to Lemma \ref{lem:zo_grad_square_diff}.

Lemma \ref{lem:grad_square_diff} (ii) together with the above inequalities yields that 
\begin{align}
&\sum_{k=k_0 + 1}^{T}\beta \mathbb{E}_{w}[\|\nabla f(x_k)\|^2]\nonumber \\
&\leq \sum_{k=k_0 + 1}^{T}\beta \mathbb{E}_{w}[\|\nabla_x F(x_k,t_k)\|^2] +  \sum_{k=k_0 + 1}^{T}\beta L_0L_1(d+3)^{3/2}\mathbb{E}_{w}[t_k]\nonumber\\
&\leq 
f(x_{k_0 + 1}) - f^* + L_0\sqrt{d}\left(t_{k_0 + 1}+\sum_{k=k_0 + 1}^T \mathbb{E}_w[|t_{k+1}-t_k|]\right) + \sum_{k=k_0 + 1}^{T} \mathbb{E}_{w} [B_k]\nonumber\\
&+
\sum_{k=k_0 + 1}^{T}\beta L_0L_1(d+3)^{3/2}\mathbb{E}_{w}[t_k]\nonumber\\
&\leq 
f(x_{k_0 + 1}) - f^* + L_0\sqrt{d}\left(t_{k_0 + 1}+\sum_{k=k_0 + 1}^T \mathbb{E}_w[|t_{k+1}-t_k|]\right)+\sum_{k=k_0 + 1}^{T}\frac{\mathbb{E}_w[t_k^2]}{4}L_1^3\beta^2(d+6)^3\nonumber\\
&+\sum_{k=k_0 + 1}^{T}L_1\beta^2(d+4) \mathbb{E}_{w}[\|\nabla f(x_k)\|^2] +\sum_{k=k_0 + 1}^{T}\beta L_0L_1(d+3)^{3/2}\mathbb{E}_{w}[t_k],
\end{align}
where the second inequality holds due to \eqref{33}, and the last inequality follows from \eqref{B_k}.
Rearrange the terms in the above inequality. Then we can get
\begin{align}
(\beta-(d+4)L_1\beta^2)\sum_{k=k_0 + 1}^{T}\mathbb{E}_{w}[\|\nabla f(x_k)\|^2] &\leq f(x_{k_0 + 1}) - f^* + L_0\sqrt{d}\left(t_{k_0 + 1}+\sum_{k=k_0 + 1}^T \mathbb{E}_w[|t_{k+1}-t_k|]\right)\nonumber\\&+\frac{L_1^3\beta^2(d+6)^3}{4}\sum_{k=k_0 + 1}^{T}\mathbb{E}_w[t_k^2]+ L_0L_1\beta(d+3)^{3/2}\sum_{k=k_0 + 1}^{T}\mathbb{E}_{w}[t_k].
\end{align}
%\memo{$\frac{t_1^2}{(1-\gamma^2)}$ need to be divided by 2?} 
Divide both sides of the above inequality by $(T - k_0)(\beta-(d+4)L_1\beta^2)$ and set the step size $\beta$ as $\frac{1}{2(d+4)L_1}$. Since $\frac{1}{\beta-(d+4)L_1\beta^2} \leq 4(d+4)L_1$ holds, we can obtain
\begin{align}
\frac{1}{T - k_0}\sum_{k=k_0 + 1}^{T}\mathbb{E}_{w}[\|\nabla f(x_k)\|^2]
&\leq \frac{4(d+4)L_1}{T - k_0} \left(f(x_{k_0 + 1}) - f^* + L_0\sqrt{d}\left(t_{k_0 + 1}+\sum_{k=k_0 + 1}^T \mathbb{E}_w[|t_{k+1}-t_k|]\right)\right.\nonumber\\
&\left.+\frac{L_1(d+6)^3}{16(d+4)^2}\sum_{k=k_0 + 1}^{T}\mathbb{E}_w[t_k^2]+ \frac{L_0(d+3)^{3/2}}{2(d+4)}\sum_{k=k_0 + 1}^{T}\mathbb{E}_{w}[t_k]\right)\nonumber\\
&= O\left(\frac{d}{T - k_0}\left(1+d \mathbb{E}_w\left[ \sum_{k=k_0 + 1}^T t_k^2 \right] + \sqrt{d} \mathbb{E}_w \left[\sum_{k=k_0 + 1}^T t_k \right]\right)\right) \nonumber
\\
&
=O\left(
    \frac{d}{T - k_0}
    \left(1 + d\gamma^{2k_0} + \sqrt{d}\gamma^{k_0} \right)
\right),
\end{align}
where the last equality follows from the update rule of $t_k$,
as shown in the proof of Theorem~\ref{iter_determin} as well.

Here, we have
$
    k_0 = O \left(
        \frac{d}{\epsilon^2}
    \right)
$
by the definition of $k_0$. Thus, by setting
$
    T
    = k_0 + O \left(
        \frac{d^2}{\epsilon^2}
    \right)
    = 
    O \left(
        \frac{d^2}{\epsilon^2}
    \right)
$,
we can obtain
$
    \mathbb{E}_{w, k'}[\|\nabla f(\hat{x})\|^2]
    =
    \frac{1}{T - k_0}\sum_{k=k_0 + 1}^{T}\mathbb{E}_{w}[\|\nabla f(x_k)\|^2]
    \leq \epsilon^2
$.
% As we have
This implies 
$
    \mathbb{E}_{w, k'}[\|\nabla f(\hat{x})\|] \leq \epsilon
$
as
$
    \mathbb{E}_{w, k'}[\|\nabla f(\hat{x})\|]^2 
    \leq
    \mathbb{E}_{w, k'}[\|\nabla f(\hat{x})\|^2]
$
follows from Jensen's inequality. Furthermore, when $\gamma$ is chosen as
$\gamma \leq d^{-\epsilon^2/2d}$
, we have 
$
    \log_\gamma d^{-1/2} = O\left( \frac{d}{\epsilon^2} \right)
$,
which implies 
$
 k_0 = \Omega \left( \log_\gamma d^{-1/2} \right)
$. Therefore, we can obtain
$
\gamma^{k_0}
=
O( d^{-1/2} )
$,
which yields the iteration complexity of $T=O\left(\frac{d}{\epsilon^2}\right)$.

\hfill $\Box$

\textbf{Proof for Theorem \ref{thm:zo_stochastic}:}
Let $\zeta_k := (\xi_k,u_k, v_k)$, $k\in [T]$ and denote $\delta_k := \Tilde{G}_x(x_k,t_k;\xi_k,u_k)-\nabla_x F(x_k,t_k)$.
%\[
%\equiv \Tilde{G}_x(x_k,t_k;\zeta_k)-\nabla_x F(x_k,t_k).
%\]
As discussed in the main paper, we have \begin{align}
    \mathbb{E}_{\xi,u} [\Tilde{G}_{x}(x,t;\xi,u)] =  \mathbb{E}_u [\mathbb{E}_\xi [\Tilde{G}_{x}(x,t;\xi,u)|u]] = \nabla_x F(x,t).
    \label{expuxi}
\end{align}
From the update rule for $x$, we can obtain
\begin{align}
F(x_{k+1},t_k) &\leq F(x_k,t_k) + \left<\nabla_x F(x_k,t_k),(x_{k+1}-x_k)\right> + \frac{L_1}{2}\|x_{k+1}-x_k\|^2\nonumber\\
&= F(x_k,t_k) - \beta \left<\nabla_x F(x_k,t_k),\Tilde{G}_x(x_k,t_k;\xi_k, u_k)\right> + \frac{L_1}{2}\beta^2\|\Tilde{G}_x(x_k,t_k;\xi_k, u_k)\|^2\nonumber\\
&= F(x_k,t_k) - \beta\|\nabla_x F(x_k,t_k)\|^2 - \beta \left<\nabla_x F(x_k,t_k),\delta_k\right> + \frac{L_1}{2}\beta^2\|\Tilde{G}_x(x_k,t_k;\xi_k, u_k)\|^2\nonumber.
\end{align}
Now, denote
\[
D_k := - \beta  \left<\nabla_x F(x_k,t_k),\delta_k\right> + \frac{L_1}{2}\beta^2\|\Tilde{G}_x(x_k,t_k;\xi_k, u_k)\|^2
\]
for simplicity. Then, we can get the upper bound for $\|\nabla_x F(x,t)\|^2$ with $D_k$:
\begin{align}
\beta \|\nabla_x F(x_k,t_k)\|^2 &\leq F(x_k,t_k) - F(x_{k+1},t_k) + D_k\nonumber\\
&= F(x_k,t_k) - F(x_{k+1},t_{k+1}) + F(x_{k+1},t_{k+1}) - F(x_{k+1},t_k) + D_k\nonumber\\
&\leq F(x_k,t_k) - F(x_{k+1},t_{k+1}) + L_0|t_{k+1}-t_k|\sqrt{d} + D_k.\nonumber
\end{align}
Sum up the above inequality for all iterations $k_0 + 1\leq k\leq T\ (T>k_0)$. Then we have
\begin{align}
&\sum_{k=k_0 + 1}^{T}\beta\|\nabla_x F(x_k,t_k)\|^2 \nonumber \\
&\leq F(x_{k_0 + 1},t_{k_0 + 1}) - F(x_{T+1},t_{T+1}) + L_0\sqrt{d}\sum_{k=k_0 + 1}^T|t_{k+1}-t_k| + \sum_{k=k_0 + 1}^{T}D_k\nonumber\\
&\leq F(x_{k_0 + 1},t_{k_0 + 1}) - f^* + L_0\sqrt{d}\sum_{k=k_0 + 1}^T|t_{k+1}-t_k| + \sum_{k=k_0 + 1}^{T}D_k\nonumber\\
&\leq f(x_{k_0 + 1}) - f^* + L_0\sqrt{d}\left(t_{k_0 + 1}+\sum_{k=k_0 + 1}^T|t_{k+1}-t_k|\right) + \sum_{k=k_0 + 1}^{T}D_k,
\label{45}
\end{align}
where the last inequality follows from Lemma \ref{lem:Lip_t}.
%\memo{Why the last inequality holds? $F(x_1,t_1)\leq f(x_1) $ holds true? The last inequality follows from \eqref{constantupper}.}
Observe from \eqref{expuxi} that
\begin{align}
    \mathbb{E}_{\zeta_k} [ \left<\nabla_x F(x_k,t_k),\delta_k\right>\mid\zeta_{[k-1]}] = 0\nonumber.
\end{align}
Thus, we have
\begin{align}
    \mathbb{E}_{\zeta_k} [D_k\mid \zeta_{[k-1]}] &= -\beta \mathbb{E}_{\zeta_k}[ \left<\nabla_x F(x_k,t_k),\delta_k\right>\mid \zeta_{[k-1]}] + \frac{L_1}{2}\beta^2\mathbb{E}_{\zeta_k}[\|\Tilde{G}_x(x_k,t_k;\xi_k, u_k)\|^2\mid \zeta_{[k-1]}] \nonumber\\
    &= \frac{L_1}{2}\beta^2\mathbb{E}_{\zeta_k}(\|\Tilde{G}_x(x_k,t_k;\xi_k, u_k)\|^2\mid \zeta_{[k-1]})\nonumber\\
    &\leq \frac{L_1}{2}\beta^2\left(\frac{\mathbb{E}_{\zeta_k}[t_k^2\mid\zeta_{[k-1]}]}{2}L_1^2(d+6)^3+2(d+4) (\mathbb{E}_{\zeta_k}[\|\nabla_x \bar{f}(x_k;\xi_k)\mid\zeta_{[k-1]}\|^2])\right)\nonumber\\
    &\leq \frac{L_1}{2}\beta^2\left(\frac{\mathbb{E}_{\zeta_k}[t_k^2\mid\zeta_{[k-1]}]}{2}L_1^2(d+6)^3+2(d+4) (\mathbb{E}_{\zeta_k}[\|\nabla f(x_k)\mid\zeta_{[k-1]}\|^2]+\sigma^2)\right),\label{47}
\end{align}
where the fist inequality follows from Lemma~\ref{lem:zo_grad_square_diff} and the last inequality holds due to Assumption~\ref{A2}~(ii).

Take the expectation for \eqref{45} with respect to $\zeta_{k_0 + 1}, \ldots, \zeta_{T}$. Together with Lemma \ref{lem:grad_square_diff} (ii), we have
\begin{align}
&\sum_{k=k_0 + 1}^{T}\beta \mathbb{E}_{\zeta}[\|\nabla f(x_k)\|^2] \nonumber \\
&\leq \sum_{k=k_0 + 1}^{T}\beta \mathbb{E}_{\zeta}[\|\nabla_x F(x_k,t_k)\|^2] + \sum_{k=k_0 + 1}^T\beta\mathbb{E}_\zeta[t_k]L_0L_1(d+3)^{3/2} \nonumber\\
&\leq 
f(x_{k_0 + 1}) - f^* + L_0\sqrt{d}\left(t_{k_0 + 1} + \sum_{k=k_0 + 1}^T\mathbb{E}_\zeta[|t_{k+1}-t_k|]\right) + \sum_{k=k_0 + 1}^T \mathbb{E}_\zeta[D_k] \nonumber\\
&+ 
\sum_{k=k_0 + 1}^T\mathbb{E}_\zeta[t_k]L_0L_1\beta(d+3)^{3/2} \nonumber\\
&\leq 
f(x_{k_0 + 1}) - f^* + L_0\sqrt{d}\left(t_{k_0 + 1} + \sum_{k=k_0 + 1}^T\mathbb{E}_\zeta[|t_{k+1}-t_k|]\right)+\sum_{k=k_0 + 1}^T\mathbb{E}_\zeta[t_k]L_0L_1\beta(d+3)^{3/2}\nonumber\\
&+ 
\sum_{k=k_0 + 1}^T \frac{\mathbb{E}_\zeta[t_k^2]}{4}L_1^3\beta^2(d+6)^3+\sum_{k=k_0 + 1}^T L_1\beta^2(d+4)\mathbb{E}_\zeta[\|\nabla f(x_k)\|^2]+L_1\beta^2(d+4)\sigma^2(T-k_0),\nonumber
\end{align}
where the last inequality holds due to \eqref{47}.
Rearrange the terms in the above inequality. Then we can get
\begin{align}
(\beta-(d+4)L_1\beta^2)\sum_{k=k_0 + 1}^{T}\mathbb{E}_{\zeta}[\|\nabla f(x_k)\|^2] &\leq f(x_{k_0 + 1}) - f^* + L_0\sqrt{d}\left(t_{k_0 + 1} + \sum_{k=k_0 + 1}^T\mathbb{E}_\zeta[|t_{k+1}-t_k|]\right)\nonumber\\
&+ \frac{L_1^3\beta^2(d+6)^3}{4}\sum_{k=k_0 + 1}^T\mathbb{E}_\zeta[t_k^2] + L_1\beta^2(d+4)\sigma^2(T - k_0) \nonumber\\
&+\sum_{k=k_0 + 1}^T\mathbb{E}_\zeta[t_k]L_0L_1\beta(d+3)^{3/2},
\label{100}
\end{align}

If the step size $\beta$ is chosen as $\mathop{\rm min} \left\{\frac{1}{2(d+4)L_1}, \frac{1}{\sqrt{(T-k_0)(d+4)}}\right\}$, then we have
\begin{align*}
    \frac{1}{\beta-(d+4)L_1\beta^2} \leq \frac{2}{\beta},\quad \frac{1}{\beta} \leq 2(d+4)L_1 + \sqrt{(T - k_0)(d+4)}.
\end{align*}
Hence,
by dividing both sides of \eqref{100} by $(T - k_0)(\beta-2(d+4)L_1\beta^2)$, we can obtain 
\begin{align}
    \nonumber
&\frac{1}{T - k_0}\sum_{k=1}^{T}\mathbb{E}_{\zeta}[\|\nabla f(x_k)\|^2] \\
\nonumber
&\leq 
\frac{f(x_{k_0 + 1}) - f^* + L_0\sqrt{d}\left(t_{k_0 + 1} + \sum_{k=k_0 + 1}^T\mathbb{E}_\zeta[|t_{k+1}-t_k|]\right) +L_0L_1(d+3)^{3/2}\beta\sum_{k=k_0 + 1}^T\mathbb{E}_\zeta[t_k]}{(T - k_0)(\beta-(d+4)L_1\beta^2)}\\
\nonumber
&+\frac{\frac{L_1^3\beta^2(d+6)^3}{4}\sum_{k=k_0 + 1}^T\mathbb{E}_\zeta[t_k^2]+L_1\beta^2(d+4)\sigma^2T}{{(T - k_0)(\beta-(d+4)L_1\beta^2)}}\\
\nonumber
&\leq \frac{2}{T - k_0}\left(f(x_{k_0 + 1}) - f^* + L_0\sqrt{d}\left(t_{k_0 + 1} + \sum_{k=k_0 + 1}^T\mathbb{E}_\zeta[|t_{k+1}-t_k|]\right)\right)\left(2(d+4)L_1+\sqrt{(T - k_0)(d+4)}\right)\\
\nonumber
&+\frac{2}{T - k_0}L_0L_1(d+3)^{3/2}\sum_{k=k_0 + 1}^{T}\mathbb{E}_{\zeta}[t_k]+\frac{L_1^3\beta(d+6)^3}{2(T - k_0)}\sum_{k=k_0 + 1}^T\mathbb{E}_\zeta[t_k^2]
+2L_1\beta(d+4)\sigma^2\\
\nonumber
&= O\left(\frac{\sqrt{d}\left(1+\sqrt{d}\sum_{k=k_0 + 1}^T\mathbb{E}_\zeta[|t_{k+1}-t_k|]\right)}{\sqrt{T - k_0}}+\frac{d\left(d\mathbb{E}_\zeta\left[ \sum_{k=k_0 + 1}^T t_k^2 \right]+ \sqrt{d}\mathbb{E}_\zeta\left[\sum_{k=k_0 + 1}^Tt_k \right]+1\right)}{T - k_0}\right)
\\
&= O\left(\frac{\sqrt{d}\left(1+\sqrt{d}\gamma^{k_0}\right)}{\sqrt{T - k_0}}+\frac{d\left(d\gamma^{2k_0}+\sqrt{d}\gamma^{k_0}+1\right)}{T - k_0}\right) \nonumber 
\end{align}
where the last equality follows from the update rule of $t_k$,
as shown in the proof of Theorem~\ref{iter_determin} as well.

Here, we have 
$
    k_0 = O \left(
        \frac{d}{\epsilon^4}
    \right)
$
by the definition of $k_0$. Thus, by setting
$
    T
    = k_0 + O \left(
        \frac{d^2}{\epsilon^4}
    \right)
    = 
    O \left(
        \frac{d^2}{\epsilon^4}
    \right)
$,
we can obtain
$
    \mathbb{E}_{\zeta, k'}[\|\nabla f(\hat{x})\|^2]
    =
    \frac{1}{T - k_0}\sum_{k=k_0 + 1}^{T}\mathbb{E}_{\zeta}[\|\nabla f(x_k)\|^2]
    \leq \epsilon^2
$.
% As we have
This implies 
$
    \mathbb{E}_{\zeta, k'}[\|\nabla f(\hat{x})\|] \leq \epsilon
$
as
$
    \mathbb{E}_{\zeta, k'}[\|\nabla f(\hat{x})\|]^2 
    \leq
    \mathbb{E}_{\zeta, k'}[\|\nabla f(\hat{x})\|^2]
$
follows from Jensen's inequality. Furthermore, when $\gamma$ is chosen as
$\gamma \leq d^{-\epsilon^4/2d}$
, we have 
$
\log_\gamma d^{-1/2} = O\left( \frac{d}{\epsilon^4} \right)
$, which implies that 
$
k_0 = \Omega\left( \log_\gamma d^{-1/2} \right)
$. Therefore, we can obtain
$
\gamma^{k_0}
=
O( d^{-1/2} )
$, which yields the iteration complexity of $T=O\left(\frac{d}{\epsilon^4}\right)$.

\hfill $\Box$\\

\section{ZOSLGH algorithm with error tolerance}
\label{sec:error_appendix}

\begin{figure}[H]
\centering
\begin{tabular}{lc}
    %\hspace{-5mm}
    \subfigure[Smooth objective function]{
    \label{fig:smooth}
    \begin{minipage}[t]{0.50\linewidth}
        \centering
        \includegraphics[width=1.5in]{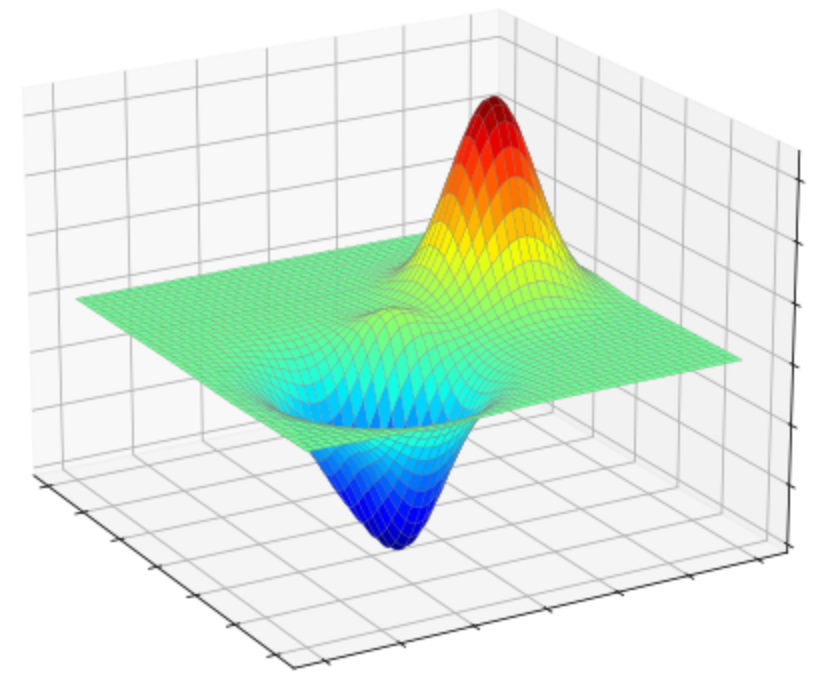}
        \vspace{1mm}
    \end{minipage}%
    }
    \subfigure[Accessible function with error]{
    \label{fig:error}
    \begin{minipage}[t]{0.50\linewidth}
        \centering
        \includegraphics[width=1.5in]{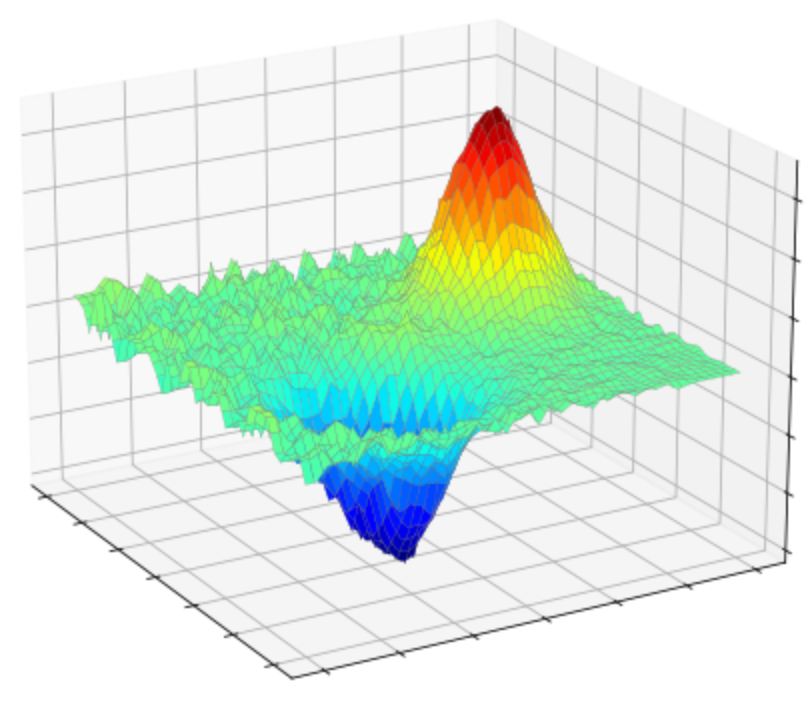}
        \vspace{1mm}
     \end{minipage}
    }
\end{tabular}
\centering
\caption{Illustration of a smooth objective function and the  accessible function that contains error.}
\label{fig:smooth_and_error}
\end{figure}

In Sections \ref{sec:first-order} and \ref{sec:zeroth-order}, we assumed that we had access to the exact function value or a gradient oracle whose variance was finite. However, in some practical cases, we will have access only to the function values containing error, and it would be impossible to obtain accurate gradient oracles of an underlying objective function. Figure \ref{fig:smooth_and_error} illustrates such a case; although the objective function $f$ (Figure \ref{fig:smooth}) is smooth, the accessible function $f'$ (Figure \ref{fig:error}) contains some error, and thus many local minima arise. In this section, we consider optimizing a smooth objective function $f$ using only the information of $f'$. We assume that the following condition holds between $f$ and $f'$.

\begin{assumption}{A3}\label{A3}$\ $
\renewcommand{\labelenumi}{(\roman{enumi})}
The supremum norm of the difference between $f$ and $f'$ is uniformly bounded:
$$\sup_{x\in\mathbb{R}^d}|f(x)-f'(x)|\leq\nu.$$
In the stochastic setting, we assume $\sup_{x\in\mathbb{R}^d}|f(x;\xi)-f'(x;\xi)|\leq\nu$ for any $\xi$.
\end{assumption}

Please note that we do not impose any other assumptions on the accessible function $f'$. Thus, $f'$ can be non-Lipschitz or even discontinuous. Even in such cases, we can develop an algorithm with a convergence guarantee because its smoothed function $F'(x,t)$ is smooth as far as $t$ is sufficiently large. In the following, we denote the Lipschitz and gradient Lipschitz constant of $F'(\cdot,t)$ as $L_0(t)$ and $L_1(t)$, respectively.

The ZOSLGH algorithm in this setting is almost the same as Algorithm \ref{alg:ZOGH}. The only difference is $\sqrt{\nu}$ rather than $\epsilon$ in the update rule of $t_{k+1}$. See the Algorithm \ref{alg:ZOGH_error} for a more detailed description. Please note that $F', \tilde{g}'_x, \tilde{G}'_{x,u}, \tilde{g}'_t, \tilde{G}'_{t,v}$ are defined in the same way as the no-error setting using $f'$.

\begin{algorithm}[H]
\caption{Deterministic/Stochastic Zeroth-Order Single Loop GH algorithm (ZOSLGH) with error tolerance}
\label{alg:ZOGH_error} 
\begin{algorithmic}
\REQUIRE Iteration number $T$, initial solution $x_1$, initial smoothing parameter $t_1$, sequence of step sizes $\{\beta_k\}$ for $x$, step size $\eta$ for $t$, decreasing factor $\gamma \in (0,1)$, error tolerance $\nu$
\FOR{$k=1$ to $T$}
    \STATE Sample $u_k$ from $\mathcal{N}(0,\mathrm{I}_d)$
    \STATE Update $x_{k}$ by
        \begin{align*}
            &x_{k+1} = x_k - \beta_k \bar{G}'_{x,u},\\ 
            &\text{where  } \bar{G}'_{x,u} = \left\{\begin{array}{cc}
                \Tilde{g}'_{x}(x_k,t_k;u_k) & (\text{deterministic}) \\
                \Tilde{G}'_{x}(x_k,t_k;\xi_k,u_k),\ \xi_k\sim P & (\text{stochastic})
            \end{array}\right.
        \end{align*}
    \STATE Sample $v_k$ from $\mathcal{N}(0,\mathrm{I}_d)$
    \STATE Update $t_k$ by
    \begin{align*}
            &t_{k+1} = \left\{\begin{array}{cc}
                \text{max} \{\gamma t_{k}, \sqrt{\nu}\} & (\text{SLGH}_{\text{r}}) \\
                \text{max}\{\text{min} \{t_{k} - \eta \bar{G}'_{t,v},
                \gamma t_{k}\},
                \sqrt{\nu}\} & (\text{SLGH}_{\text{d}})
            \end{array},\right.\\ 
            &\text{where  } \bar{G}'_{t,v} = \left\{\begin{array}{cc}
                \Tilde{g}'_{t}(x_k,t_k;v_k) & (\text{deterministic}) \\
                \Tilde{G}'_{t}(x_k,t_k;\xi_k,v_k),\ \xi_k\sim P & (\text{stochastic})
            \end{array}\right.
    \end{align*}
\ENDFOR
%\ENSURE ~~ $x_T$
\end{algorithmic}
\end{algorithm}

We provide the convergence analyses in the following theorems.
The definitions of $\hat{x}$ in the deterministic and stochastic settings are given in Appendix \ref{subsec:proof_error_deterministic} and \ref{subsec:proof_error_stochastic}, respectively.
\begin{theorem}[\textbf{Convergence of ZOSLGH with error tolerance, Deterministic setting}]
    \label{thm:error_deterministic}
    Suppose Assumptions~\ref{A1} and \ref{A3} hold.
    
    Take $k_1 := \Theta(d/\epsilon^2)$ and $k_2 := O\left(\log_{\gamma} 1/d \right)$ and define $k_0 = \min \{ k_1, k_2 \}$. Let $\hat{x}:=x_{k'}$, where $k'$ is chosen from a uniform distribution over $\{ k_0+1, k_0+2, \ldots, T \}$. Set the stepsize for $x$ at iteration $k$ as $\beta_k=\frac{1}{16(d+4)L_1(t_k)}\ ,k\in[T]$.
    % If the error level $\nu$ satisfies $\nu=O(t_k^2)\ (k\in[T])$, the output sequence of the ZOSLGH algorithm $\{x_k\}_{k=1}^T$ will then satisfy 
    % $
    %     \frac{1}{T}\sum_{k=1}^T \mathbb{E}[\|\nabla f(x_k)\|^2]
    %     = O\left(\frac{d}{T}\left(1+\sqrt{d}\sum_{k=1}^T\mathbb{E}_u[|t_{k+1}-t_k|] + d^2\sum_{k=1}^T\mathbb{E}_u[t_k^2]\right)\right)
    % $,
    % where the expectation is taken w.r.t.~random vectors $\{u_k\}$.
    % Consequently,
    % if we update $t_k$ as in Algorithm~\ref{alg:ZOGH_error} with an arbitrary constant $\gamma \in (0, 1)$,
    Then,
    for any setting of the parameter $\gamma$,
    if the error level $\nu$ satisfies $\nu=O(\epsilon^2/d^3)$,
    $\hat{x}$ satisfies
    $\mathbb{E}[ \|\nabla f(\hat{x})\| ] \leq\epsilon$ with the iteration complexity of 
    $T = O({d^{3}}/{\epsilon^2})$,
    where the expectation is taken w.r.t.~random vectors $\{u_k\}$ and $\{v_k\}$.
    Further,
    if we choose $\gamma \leq d^{-\Omega(\epsilon^2 / d)}$,
    the iteration complexity can be bounded as
    $T = O({d}/{\epsilon^2})$.
\end{theorem}

\begin{theorem}[\textbf{Convergence of ZOSLGH with error tolerance, Stochastic setting}]
    \label{thm:error_stochastic}
    Suppose Assumptions~\ref{A1}, \ref{A2} and \ref{A3} hold. 
    Take $k_1 := \Theta(d/\epsilon^4)$ and $k_2 := O\left(\log_{\gamma} 1/d \right)$ and define $k_0 = \min \{ k_1, k_2 \}$. Let $\hat{x}:=x_{k'}$, where $k'$ is chosen from a uniform distribution over $\{ k_0+1, k_0+2, \ldots, T \}$. Set the stepsize for $x$ at iteration $k$ as $\beta_k=\min\left\{\frac{1}{16(d+4)L_1(t_k)}, \frac{1}{\sqrt{(T - k_0)(d+4)}}\right\}$.
    Then, for any setting of the parameter $\gamma$, if the error level $\nu$ satisfies $\nu=O(\epsilon^2/d^3)$, $\hat{x}$ satisfies $\mathbb{E}[\|\nabla f(\hat{x})\|]\leq\epsilon$ with the iteration complexity of $T=O(d^2/\epsilon^4+d^3/\epsilon^2)$, where the expectation is taken w.r.t. random vectors $\{u_k\}, \{v_k\}$ and $\{\xi_k\}$. Further, if we choose $\gamma \leq d^{-\Omega(\epsilon^4 / d)}$, the iteration complexity can be bounded as $T=O(d/\epsilon^4)$.
\end{theorem}

\subsection{Proofs for technical lemmas}
\label{subsec:tech_lemmas}
We introduce several lemmas before going to the convergence analysis. All of them describe properties of the function with error $f'$ and its Gaussian smoothing $F'$. Throughout this subsection, we assume that $f$ is $L_0$-Lipschitz and $L_1$-smooth function. We also suppose that the function pair $(f, f')$ satisfies $\mathop{\rm sup}\limits_{x\in\mathbb{R}^d}|f(x)-f'(x)|\leq\nu$.

\begin{lemma}
    For any $x\in \mathbb{R}^d$ and $t>0$, we have
    \begin{align*}
        \mathbb{E}_u\left[\frac{1}{t^2}(f'(x+tu)-f'(x))^2\|u\|^2\right] \leq 4(d+4)\|\nabla f(x)\|^2 + t^2L_1^2(d+6)^3 +  8d\frac{\nu^2}{t^2}.
    \end{align*}
\label{lem_error:zo_grad_square_diff}
\end{lemma}
\textbf{Proof:} 
\begin{align*}
    \mathbb{E}_u\left[\frac{1}{t^2}(f'(x+tu)-f'(x))^2\|u\|^2\right] &= \mathbb{E}_u\left[\frac{1}{t^2}(f(x+tu)-f(x) + (f'-f)(x+tu)-(f'-f)(x))^2\|u\|^2\right]\\
    &\leq 2\mathbb{E}_u\left[\frac{1}{t^2}(f(x+tu)-f(x))^2\|u\|^2\right] + 2\mathbb{E}_u\left[\frac{1}{t^2}(2\nu)^2\|u\|^2\right]\\
    &\leq 4(d+4)\|\nabla f(x)\|^2 + t^2L_1^2(d+6)^3 +  8d\frac{\nu^2}{t^2},
\end{align*}
where the last inequality holds due to Lemma~\ref{lem:gauss_norm} and Lemma~\ref{lem:zo_grad_square_diff}.

\begin{lemma}
    For any $x\in \mathbb{R}^d$ and $t>0$, we have
    \begin{align*}
        \mathbb{E}_\zeta\left[\frac{1}{t^2}(\bar{f}'(x+tu;\xi)-\bar{f}'(x;\xi))^2\|u\|^2\right] \leq 4(d+4)(\|\nabla f(x)\|^2+\sigma^2) + t^2L_1^2(d+6)^3 + 8d\frac{\nu^2}{t^2}.
    \end{align*}
\label{lem_error:zo_stochastic_grad_square_diff}
\end{lemma}
\textbf{Proof:} 
\begin{align*}
    &\mathbb{E}_\zeta\left[\frac{1}{t^2}(\bar{f}'(x+tu;\xi)-\bar{f}'(x;\xi))^2\|u\|^2\right]\\ 
    &= \mathbb{E}_\xi\left[\mathbb{E}_u\left[\frac{1}{t^2}(\bar{f}(x+tu;\xi)-\bar{f}(x;\xi) + (\bar{f}'-\bar{f})(x+tu;\xi)-(\bar{f}'-\bar{f})(x;\xi))^2\|u\|^2\right]\right]\\
    &\leq 2 \mathbb{E}_\xi\left[\mathbb{E}_u\left[\frac{1}{t^2}(\bar{f}(x+tu;\xi)-\bar{f}(x;\xi))^2\|u\|^2\right]\right] + \frac{2}{t^2}\mathbb{E}_\xi[\mathbb{E}_u[(2\nu)^2\|u\|^2]]\\
    &\leq 2\mathbb{E}_\xi\left[\frac{t^2}{2}L_1^2(d+6)^3+2(d+4)\|\nabla \bar{f}(x;\xi)\|^2\right] + 8d\frac{\nu^2}{t^2}\\
    &\leq 4(d+4)(\|\nabla f(x)\|^2+\sigma^2) + t^2L_1^2(d+6)^3 + 8d\frac{\nu^2}{t^2},
\end{align*}
where the second inequality follows from Lemma~\ref{lem:gauss_norm} and Lemma~\ref{lem:zo_grad_square_diff},
 and the last inequality holds due to Assumption~\ref{A2}~(ii).
\begin{lemma}
\label{lem_error:Lip_t}
For any $x\in \mathbb{R}^d$ and for any $t_1,t_2\in\mathcal{T}$, we have
\begin{align*}
    |F'(x,t_1)-F'(x,t_2)|\leq L_0|t_1-t_2|\sqrt{d}+2\nu.
    \end{align*}
\end{lemma}
\textbf{Proof:} 
\begin{align*}
    |F'(x,t_1)-F'(x,t_2)| &= |F(x,t_1)-F(x,t_2) + (F'-F)(x,t_1)-(F'-F)(x,t_2)|\\
    &\leq |F(x,t_1)-F(x,t_2)| + |\mathbb{E}_u[(f'-f)(x+t_1u)]|+|\mathbb{E}_u[(f'-f)(x+t_2u)]|\\
    &\leq |F(x,t_1)-F(x,t_2)| + \mathbb{E}_u[|(f'-f)(x+t_1u)|]+\mathbb{E}_u[|(f'-f)(x+t_2u)|]\\
    &\leq |F(x,t_1)-F(x,t_2)|+ 2\nu\\
    &\leq L_0|t_1-t_2|\sqrt{d}+2\nu,
\end{align*}
where the last inequality holds due to Lemma \ref{lem:Lip_t}.

\begin{lemma}[\textbf{Lemma 30 in \cit{jin2018local}}]\label{lem:grad_diff} For any $x\in \mathbb{R}^d$ and for any $t_1,t_2\in\mathcal{T}$, we have
$$\|\nabla_x (F'-F)(x,t)\|\leq\sqrt{\frac{2}{\pi}}\frac{\nu}{t}.$$
\label{lem_error:grad_diff}
\end{lemma}

\begin{lemma}\label{lem:F'_lip}\quad\\
\textbf{(i)} $F'(x,t)$ is $L_0+\sqrt{\frac{2}{\pi}}\frac{\nu}{t}$-Lipschitz in terms of $x$.\\
\textbf{(ii) (\textbf{Lemma 20 in \cit{jin2018local}})} $F'(x,t)$ is $L_1+\frac{2\nu}{t^2}$-smooth in terms of $x$.
\label{lem_error:Lip}
\end{lemma}
\textbf{Proof for (i):}
\begin{align*}
    \|\nabla_x F'(x,t)\|&\leq \|\nabla_x F(x,t)\| + \|\nabla_x (F'-F)(x,t)\|\\
    &\leq L_0 + \sqrt{\frac{2}{\pi}}\frac{\nu}{t},\\
\end{align*}
where the last inequality holds due to Lemma \ref{lem:Lip} and Lemma \ref{lem:grad_diff}.

\begin{lemma}\label{lem_error:grad_square_diff} For any $x\in \mathbb{R}^d$ and $t>0$, we have
    \begin{align*}
        \|\nabla f(x)\|^2 \leq 4\|\nabla_x F'(x,t)\|^2 + \frac{t^2}{2}L_1^2(d+6)^3 + \frac{8}{\pi}\frac{\nu^2}{t^2}.
    \end{align*}
\end{lemma}
\textbf{Proof:}
We have
\begin{align*}
    \|\nabla f(x)\|^2 &= \|\mathbb{E}_u[\langle\nabla f(x), u\rangle u]\|^2\\
    &= \left\|\frac{1}{t}\mathbb{E}_u[(f(x+tu)-f(x) - [f(x+tu)-f(x)-t\langle \nabla f(x), u\rangle])u]\right\|^2\\
    &\leq \left\|\nabla_x F(x,t) -\frac{1}{t}\mathbb{E}_u[(f(x+tu)-f(x)-t\langle \nabla f(x), u\rangle)u]\right\|^2\\
    &\leq 2\|\nabla_x F(x,t)\|^2 + \frac{2}{t^2}\left\|\mathbb{E}_u[(f(x+tu)-f(x)-t\langle \nabla f(x), u\rangle)u]\right\|^2\\
    &\leq 2\|\nabla_x F(x,t)\|^2 + \frac{2}{t^2}\mathbb{E}_u[|f(x+tu)-f(x)-t\langle \nabla f(x), u\rangle|^2\|u\|^2]\\
    &\leq 2\|\nabla_x F(x,t)\|^2 + \frac{t^2L_1^2}{2}\mathbb{E}_u[\|u\|^6]\\
    &\leq 2\|\nabla_x F(x,t)\|^2 + \frac{t^2L_1^2}{2}(d+6)^3\\
    &\leq 2(2\|\nabla_x (F-F')(x,t)\|^2 + 2\|\nabla_x F'(x,t)\|^2) + \frac{t^2L_1^2}{2}(d+6)^3\\
    &\leq 4 \|\nabla_x F'(x,t)\|^2 + \frac{t^2L_1^2}{2}(d+6)^3 + \frac{8}{\pi}\frac{\nu^2}{t^2},
\end{align*}
where the third last inequality holds due to Lemma \ref{lem:gauss_norm}, and the last inequality holds due to Lemma \ref{lem:grad_diff}.

\subsection{Proof for the deterministic setting}
\label{subsec:proof_error_deterministic}

\textbf{Proof for Theorem \ref{thm:error_deterministic}:}
Let $w_k:=(u_k, v_k)$ and denote $\delta_k := \Tilde{g}'_x(x_k,t_k;u_k)-\nabla_x F'(x_k,t_k)$. Utilize the updating rule for $x$ and $L_1(t)$-smoothness of $F'(\cdot,t)$. Then we have
\begin{align}
F'(x_{k+1},t_k) &\leq F'(x_k,t_k) + \left<\nabla_x F'(x_k,t_k),(x_{k+1}-x_k)\right> + \frac{L_1(t_k)}{2}\|x_{k+1}-x_k\|^2\nonumber\\
&= F'(x_k,t_k) - \beta_k \left<\nabla_x F'(x_k,t_k),\Tilde{g}'_x(x_k,t_k;u_k)\right> + \frac{L_1(t_k)}{2}\beta_k^2\|\Tilde{g}'_x(x_k,t_k;u_k)\|^2\nonumber\\
&= F'(x_k,t_k) - \beta_k\|\nabla_x F'(x_k,t_k)\|^2 - \beta_k \left<\nabla_x F'(x_k,t_k),\delta_k\right> + \frac{L_1(t_k)}{2}\beta_k^2\|\Tilde{g}'_x(x_k,t_k;u_k))\|^2.
\label{error_1}
\end{align}
Denote
\[E_k := - \beta_k \left<\nabla_x F'(x_k,t_k),\delta_k\right> + \frac{L_1(t_k)}{2}\beta_k^2\|\Tilde{g}'_x(x_k,t_k;u_k)\|^2
  \]
  for simplicity. From Lemma \ref{lem_error:Lip_t} and \eqref{error_1}, we get the upper bound for $\|\nabla_x F'(x,t)\|^2$ as 
\begin{align}
\beta_k \|\nabla_x F'(x_k,t_k)\|^2 &\leq F'(x_k,t_k) - F'(x_{k+1},t_k) + E_k\nonumber\\
&= F'(x_k,t_k) - F'(x_{k+1},t_{k+1}) + F'(x_{k+1},t_{k+1}) - F'(x_{k+1},t_k) + E_k\nonumber\\
&\leq F'(x_k,t_k) - F'(x_{k+1},t_{k+1}) + L_0|t_{k+1}-t_k|\sqrt{d} + 2\nu + E_k\nonumber.
\end{align}
Now, sum up the above inequality for all iterations $k_0+1\leq k\leq T\ (T>k_0)$. Then we have
\begin{align}
&\sum_{k=k_0 + 1}^{T}\beta_k\|\nabla_x F'(x_k,t_k)\|^2\nonumber\\
&\leq F'(x_{k_0 + 1},t_{k_0 + 1}) - F'(x_{T+1},t_{T+1}) + L_0\sum_{k=k_0 + 1}^T |t_{k+1}-t_k|\sqrt{d} + 2\nu (T - k_0) + \sum_{k=k_0 + 1}^{T}E_k\nonumber\\
&\leq F'(x_{k_0 + 1},t_{k_0 + 1}) - f^* +\nu + L_0\sqrt{d}\sum_{k=k_0 + 1}^T |t_{k+1}-t_k| + 2\nu (T - k_0) + \sum_{k=k_0 + 1}^{T}E_k\nonumber.\\
&\leq f'(x_{k_0 + 1}) - f^* +3\nu + L_0\sqrt{d}\left(t_{k_0 + 1}+\sum_{k=k_0 + 1}^T |t_{k+1}-t_k|\right) + 2\nu (T - k_0) + \sum_{k=k_0 + 1}^{T}E_k\nonumber\\
&\leq f(x_{k_0 + 1}) - f^* +4\nu + L_0\sqrt{d}\left(t_{k_0 + 1}+\sum_{k=k_0 + 1}^T |t_{k+1}-t_k|\right) + 2\nu (T - k_0) + \sum_{k=k_0 + 1}^{T}E_k,
\label{error_9}
\end{align}
where the third inequality holds due to Lemma~\ref{lem_error:Lip_t}. We can bound the conditional expectation of $E_k$ as
\begin{align*}
    &\mathbb{E}_{w_k} [E_k\mid w_{[k-1]}] \\ 
    &= -\beta_k \mathbb{E}_{w_k}[ \left<\nabla_x F'(x_k,t_k),\delta_k\right>\mid w_{[k-1]}] + \frac{\mathbb{E}_{w_k} [L_1(t_k)\mid w_{[k-1]}]}{2}\beta_k^2\mathbb{E}_{w_k}[\|\Tilde{g}'_x(x_k,t_k;u_k)\|^2\mid w_{[k-1]}]\\
    &\leq \frac{\mathbb{E}_{w_k} [L_1(t_k)\mid w_{[k-1]}]}{2}\beta_k^2\mathbb{E}_{w_k}[\|\Tilde{g}'_x(x_k,t_k;u_k)\|^2\mid w_{[k-1]}]\\
    &\leq \frac{\mathbb{E}_{w_k} [L_1(t_k)\mid w_{[k-1]}]}{2}\beta_k^2\left(4(d+4)\mathbb{E}_{w_k} [\|\nabla f(x_k)\|^2\mid w_{[k-1]}]+L_1^2(d+6)^3\mathbb{E}_{w_k}[t_k^2\mid w_{[k-1]}]\right.\\
    &\left.\hspace{40mm} 
    + 8d\mathbb{E}_{w_k}\left[\nu^2/t_k^2\mid w_{[k-1]}\right]\right),
\end{align*} % TODO: may cause losing form
where the first inequality holds since we have $\E_{w_k}[\delta_k\mid w_{[k-1]}]=\E_{u_k}[\delta_k\mid u_{[k-1]}]=0$, and the last inequality holds due to Lemma~\ref{lem_error:zo_grad_square_diff}. Take the expectations of \eqref{error_9} w.r.t. random vectors $\{w_{k_0+1},...,w_T\}$. Then we can get
\begin{align}
    & \sum_{k=k_0 + 1}^{T}\beta_k \mathbb{E}_{w}[\|\nabla_x F'(x_k,t_k)\|^2] \nonumber\\ 
    & \leq f(x_{k_0 + 1}) - f^* +4\nu + L_0\sqrt{d}\left(t_{k_0 + 1}+\sum_{k=k_0 + 1}^T \mathbb{E}_w[|t_{k+1}-t_k|]\right) + 2\nu (T - k_0) \nonumber \\ 
    & + \frac{1}{2}\left(4(d+4)\sum_{k=k_0 + 1}^{T}\beta_k^2\mathbb{E}_{w} [L_1(t_k)\|\nabla f(x_k)\|^2]+L_1^2(d+6)^3\sum_{k=k_0 + 1}^{T}\beta_k^2\mathbb{E}_{w}[L_1(t_k)t_k^2]\right.\nonumber\\
    &\left.\hspace{8.5mm}
     + 8d\sum_{k=k_0 + 1}^{T}\beta_k^2\mathbb{E}_{w}\left[L_1(t_k)\frac{\nu^2}{t_k^2}\right]\right).
    \label{300}
\end{align}
Lemma~\ref{lem_error:grad_square_diff} together with \eqref{300} yields 
\begin{align*}
    & \sum_{k=k_0 + 1}^{T}\beta_k \mathbb{E}_{w}[\|\nabla f(x_k)\|^2] \\
    & \leq 4\sum_{k=k_0 + 1}^{T}\beta_k \mathbb{E}_{w}[\|\nabla_x F'(x_k,t_k)\|^2] + \frac{1}{2}\sum_{k=k_0 + 1}^{T}\beta_k\mathbb{E}_{w}[t_k^2]L_1^2(d+6)^3 + \frac{8}{\pi}\sum_{k=k_0 + 1}^{T}\beta_k\mathbb{E}_{w}\left[\frac{\nu^2}{t_k^2}\right]\\
    &\leq 4\left(f(x_{k_0 + 1}) - f^* +4\nu + L_0\sqrt{d}\left(t_{k_0 + 1}+\sum_{k=k_0 + 1}^T \mathbb{E}_w[|t_{k+1}-t_k|]\right) + 2\nu (T - k_0)\right)\\
    &+2\left(4(d+4)\sum_{k=k_0 + 1}^{T}\beta_k^2\mathbb{E}_{w} [L_1(t_k)\|\nabla f(x_k)\|^2]+L_1^2(d+6)^3\sum_{k=k_0 + 1}^{T}\beta_k^2\mathbb{E}_{w}[L_1(t_k)t_k^2]\right.\\
    &\left.\hspace{8.5mm}
    + 8d\sum_{k=k_0 + 1}^{T}\beta_k^2\mathbb{E}_{w}\left[L_1(t_k)\frac{\nu^2}{t_k^2}\right]\right)\\
    &+ \frac{1}{2}\sum_{k=k_0 + 1}^{T}\beta_k\mathbb{E}_{w}[t_k^2]L_1^2(d+6)^3 + \frac{8}{\pi}\sum_{k=k_0 + 1}^{T}\beta_k\mathbb{E}_{w}\left[\frac{\nu^2}{t_k^2}\right].
\end{align*}

By rearranging the terms, we obtain
\begin{align}
&\sum_{k=k_0 + 1}^{T}\left(\beta_k\mathbb{E}_{w}[\|\nabla f(x_k)\|^2]-8(d+4)\beta_k^2\mathbb{E}_{w}[L_1(t_k)\|\nabla f(x_k)\|^2]\right)\nonumber\\
&\leq 4\left(f(x_{k_0 + 1}) - f^* +4\nu + L_0\sqrt{d}\left(t_{k_0 + 1}+\sum_{k=k_0 + 1}^T \mathbb{E}_w[|t_{k+1}-t_k|]\right) + 2\nu (T - k_0)\right)\nonumber\\
&+ 2\left(L_1^2(d+6)^3\sum_{k=k_0 + 1}^T\beta_k^2\mathbb{E}_w[L_1(t_k)t_k^2] + 8d\sum_{k=k_0 + 1}^T\beta_k^2\mathbb{E}_w\left[L_1(t_k)\frac{\nu^2}{t_k^2}\right]\right)\nonumber\\
&+ \frac{1}{2}\sum_{k=k_0 + 1}^{T}\beta_k\mathbb{E}_{w}[t_k^2]L_1^2(d+6)^3 + \frac{8}{\pi}\sum_{k=k_0 + 1}^{T}\beta_k\mathbb{E}_{w}\left[\frac{\nu^2}{t_k^2}\right].
\label{error_5}
\end{align}

If we update $t_k\ (k\in[T])$ as in Algorithm~\ref{alg:ZOGH_error}, we have $\nu=O(t_k^2)$, which yields $L_1(t_k)=O(1)$ from Lemma \ref{lem_error:Lip}. Hence, by setting the step size $\beta_k$ as $\frac{1}{16(d+4)L_1(t_k)}\ (k\in[T])$, we can obtain 

\begin{align*}
    \frac{1}{T - k_0}\sum_{k=k_0 + 1}^T \mathbb{E}_{w}[\|\nabla f(x_k)\|^2]
    &= O\left(\frac{d}{T - k_0}\left(1+\sqrt{d}\sum_{k=k_0 + 1}^T\mathbb{E}_w[|t_{k+1}-t_k|] + d^2\sum_{k=k_0 + 1}^T\mathbb{E}_w[t_k^2]\right)\right)
\end{align*}

in the same way as before. We can also get
$\sum_{k=k_0 + 1}^T |t_{k+1} - t_k| = \sum_{k=k_0 + 1}^T ( t_k - t_{k+1}) = t_{k_0 + 1} - t_{T+1} = t_{k_0 + 1} = O(\gamma^{k_0}) $.
Further,
we have
\begin{align*}
    \sum_{k=k_0 + 1}^T t_k^2
    \leq
    \sum_{k=k_0 + 1}^T \max \{  t_1^2 \gamma^{2(k-1)}, \nu \}
    \leq
    \sum_{k=k_0 + 1}^T 
    \left(
        t_1^2 \gamma^{2(k-1)} +
        \nu 
    \right)
    &\leq
    \frac{t_1^2 \gamma^{2k_0}}{1 - \gamma^2}
    +
    \nu (T - k_0) \\
    &=
    O( \gamma^{2k_0} + \nu (T - k_0) ),
\end{align*}
where the first inequality follows from the update rule of $t_k$ in Algorithm~\ref{alg:ZOGH_error}.
Hence,
we obtain
\begin{align}
    \frac{1}{T - k_0}\sum_{k=k_0 + 1}^T \mathbb{E}_{w}[\|\nabla f(x_k)\|^2]
    &=
    O\left(
        \frac{d}{T - k_0}
        \left( 1 + \sqrt{d}\gamma^{k_0} + d^2 ( \gamma^{2k_0} + \nu (T - k_0)) \right)
    \right) \nonumber \\
    &=
    O\left(
        \frac{d(1+d^2\gamma^{2k_0})}{T - k_0}
        +
        d^3 \nu
    \right)
    =
    O\left(
        \frac{d(1 + d^2\gamma^{2k_0})}{T - k_0}
        +
        \epsilon^2
    \right) \nonumber,
\end{align}
where the last equality follows from the assumption of $\nu = O(\epsilon^2 / d^3)$.

Here, we have
$
    k_0 = O \left(
        \frac{d}{\epsilon^2}
    \right)
$
by the definition of $k_0$. Thus, by setting
$
    T
    = k_0 + O \left(
        \frac{d^3}{\epsilon^2}
    \right)
    = 
    O \left(
        \frac{d^3}{\epsilon^2}
    \right)
$,
we can obtain
$
    \mathbb{E}_{w, k'}[\|\nabla f(\hat{x})\|^2]
    =
    \frac{1}{T - k_0}\sum_{k=k_0 + 1}^{T}\mathbb{E}_{w}[\|\nabla f(x_k)\|^2]
    \leq \epsilon^2
$.
% As we have
This implies 
$
    \mathbb{E}_{w, k'}[\|\nabla f(\hat{x})\|] \leq \epsilon
$
as
$
    \mathbb{E}_{w, k'}[\|\nabla f(\hat{x})\|]^2 
    \leq
    \mathbb{E}_{w, k'}[\|\nabla f(\hat{x})\|^2]
$
follows from Jensen's inequality. Furthermore, when $\gamma$ is chosen as
$\gamma \leq d^{-\epsilon^2/d}$
, we have 
$
    \log_\gamma d^{-1} = O\left( \frac{d}{\epsilon^2} \right)
$,
which implies 
$
 k_0 = \Omega \left( \log_\gamma d^{-1} \right)
$. Therefore, we can obtain
$
\gamma^{k_0}
=
O( d^{-1} )
$,
which yields the iteration complexity of $T=O\left(\frac{d}{\epsilon^2}\right)$.

\hfill $\Box$

\subsection{Proof for the stochastic setting}
\label{subsec:proof_error_stochastic}

\textbf{Proof for Theorem \ref{thm:error_stochastic}:}

Let $\zeta_k := (\xi_k,u_k, v_k)$, $k\in [T]$ and denote $\delta_k := \Tilde{G}'_x(x_k,t_k;\xi_k,u_k)-\nabla_x F'(x_k,t_k).$
Since $\Tilde{G}'_x(x,t;\xi, u)$ is an unbiased estimator of $\nabla_x F'(x,t)$, we have
\begin{align}
F'(x_{k+1},t_k) &\leq F'(x_k,t_k) + \left<\nabla_x F'(x_k,t_k),(x_{k+1}-x_k)\right> + \frac{L_1(t_k)}{2}\|x_{k+1}-x_k\|^2\nonumber\\
&= F'(x_k,t_k) - \beta_k \left<\nabla_x F'(x_k,t_k),\Tilde{G}'_x(x_k,t_k;\xi_k,u_k)\right> + \frac{L_1(t_k)}{2}\beta_k^2\|\Tilde{G}'_x(x_k,t_k;\xi_k,u_k)\|^2\nonumber\\
&= F'(x_k,t_k) - \beta_k\|\nabla_x F'(x_k,t_k)\|^2 - \beta_k \left<\nabla_x F'(x_k,t_k),\delta_k\right> + \frac{L_1(t_k)}{2}\beta_k^2\|\Tilde{G}'_x(x_k,t_k;\xi_k,u_k)\|^2\nonumber.
\end{align}
Now, denote
\[
I_k := - \beta_k  \left<\nabla_x F'(x_k,t_k),\delta_k\right> + \frac{L_1(t_k)}{2}\beta_k^2\|\Tilde{G}'_x(x_k,t_k;\xi_k,u_k)\|^2
\]
for simplicity. Then, we can get the upper bound for $\|\nabla_x F(x,t)\|^2$ with $I_k$:
\begin{align}
\beta_k \|\nabla_x F'(x_k,t_k)\|^2 &\leq F'(x_k,t_k) - F'(x_{k+1},t_k) + I_k\nonumber\\
&= F'(x_k,t_k) - F'(x_{k+1},t_{k+1}) + F'(x_{k+1},t_{k+1}) - F'(x_{k+1},t_k) + I_k\nonumber\\
&\leq F'(x_k,t_k) - F'(x_{k+1},t_{k+1}) + L_0|t_{k+1}-t_k|\sqrt{d} + 2\nu + I_k,\nonumber
\end{align}
where the last inequality follows from Lemma~\ref{lem_error:Lip_t}. Sum up the above inequality for all iterations $k_0 + 1\leq k\leq T$. Then we have
\begin{align}
&\sum_{k=k_0 + 1}^{T}\beta_k\|\nabla_x F'(x_k,t_k)\|^2\nonumber \\
&\leq F'(x_{k_0 + 1},t_{k_0 + 1}) - F'(x_{T+1},t_{T+1}) + L_0\sqrt{d}\sum_{k=k_0 + 1}^T|t_{k+1}-t_k| + 2\nu (T - k_0) + \sum_{k=k_0 + 1}^{T}I_k\nonumber\\
&\leq f(x_{k_0 + 1}) - f^* + 4\nu + L_0\sqrt{d}\left(t_{k_0 + 1} +\sum_{k=k_0 + 1}^T|t_{k+1}-t_k|\right) + 2\nu (T - k_0) + \sum_{k=k_0 + 1}^{T}I_k.
\label{error_6}
\end{align}
%\memo{Why the last inequality holds? $F(x_1,t_1)\leq f(x_1) $ holds true? The last inequality follows from \eqref{constantupper}.}
We can also obtain
\begin{align*}
    & \mathbb{E}_{\zeta_k} [I_k\mid \zeta_{[k-1]}]\\
    & = -\beta_k \mathbb{E}_{\zeta_k}[ \left<\nabla_x F'(x_k,t_k),\delta_k\right>\mid \zeta_{[k-1]}] + \frac{\mathbb{E}_{\zeta_k}[L_1(t_k)\mid \zeta_{[k-1]}]}{2}\beta_k^2\mathbb{E}_{\zeta_k}[\|\Tilde{G}'_x(x_k,t_k;\xi_k, u_k)\|^2\mid \zeta_{[k-1]}] \\
    &= \frac{\mathbb{E}_{\zeta_k}[L_1(t_k)\mid \zeta_{[k-1]}]}{2}\beta_k^2\mathbb{E}_{\zeta_k}[\|\Tilde{G}'_x(x_k,t_k;\xi_k, u_k)\|^2\mid \zeta_{[k-1]}]\\
    &\leq \frac{\mathbb{E}_{\zeta_k}[L_1(t_k)\mid \zeta_{[k-1]}]s}{2}\beta_k^2\left(4(d+4)(\mathbb{E}_{\zeta_k}[\|\nabla f(x_k)\|^2\mid \zeta_{[k-1]}]+\sigma^2) + \mathbb{E}_{\zeta_k}[t_k^2\mid \zeta_{[k-1]}]L_1^2(d+6)^3\right. \\
    &\left.
    \hspace{39.8mm} + 8d\mathbb{E}_{\zeta_k}\left[\nu^2/t_k^2\mid \zeta_{[k-1]}\right]\right),  %\label{error_7}
\end{align*} % TODO: may cause losing form
where the last inequality holds due to Lemma~\ref{lem_error:zo_stochastic_grad_square_diff}.

Take the expectation of \eqref{error_6} with respect to $\zeta_{k_0 + 1}, \ldots, \zeta_{T}$. Then we have
\begin{align*}
& \sum_{k=k_0 + 1}^{T}\beta_k \mathbb{E}_{\zeta}[\|\nabla_x F'(x_k,t_k)\|^2]  \\
& \leq f(x_{k_0 + 1}) - f^* + 4\nu + L_0\sqrt{d}\left(t_{k_0 + 1} + \sum_{k=k_0 + 1}^T\mathbb{E}_\zeta[|t_{k+1}-t_k|]\right) + 2\nu (T - k_0) + \sum_{k=k_0 + 1}^T \mathbb{E}_\zeta[I_k] \\
&\leq f(x_{k_0 + 1}) - f^* + 4\nu + L_0\sqrt{d}\left(t_{k_0 + 1} + \sum_{k=k_0 + 1}^T\mathbb{E}_\zeta[|t_{k+1}-t_k|]\right) + 2\nu (T - k_0) +\\
&+  \frac{1}{2}\left(4(d+4)\sum_{k=k_0 + 1}^T\beta_k^2(\mathbb{E}_{\zeta}[L_1(t_k)\|\nabla f(x_k)\|^2]+\sigma^2) + L_1^2(d+6)^3\sum_{k=k_0 + 1}^T\beta_k^2\mathbb{E}_{\zeta}[L_1(t_k)t_k^2] \right.\\
& \left. 
\hspace{8.5mm} + 8d\sum_{k=k_0 + 1}^T\beta_k^2\mathbb{E}_\zeta\left[L_1(t_k)\frac{\nu^2}{t_k^2}\right]\right),
\end{align*}
From Lemma~\ref{lem_error:grad_square_diff}~(ii), we have
\begin{align}
& \sum_{k=k_0 + 1}^{T}\beta_k \mathbb{E}_{\zeta}[\|\nabla f(x_k)\|^2]  \nonumber \\ 
&\leq 4\sum_{k=k_0 + 1}^{T}\beta_k \mathbb{E}_{\zeta}[\|\nabla_x F'(x_k,t_k)\|^2] + \frac{L_1^2(d+6)^3}{2}\sum_{k=k_0 + 1}^{T}\beta_k \mathbb{E}_{\zeta}[t_k^2] +  \frac{8}{\pi}\sum_{k=k_0 + 1}^{T}\beta_k \mathbb{E}_{\zeta}\left[\frac{\nu^2}{t_k^2}\right]\nonumber\\
&\leq 4\left(f(x_{k_0 + 1}) - f^* + 4\nu + L_0\sqrt{d}\left(t_{k_0 + 1} + \sum_{k=k_0 + 1}^T\mathbb{E}_\zeta[|t_{k+1}-t_k|]\right) + 2\nu (T - k_0)\right)\nonumber\\
&+ 2\left(4(d+4)\sum_{k=k_0 + 1}^T(\beta_k^2\mathbb{E}_{\zeta}[L_1(t_k)(\|\nabla f(x_k)\|^2+\sigma^2)]) + L_1^2(d+6)^3\sum_{k=k_0 + 1}^T\mathbb{E}_{\zeta}\beta_k^2[L_1(t_k)t_k^2]\right.\nonumber\\
& \left. 
\hspace{7.5mm} + 8d\sum_{k=k_0 + 1}^T\beta_k^2\mathbb{E}_\zeta\left[L_1(t_k)\frac{\nu^2}{t_k^2}\right]\right) \nonumber \\
&+\frac{L_1^2(d+6)^3}{2}\sum_{k=k_0 + 1}^{T}\beta_k \mathbb{E}_{\zeta}[t_k^2] + \frac{8}{\pi}\sum_{k=k_0 + 1}^{T}\beta_k \mathbb{E}_{\zeta}\left[\frac{\nu^2}{t_k^2}\right].
\label{error_8}
\end{align}
By rearranging the terms, we obtain

\begin{align}
&\sum_{k=k_0 + 1}^{T}\left(\beta_k\mathbb{E}_{\zeta}[\|\nabla f(x_k)\|^2]-8(d+4)\beta_k^2\mathbb{E}_{\zeta}[L_1(t_k)\|\nabla f(x_k)\|^2]\right)\nonumber\\
&\leq 4\left(f(x_{k_0 + 1}) - f^* +4\nu + L_0\sqrt{d}\left(t_{k_0 + 1}+\sum_{k=k_0 + 1}^T \mathbb{E}_\zeta[|t_{k+1}-t_k|]\right) + 2\nu (T - k_0) \right)\nonumber\\
&+ 2\left(4(d+4)\sigma^2\sum_{k=k_0 + 1}^T\beta_k^2\mathbb{E}_\zeta[L_1(t_k)]+L_1^2(d+6)^3\sum_{k=k_0 + 1}^T\beta_k\mathbb{E}_\zeta[L_1(t_k)t_k^2] + 8d\sum_{k=k_0 + 1}^T\beta_k\mathbb{E}_\zeta\left[L_1(t_k)\frac{\nu^2}{t_k^2}\right]\right)\nonumber\\
&+ \frac{L_1^2(d+6)^3}{2}\sum_{k=k_0 + 1}^{T}\beta_k\mathbb{E}_{\zeta}[t_k^2] + \frac{8}{\pi}\sum_{k=k_0 + 1}^{T}\beta_k\mathbb{E}_{\zeta}\left[\frac{\nu^2}{t_k^2}\right].
\end{align}

If we update $t_k\ (k\in[T])$ as in Algorithm~\ref{alg:ZOGH_error}, we have $\nu=O(t_k^2)$, which yields $L_1(t_k)=O(1)$ from Lemma \ref{lem_error:Lip}. Furthermore, if we set the step size $\beta_k$ as $\min\left\{\frac{1}{16(d+4)L_1(t_k)}, \frac{1}{\sqrt{(T - k_0)(d+4)}}\right\}\ (k\in[T])$, then we have
\begin{align*}
    \frac{1}{\beta_k-8(d+4)L_1(t_k)\beta_k^2} \leq \frac{2}{\beta_k},
\end{align*}
\begin{align*}
    \frac{1}{\beta_k} \leq 16(d+4)L_1(t_k) + \sqrt{(T - k_0)(d+4)}.
\end{align*}
for all $k\in[T]$. Using the above inequalities, we can obtain 
\begin{align*}
    \frac{1}{T - k_0}\sum_{k=k_0 + 1}^T \mathbb{E}_{\zeta}[\|\nabla f(x_k)\|^2]
    &= O\left(\frac{\sqrt{d}\left(1+\sqrt{d}\sum_{k=k_0 + 1}^T\mathbb{E}_\zeta[|t_{k+1}-t_k|]\right)}{\sqrt{T - k_0}}+\frac{d\left(1+d^2\sum_{k=k_0 + 1}^T \mathbb{E}_\zeta[t_k^2]\right)}{T - k_0}\right)
    \\
    &
    = O\left(
        \frac{\sqrt{d}+d\gamma^{k_0}}{\sqrt{T - k_0}}
        +
        \frac{d+d^3\gamma^{2k_0}}{T - k_0}
        +
        d^3\nu
    \right)
    \\
    &
    = O\left(
        \frac{\sqrt{d}+d\gamma^{k_0}}{\sqrt{T - k_0}}
        +
        \frac{d+d^3\gamma^{2k_0}}{T - k_0}
        +
        \epsilon^2
    \right),
\end{align*}
where the second and last equality can be shown via a similar way as in the proof of Theorem~\ref{thm:error_deterministic}.

Here, we have 
$
    k_0 = O \left(
        \frac{d}{\epsilon^4}
    \right)
$
by the definition of $k_0$. Thus, by setting
$
    T
    = 
    O\left(\frac{d^3}{\epsilon^2} + \frac{d^2}{\epsilon^4}\right)
    = 
    O\left(\frac{d^3}{\epsilon^2} + \frac{d^2}{\epsilon^4}\right)
$,
we can obtain
$
    \mathbb{E}_{\zeta, k'}[\|\nabla f(\hat{x})\|^2]
    =
    \frac{1}{T - k_0}\sum_{k=k_0 + 1}^{T}\mathbb{E}_{\zeta}[\|\nabla f(x_k)\|^2]
    \leq \epsilon^2
$.
% As we have
This implies 
$
    \mathbb{E}_{\zeta, k'}[\|\nabla f(\hat{x})\|] \leq \epsilon
$
as
$
    \mathbb{E}_{\zeta, k'}[\|\nabla f(\hat{x})\|]^2 
    \leq
    \mathbb{E}_{\zeta, k'}[\|\nabla f(\hat{x})\|^2]
$
follows from Jensen's inequality. Furthermore, when $\gamma$ is chosen as
$\gamma \leq d^{-\epsilon^4/d}$
, we have 
$
\log_\gamma d^{-1} = O\left( \frac{d}{\epsilon^4} \right)
$, which implies that 
$
k_0 = \Omega\left( \log_\gamma d^{-1} \right)
$. Therefore, we can obtain
$
\gamma^{k_0}
=
O( d^{-1} )
$, which yields the iteration complexity of $T=O\left(\frac{d}{\epsilon^4}\right)$.

\hfill $\Box$

\section{Optimization of test functions}
\label{sec:test_funcs}
In the first three subsections, let us compare the performance of our SLGH algorithms with GD-based algorithms and double loop GH algorithms using highly-non-convex test functions for optimization: the Ackley function \cit{molga2005test}, Rosenbrock function, and Himmelblau function \cit{Andrei08anunconstrained}. We implemented the following five types of algorithms: (ZOS)GD, (ZO)GradOpt, in which the factor for decreasing the smoothing parameter was 0.5 or 0.8, $\text{(ZO)SLGH}_{\text{r}}$ with $\gamma=0.995$ or $\gamma=0.999$.

\subsection{Ackley Function}
\begin{figure}[H]
\centering
 \label{fig:ackley_and_smoothed}
\begin{tabular}{cc}
    \subfigure[Ackley function]{
    \begin{minipage}[t]{0.45\linewidth}
    \label{fig:ackley}
    \centering
        \includegraphics[width=2.2in]{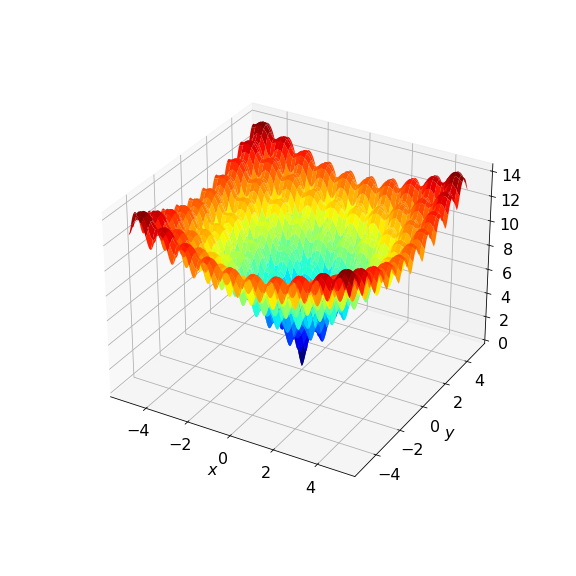}
    \end{minipage}%
    }
    \subfigure[Gaussian smoothed function with parameter $t=1$]{
    \begin{minipage}[t]{0.45\linewidth}
    \label{fig:ackley_smoothed}
    \centering
        \includegraphics[width=2.2in]{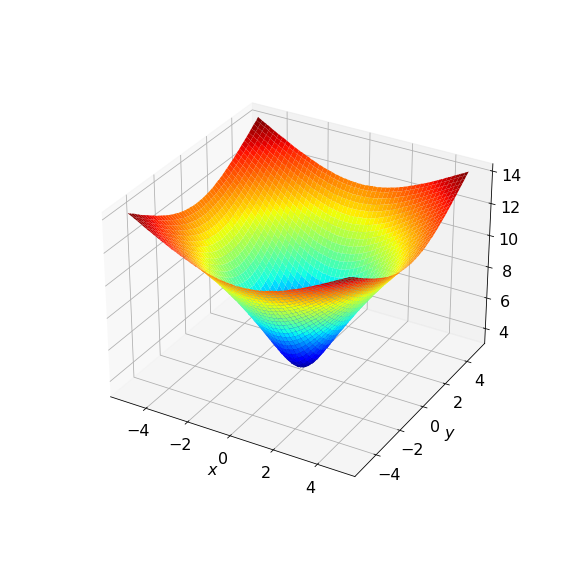}
     \end{minipage}
    }
\end{tabular}
\centering
\caption{Visualization of the Ackley function and its Gaussian smoothed function.}
\end{figure}

The Ackley function is defined as
\begin{align*}
f(x, y)=-20 \exp \left[-0.2 \sqrt{0.5\left(x^{2}+y^{2}\right)}\right]-\exp [0.5(\cos 2 \pi x+\cos 2 \pi y)]+e+20, 
\end{align*}
whose global optimum is $f(0,0)=0$. As shown in Figure~\ref{fig:ackley}, it has numerous small local minima due to cosine functions which are included in the second term. We ran the aforementioned five types of zeroth-order algorithms with the stepsize $\beta=0.1$ for $T=1000$ iterations. The initial smoothing parameter for the GH algorithms (ZOGradOpt and $\text{ZOSLGH}_{\text{r}}$) was set to $t_1=1$, where local minima of the smoothed function almost disappeared (Figure~\ref{fig:ackley_smoothed}). The smoothing parameter for ZOSGD was chosen as $t=0.005$. We set the initial point for the optimization as $(x,y)=(5,5)$.

\begin{figure}[H]
\centering
\begin{tabular}{cc}
    \subfigure[Function value $f(x,y)$ versus iterations.]{
    \begin{minipage}[t]{0.50\linewidth}
    \label{fig:fvalue_ackley}
    \centering
        \includegraphics[width=3in]{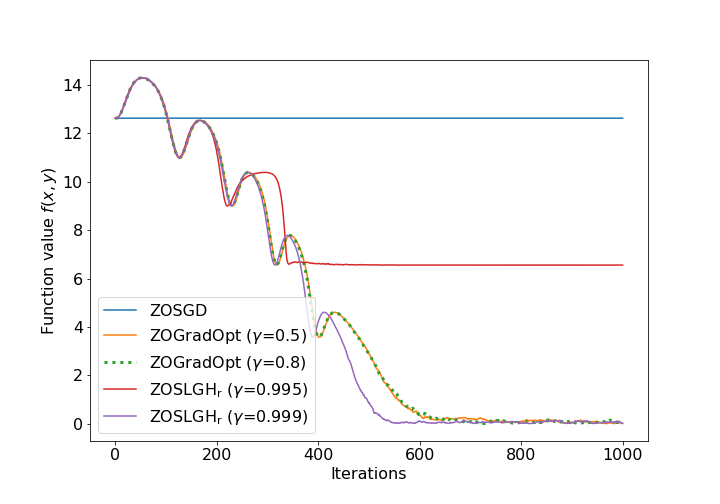}
    \end{minipage}%
    }
    \subfigure[Smoothing parameter $t$ versus iterations.]{
    \begin{minipage}[t]{0.50\linewidth}
    \label{fig:t_ackley}
    \centering
        \includegraphics[width=3in]{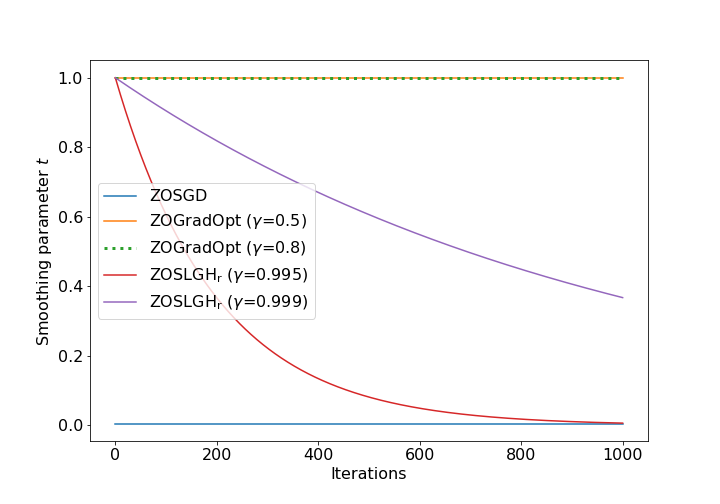}
     \end{minipage}
    }
\end{tabular}
\centering
\caption{Plots of the function value and the smoothing parameter during optimization of the Ackley function.}
\label{fig:results_ackley}
\end{figure}

We illustrate the optimization results in Table~\ref{table:results_ackley} and Figure~\ref{fig:results_ackley}. The GH methods successfully reach near the optimal solution $(0,0)$ when the decreasing speed of $t$ is not so fast, while ZOSGD is stuck in a local minimum in the immediate vicinity of the initial point $(5,5)$. Please note that GradOpt succeeds in optimization without decreasing the smoothing parameter since the optimal solution of the smoothed function with $t=1$ almost matches that of the original target function.

\begin{table}[H]
\centering
  \caption{Optimization results of the Ackley function. The global optimum is $f(0,0)=0$.}
\begin{tabular}{cc|c|c}
\toprule
&Methods & $(x,y)$ & $f(x,y)$ \\ \hline
SGD algo. &ZOSGD & $(4.99, 4.99)$ & $12.63$\\\hline
GH algo.&ZOGradOpt $(\gamma=0.5)$& $(4.2\times 10^{-3}, 1.9\times 10^{-3})$ & $\mathbf{1.4\times 10^{-2}}$\\
&ZOGradOpt $(\gamma=0.8)$& $(-2.2\times 10^{-3}, 6.7\times 10^{-3})$ & $\mathbf{8.1\times 10^{-2}}$\\ 
&$\text{ZOSLGH}_{\text{r}}\ (\gamma=0.995)$ & $(1.97, 1.97)$ & $6.56$\\
&$\text{ZOSLGH}_{\text{r}}\ (\gamma=0.999)$ & $(-3.6\times 10^{-3}, -4.6\times 10^{-3})$ & $\mathbf{1.7\times 10^{-2}}$\\
\bottomrule
\end{tabular}%
\label{table:results_ackley}%
\end{table}%

\subsection{Rosenbrock Function}
\begin{figure}[H]
\label{fig:rosenbrock_and_smoothed}
\centering
\begin{tabular}{cc}
    \subfigure[Rosenbrock function]{
    \begin{minipage}[t]{0.45\linewidth}
    \label{fig:rosenbrock}
    \centering
        \includegraphics[width=2.2in]{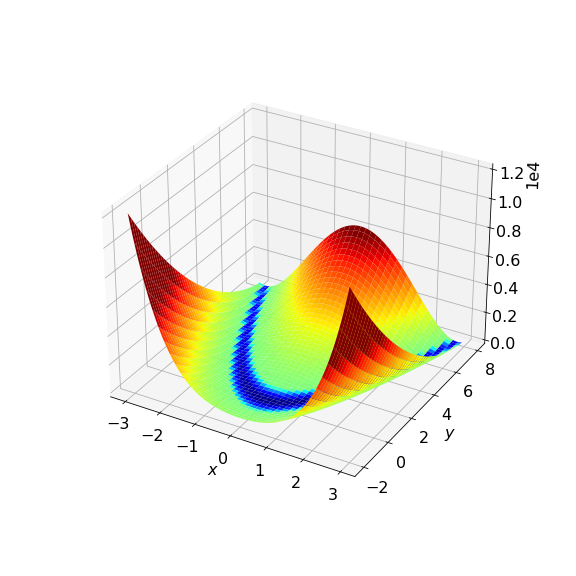}
    \end{minipage}%
    }
    \subfigure[Smoothed function with parameter $t=1.5$]{
    \begin{minipage}[t]{0.45\linewidth}
    \label{fig:rosenbrock_smoothed}
    \centering
        \includegraphics[width=2.2in]{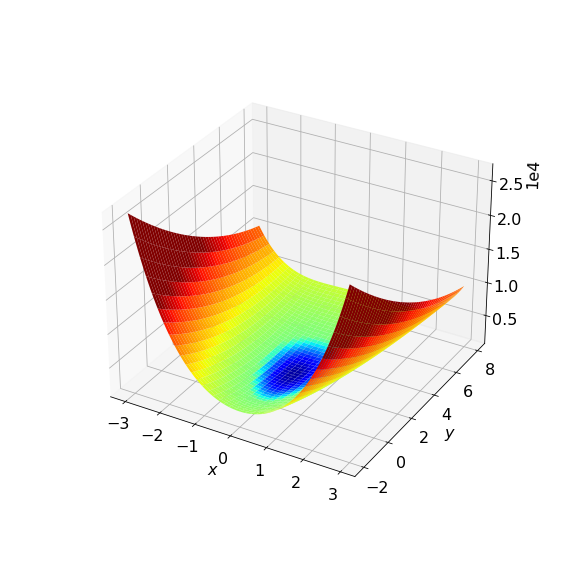}
     \end{minipage}
    }
\end{tabular}
\centering
\caption{Visualization of the Rosenbrock function and its Gaussian smoothed function.}
\end{figure}

Let us define the Rosenbrock function in 2D as
\begin{align*}
f(x,y)=100\left(y-x^{2}\right)^{2}+\left(1-x\right)^{2},
\end{align*}
whose global optimum is $f(1,1)=0$. This function is difficult to optimize because the global optimum lies inside a flat parabolic shaped valley with low function value (Figure~\ref{fig:rosenbrock}). Since this function is polynomial, we can calculate the GH smoothed function analytically (see \cit{mobahi2012gaussian}):
\begin{align*}
    F(x,y,t) &:= \mathbb{E}_{u_x, u_y}[f(x+tu_x, y+tu_y)],\quad \left(u_x, u_y\sim\mathcal{N}(0, 1)\right)\\
    &= 100x^4+(-200y+600t^2+1)x^2-2x+100y^2-200t^2y+(300t^4+101t^2+1).
\end{align*}

Thus, we applied first-order methods to this function. The stepsize and iteration number were set to $\beta=1\times 10^{-4}$ and $T=20000$, respectively. The initial smoothing parameter for the GH algorithms (GradOpt and $\text{SLGH}_{\text{r}}$) was set to $t_1=1.5$, where the smoothed function became almost convex around the optimal solution (Figure~\ref{fig:rosenbrock_smoothed}). We set the initial point for the optimization as $(x,y)=(-3,2)$.

\begin{figure}[H]
\centering
\begin{tabular}{cc}
    \subfigure[Function value $f(x,y)$ versus iterations.]{
    \begin{minipage}[t]{0.50\linewidth}
    \label{fig:fvalue_rosenbrock}
    \centering
        \includegraphics[width=3in]{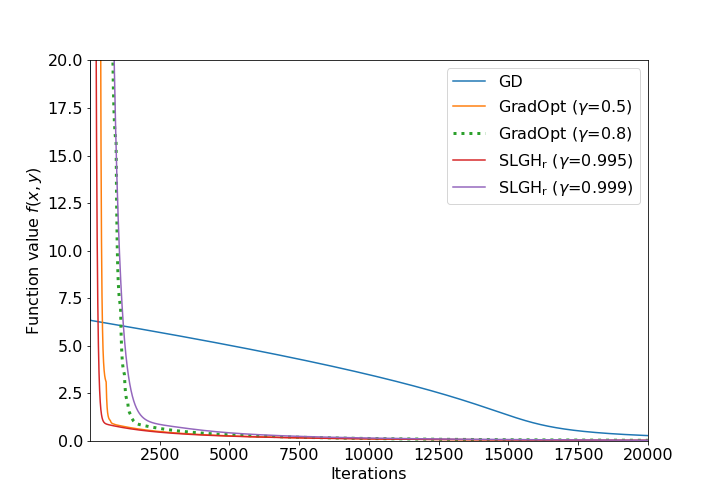}
    \end{minipage}%
    }
    \subfigure[Smoothing parameter $t$ versus iterations.]{
    \begin{minipage}[t]{0.50\linewidth}
    \label{fig:t_rosenbrock}
    \centering
        \includegraphics[width=3in]{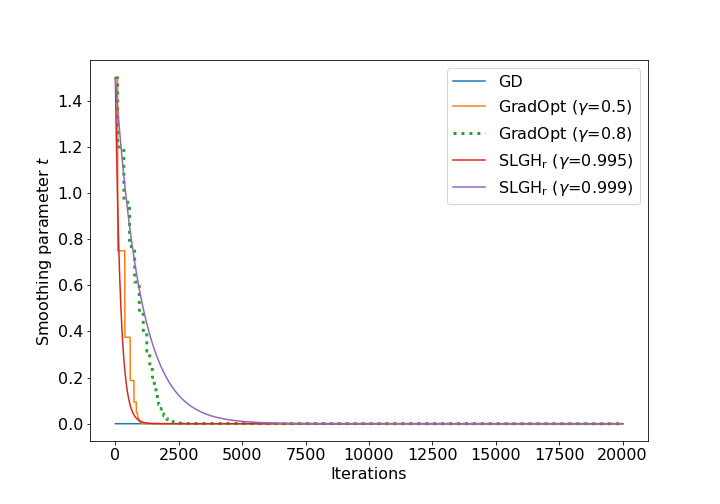}
     \end{minipage}
    }
\end{tabular}
\caption{Plots of the function value and the smoothing parameter during optimization of the Rosenbrock function.}
\label{fig:results_rosenbrock}
\end{figure}

\begin{figure}[H]
\centering
\includegraphics[width=3in]{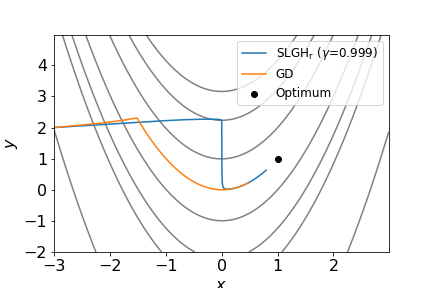}
%\caption{Comparison of output sequences between GD and $\text{SLGH}_{\text{r}}\ (\gamma=0.999)$ with contours of the Rosenbrock function.}
\caption{\begin{tabular}{l}
    Comparison of output sequences between GD and $\text{SLGH}_{\text{r}}\  (\gamma=0.999)$ with contours \\of the Rosenbrock function.
\end{tabular}}
\label{fig:x_rosenbrock}
\end{figure}

We illustrate the optimization results in Table~\ref{table:results_rosenbrock}, Figure~\ref{fig:results_rosenbrock} and Figure~\ref{fig:x_rosenbrock}. The GH methods can decrease the function value much faster than GD. This is because the smoothed function is much easier to optimize than the original function while its optimal solution is close to that of the original one. In the early stage of optimization, the GH methods reach near a point $(0,2)$, which is a good initial point for optimization, while GD falls into a point in the flat valley, which is far from the optimal solution.  (Figure~\ref{fig:x_rosenbrock}).

\begin{table}[H]
\centering
  \caption{Optimization results of the Rosenbrock function. The global optimum is $f(1,1)=0$.}
\begin{tabular}{cc|c|c}
\toprule
&Methods & $(x,y)$ & $f(x,y)$ \\ \hline
GD algo. &GD & $(0.468, 0.216)$ & $0.284$\\\hline
GH algo.&GradOpt $(\gamma=0.5)$& $(0.817, 0.667)$ & $\mathbf{3.36\times 10^{-2}}$\\
&GradOpt $(\gamma=0.8)$& $(0.808, 0.652)$ & $\mathbf{3.70\times 10^{-2}}$\\ 
&$\text{SLGH}_{\text{r}}\ (\gamma=0.995)$ & $(0.819, 0.670)$ & $\mathbf{3.27\times 10^{-2}}$\\
&$\text{SLGH}_{\text{r}}\ (\gamma=0.999)$ & $(0.795, 0.631)$ & $\mathbf{4.19\times 10^{-2}}$\\
\bottomrule
\end{tabular}%
\label{table:results_rosenbrock}%
\end{table}%

\subsection{Himmelblau Function}
\begin{figure}[H]
\centering
 \label{fig:himmelblau_and_smoothed}
\begin{tabular}{cc}
    \subfigure[Himmelblau function]{
    \begin{minipage}[t]{0.45\linewidth}
    \label{fig:himmelblau}
    \centering
        \includegraphics[width=2.4in]{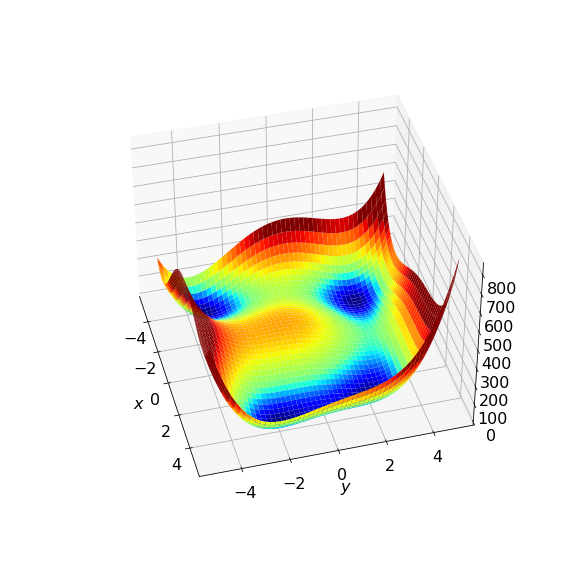}
    \end{minipage}
    }
    \subfigure[Gaussian smoothed function with parameter $t=2$]{
    \begin{minipage}[t]{0.45\linewidth}
    \label{fig:himmelblau_smoothed}
    \centering
        \includegraphics[width=2.4in]{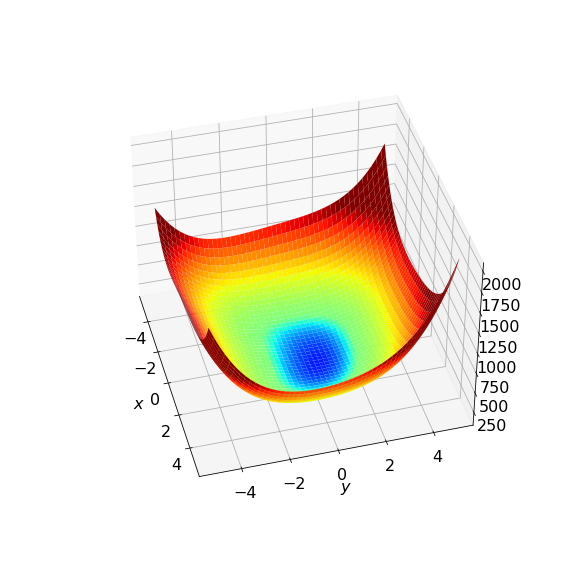}
     \end{minipage}
    }
\end{tabular}
\centering
\caption{Visualization of the Himmelblau function and its Gaussian smoothed function.}
\end{figure}

The Himmelblau function is defined as
\begin{align*}
f(x, y)= (x^2+y-11)^2+(x+y^2-7)^2. 
\end{align*}
It has four minimum points in the vicinity of $(x,y)=(3.000, 2.000), (-2.805, 3.131), (-3.779, -3.283)$, $(3.584, -1.848)$ and one maximum point in the vicinity of $(x,y)=(-0.271, -0.923)$. It takes the optimal value $0$ at the four points. Since this function is also polynomial, we can calculate the GH smoothed function analytically:
\begin{align*}
    & F(x,y,t):= \mathbb{E}_{u_x, u_y}[f(x+tu_x, y+tu_y)],\quad \left(u_x, u_y\sim\mathcal{N}(0, 1)\right)\\
    &= x^4+(2y+6t^2-21)x^2+(2y^2+2t^2-14)x+y^4+(6t^2-13)y^2+(2t^2-22)y + (6t^4-34t^2+170).
\end{align*}
Thus, we applied first-order methods to this function. The stepsize and iteration number were set to $\beta=1\times 10^{-4}$ and $T=2000$, respectively. The initial smoothing parameter for GH algorithms was set to $t_1=2$, where the smoothed function became almost convex around the optimal solution (Figure~\ref{fig:himmelblau_smoothed}). We set the initial point for the optimization as $(x,y)=(5,5)$.

\begin{figure}[H]
\centering
\begin{tabular}{cc}
    \subfigure[Function value $f(x,y)$ versus iterations.]{
    \begin{minipage}[t]{0.50\linewidth}
    \label{fig:fvalue_himmelblau}
    \centering
        \includegraphics[width=3in]{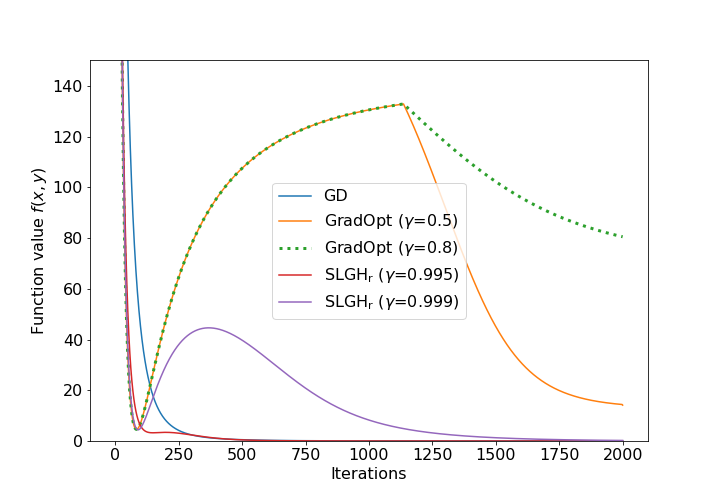}
    \end{minipage}%
    }
    \subfigure[Smoothing parameter $t$ versus iterations.]{
    \begin{minipage}[t]{0.50\linewidth}
    \label{fig:t_himmelblau}
    \centering
        \includegraphics[width=3in]{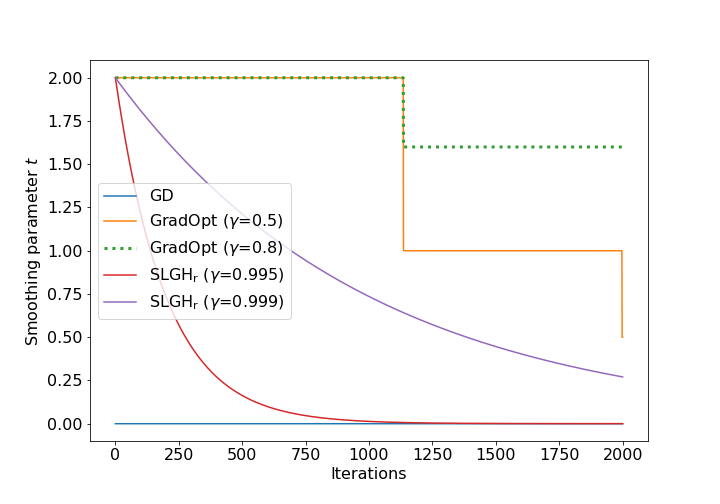}
     \end{minipage}
    }
\end{tabular}
\centering
\caption{Plots of the function value and the smoothing parameter during optimization of the Himmelblau function.}
\label{fig:results_himmelblau}
\end{figure}

\begin{figure}[H]
\centering
\includegraphics[width=3in]{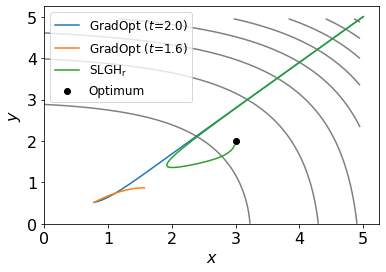}
\caption{Comparison of output sequences of GradOpt, in which the factor for decreasing the smoothing parameter is $0.8$, and $\text{SLGH}_{\text{r}}\ (\gamma=0.999)$ with contours of the {\it smoothed} Himmelblau function. The blue and orange lines represent output sequences of GradOpt when $t=2.0$ and $t=1.6$, respectively.}
\label{fig:x_himmelblau}
\end{figure}

Table~\ref{table:results_himmelblau}, Figure~\ref{fig:results_himmelblau}, and Figure~\ref{fig:x_himmelblau} show the optimization results. GD and our SLGH algorithms successfully reach near the global optimum, while GradOpt fails to decrease the function value. This is because the optimal solution of the smoothed function when $t=2$ lies near the maximum point of the original Himmelblau function $(-0.271, -0.923)$. Figure~\ref{fig:x_himmelblau} describes detailed optimization process. Our SLGH algorithm succeeds in returning to the optimal solution once it has passed by reducing $t$. In contrast, GradOpt reaches the vicinity of a minimum of the smoothed function without knowing the detailed shape of the original function; as a result, it is stuck around a local maximum of the original function.

\begin{table}[H]
\centering
  \caption{Results of optimization of the Himmelblau function. It has a global optimum $f(3,2)=0$.}
\begin{tabular}{cc|c|c}
\toprule
&Methods & $(x,y)$ & $f(x,y)$ \\ \hline
GD algo. &GD & $(2.998, 2.003)$ & $\mathbf{1.6\times 10^{-4}}$\\\hline
GH algo.&GradOpt $(\gamma=0.5)$& $(2.575, 1.437)$ & $14.14$\\
&GradOpt $(\gamma=0.8)$& $(1.573, 0.868)$ & $80.51$\\ 
&$\text{SLGH}_{\text{r}}\ (\gamma=0.995)$ & $(2.999, 2.002)$ & $\mathbf{6.9\times 10^{-5}}$\\
&$\text{SLGH}_{\text{r}}\ (\gamma=0.999)$ & $(2.983, 1.897)$ & $\mathbf{0.21}$\\
\bottomrule
\end{tabular}%
\label{table:results_himmelblau}%
\end{table}%

\subsection{Additional Toy Example}
\label{subsec:toy}
At the end of this section, let us present a toy example problem in which $\text{SLGH}_{\text{d}}$, which utilizes the derivative $\frac{\partial F}{\partial t}$ for the update of $t$, outperforms $\text{SLGH}_{\text{d}}$. Let us consider the following artificial non-convex function:
\begin{align*}
    f(x, y)=\left\{\begin{array}{cc}
x^2 - 150\times 1.1^{-((x-10)^2+y^2)}& (x\geq0) \\
x^2/50 - 150\times 1.1^{-((x-10)^2+y^2)}& (x<0)\\
\end{array}\right..
\end{align*}

The second term creates a hole around $(x,y)=(10,0)$ (see Figure\ref{fig:toy}), and this function has an optimum in the vicinity of $f(9.319, 0)\simeq -56.670$. This function is difficult to optimize for GH methods since the hole around the optimum disappears when the smoothing parameter $t$ is large (Figure\ref{fig:toy_smoothed}). 

We ran $\text{SLGH}_{\text{r}}\ (\gamma=0.995 \text{ or } 0.999)$ and $\text{SLGH}_{\text{d}}\ (\gamma=0.999)$ with the stepsize (for $x$) $\beta=0.01$ for $T=1000$ iterations. The initial point and initial smoothing parameter were set to $(x,y)=(15,0)$ and $t_1=5$, respectively. We set the stepsize for $t$ as $0.01$.

\begin{figure}[H]
\centering
\begin{tabular}{cc}
    \subfigure[Function value $f(x,y)$ versus iterations.]{
    \begin{minipage}[t]{0.50\linewidth}
    \label{fig:fvalue_toy}
    \centering
        \includegraphics[width=3in]{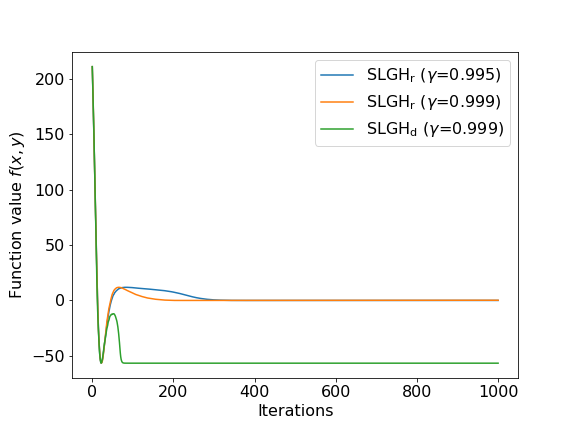}
    \end{minipage}%
    }
    \subfigure[Smoothing parameter $t$ versus iterations.]{
    \begin{minipage}[t]{0.50\linewidth}
    \label{fig:t_toy}
    \centering
        \includegraphics[width=3in]{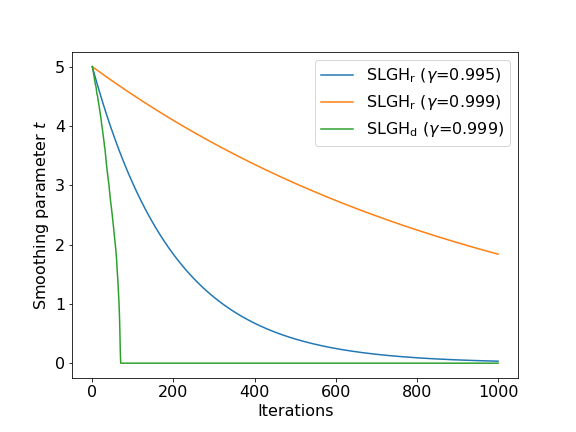}
     \end{minipage}
    }
\end{tabular}
\centering
\caption{\begin{tabular}{l}Plots of the function value and the smoothing parameter during optimization\\ of the artificial non-convex function.\end{tabular}}
\label{fig:results_toy}
\end{figure}

\begin{figure}[H]
    \centering
     \label{fig:toy_and_smoothed}
    \begin{tabular}{cc}
        \subfigure[Artificial non-convex function]{
        \begin{minipage}[t]{\linewidth}
        \label{fig:toy}
        \centering
            \includegraphics[width=6in]{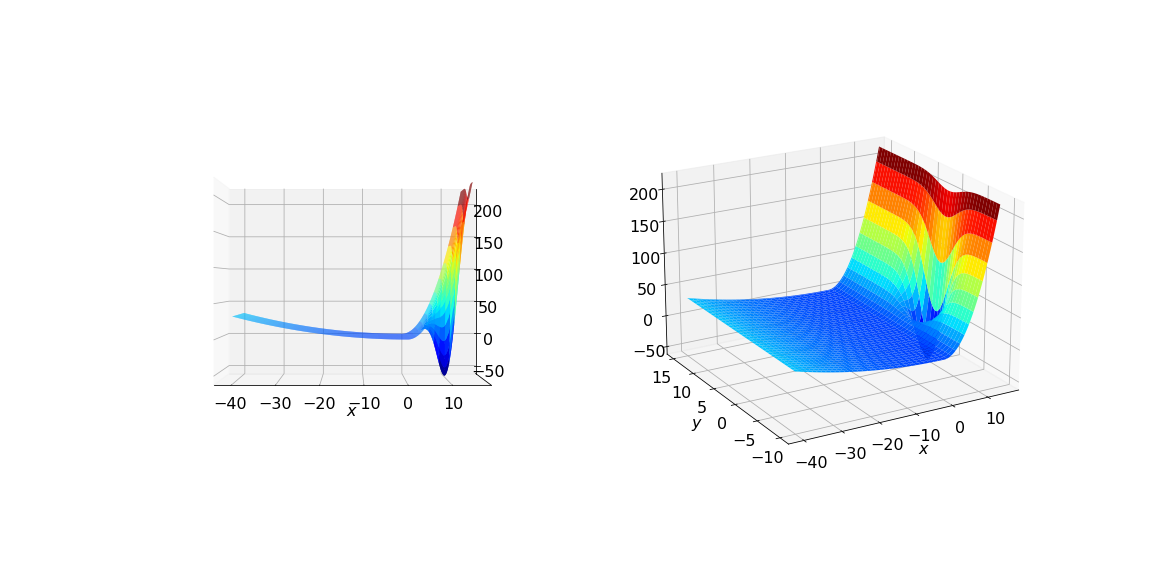}
        \end{minipage}%
        }
        \\
        \subfigure[Gaussian smoothed function with parameter $t=5$]{
        \begin{minipage}[t]{\linewidth}
        \label{fig:toy_smoothed}
        \centering
            \includegraphics[width=6in]{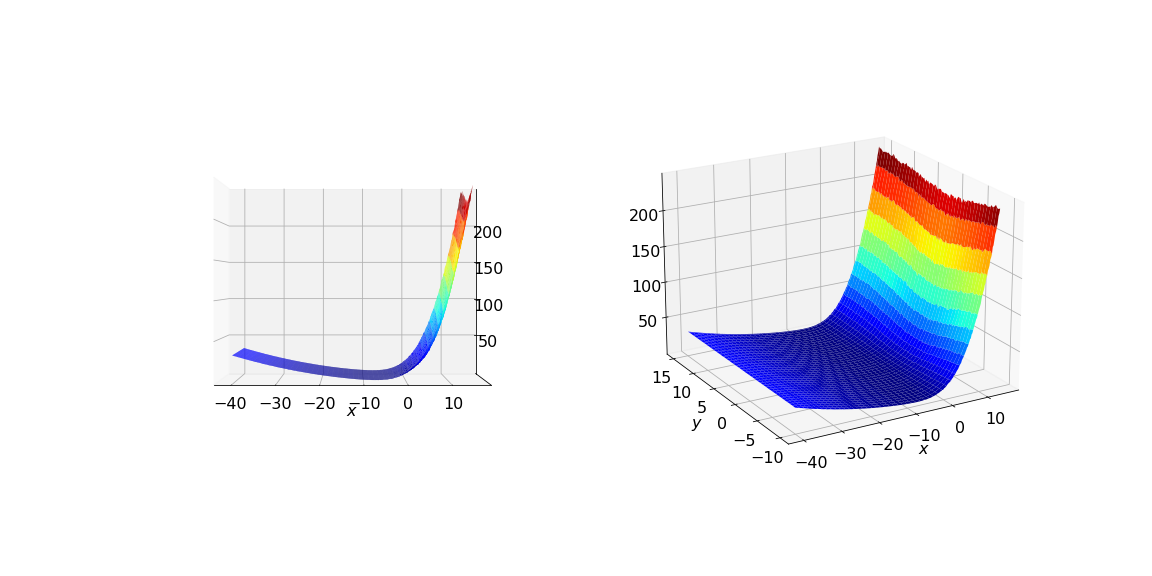}
         \end{minipage}
        }
    \end{tabular}
    \caption{Visualization of the artificial non-convex function and its Gaussian smoothed function.}
    \end{figure}

Table~\ref{table:results_toy} and Figure~\ref{fig:results_toy} show the optimization results. We can see that only $\text{SLGH}_{\text{d}}$ can decrease $t$ around the hole adaptively, and thus successfully can find the optimal solution.

\begin{table}[H]
\centering
  \caption{Optimization results of the artificial non-convex function. It has a global optimum in the vicinity of $f(9.319, 0)\simeq -56.670$.}
\begin{tabular}{cc|c|c}
\toprule
&Methods & $(x,y)$ & $f(x,y)$ \\ \hline
GH algo.&$\text{SLGH}_{\text{r}}\ (\gamma=0.995)$ & $(-0.248, 2.38\times 10^{-2})$ & $-5.52\times 10^{-3}$\\
&$\text{SLGH}_{\text{r}}\ (\gamma=0.999)$ & $(-2.959, -2.18\times 10^{-3})$ & $0.175$\\
&$\text{SLGH}_{\text{d}}\ (\gamma=0.999)$ & $(9.319, 8.33\times 10^{-3})$ & $\mathbf{-56.670}$\\ 
\bottomrule
\end{tabular}%
\label{table:results_toy}%
\end{table}%

\newpage

\section{Black-box adversarial attack}
\label{sec:black_box}
\subsection{Experimental Setup}
We used well-trained $\mathop{\rm DNNs}$\footnote{https://github.com/carlini/nn\_robust\_attacks} for CIFAR10 and MNIST classification tasks as target models, respectively. We adopt the implementation\footnote{https://github.com/KaidiXu/ZO-AdaMM} in \cit{chen2019zo} for ZOSGD and ZOAdaMM. GradOpt \cit{hazan2016graduated} in our implementation adopts the same random gradient-free oracles \cit{nesterov2017random} as with our ZOSLGH methods, rather than their smoothed gradient oracle, where random variables are sampled from the unit sphere. Moreover, we set the stepsize in its inner loop as a constant instead of $\Theta(1/k)$, where $k$ denotes an iteration number in the inner loop, due to less efficiency of the original setting. Therefore, the essential difference between GradOpt and $\text{ZOSLGH}_\text{r}$ is whether or not the structure of algorithms is single loop.

As recommended in their work, we set the parameter for ZOAdaMM as ${v}_0 = 10^{-5}$, $\beta_1=0.9$, and $\beta_2=0.3$. The factor for decreasing the smoothing parameter in ZOGradOpt was set to $0.5$. For all algorithms, we chose the regularization parameter $\lambda$ as $\lambda=10$ and set attack confidence $\kappa=1e-10$. We chose minibatch size as $M=10$ to stabilize estimation of values and gradients of the smoothed function. The initial adversarial perturbation was chosen as $x_0=0$, and the initial smoothing parameter $t_0$ was $10$ for GH methods and $0.005$ for the others. The decreasing factor for $t$ in the ZOSLGH algorithm was set to $\gamma=0.999$ for both of $\text{ZOSLGH}_\text{r}$ and $\text{ZOSLGH}_\text{d}$, unless otherwise noted. Other parameter settings are described in Table \ref{tab:parameter_setting}. We used different step sizes for ZOAdaMM because it adaptively penalizes the step size using the information of past gradients \cit{chen2019zo}. 
%For fair comparison, we multiply $\beta$ by average of gradients of the loss functions obtained by running ZOSGD for 100 iterations, only when ZOAdaMM is used.
%\memo{HI}{Is the last sentence necessary?}

\begin{table}[H]
\centering
\caption{Parameter settings in the adversarial attack problems. $T$ represents the iteration number. $\beta$ is the step size for $x$, and $\eta$ is the step size for $t$. $N_0$ and $\epsilon_0$ are used to determine termination condition of the inner loop in ZOGradOpt: we stop the inner loop and decrease $t$ if the condition $|\frac{1}{M}\sum_{i=1}^Mf(x_{k+1}+tu_i)-\frac{1}{M}\sum_{i=1}^Mf(x_{k}+tu'_i)|\leq\epsilon_0$ is satisfied $N_0$ times, where $u_i$ and $u_i'\ (i=1, ..., M)$ are sampled from $\mathcal{N}(0, \mathrm{I}_d)$. Each of ``3072'' and ``784'' is the dimension of images in CIFAR-10 and MNIST.}
\label{tab:parameter_setting}
\begin{tabular}{c|c|c|c|c|c}
\toprule
& $T$ & \begin{tabular}{c} $\beta$ \\ (other than\\ ZOAdaMM) \end{tabular} & \begin{tabular}{c} $\beta$ \\ (for ZOAdaMM) \end{tabular} & $\eta$  & $(N_0, \epsilon_0)$ 
\\ \hline
CIFAR-10 & $10000$ & $0.01/3072$ & $0.5/3072$ & $1\times10^{-4}/3072$ & $(100, 5\times10^{-3})$\\
MNIST & $20000$ & $1/784$ & $100/784$ & $0.1/784$ & $(100, 1\times10^{-3})$\\ 
\bottomrule
\end{tabular}%

% \begin{tabular}{c|c|c|c}
% \toprule
% & $T$ & $\beta$ (other than ZOAdaMM) & $\beta$ (for ZOAdaMM) \\ \hline
% CIFAR-10 & 10000 & 0.01/3072 & 0.5/3072 \\
% MNIST & 20000 & 1/784 & 100/784 \\ 
% \bottomrule
% \end{tabular}%
% 
% 
% \begin{tabular}{c|c|c|c}
%     \toprule
%     & $\eta$ & $(N_0, \epsilon_0)$ 
%     \\ \hline
%     CIFAR-10 & $1\times10^{-4}/3072$  & $(100, 5\times10^{-3})$\\
%     MNIST & 0.1/784 & $(100, 1\times10^{-3})$\\ 
%     \bottomrule
% \end{tabular}%

\end{table}%

\subsection{CIFAR-10}
\label{sec:cifar10}
\begin{figure}[H]
    \vspace{-10mm}
\subfigure[CIFAR-10, Image ID = 1]{
    \label{fig:easy01_cifar10}
    \begin{minipage}[t]{0.47\linewidth}
        \centering
        \includegraphics[width=2.0in]{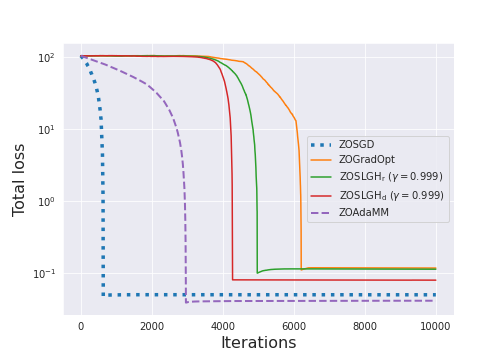}
     \end{minipage}
    }
    \vspace{-1mm}
    \subfigure[CIFAR-10, Image ID = 16]{
    \label{fig:easy1_cifar10}
    \begin{minipage}[t]{0.50\linewidth}
        \centering
        \includegraphics[width=2.0in]{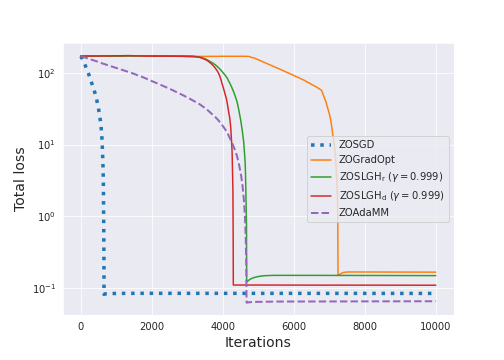}
    \end{minipage}%
    }
    \vspace{-1mm}
    \subfigure[CIFAR-10, Image ID = 37]{
    \label{fig:easy2_cifar10}
    \begin{minipage}[t]{0.47\linewidth}
        \centering
        \includegraphics[width=2.0in]{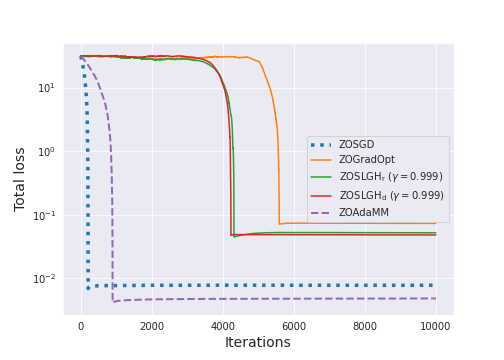}
     \end{minipage}
    }
\vspace{-1mm}
    \subfigure[CIFAR-10, Image ID = 7]{
    \label{fig:easy02_cifar10}
    \begin{minipage}[t]{0.50\linewidth}
        \centering
        \includegraphics[width=2.0in]{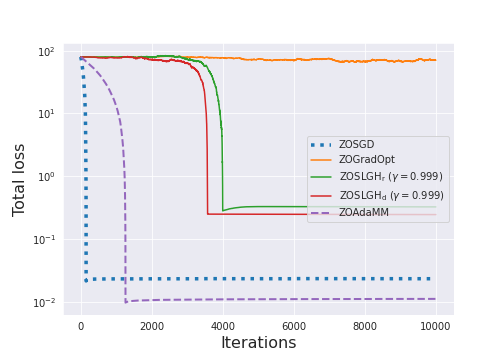}
     \end{minipage}
    }
\vspace{-1mm} 
    \subfigure[CIFAR-10, Image ID = 51]{
    \label{fig:mid1_cifar10}
    \begin{minipage}[t]{0.47\linewidth}
        \centering
        \includegraphics[width=2.0in]{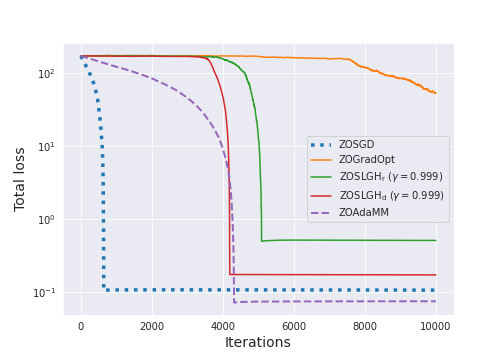}
    \end{minipage}%
    }
    \vspace{-1mm}
    \subfigure[CIFAR-10, Image ID = 96]{
    \label{fig:mid2_cifar10}
    \begin{minipage}[t]{0.50\linewidth}
        \centering
        \includegraphics[width=2.0in]{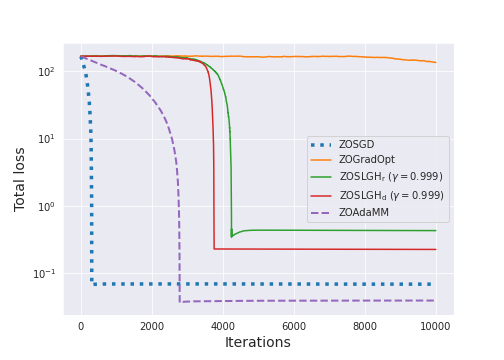}
     \end{minipage}
    }
    
    \subfigure[CIFAR-10, Image ID = 39]{
    \label{fig:hard1_cifar10}
    \begin{minipage}[t]{0.47\linewidth}
        \centering
        \includegraphics[width=2.0in]{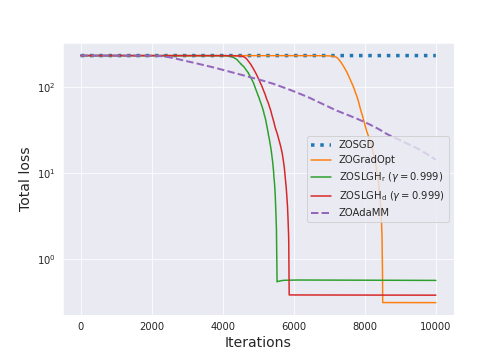}
    \end{minipage}%
    }
    \vspace{-1mm}
    \subfigure[CIFAR-10, Image ID = 104]{
    \label{fig:hard2_cifar10}
    \begin{minipage}[t]{0.50\linewidth}
        \centering
        \includegraphics[width=2.0in]{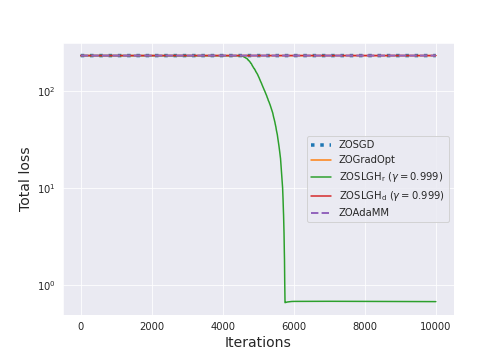}
     \end{minipage}
    }
     
    \subfigure[CIFAR-10, Image ID = 14]{
    \label{fig:veryhard1_cifar10}
    \begin{minipage}[t]{0.45\linewidth}
        \begin{center}
        \includegraphics[width=1.9in]{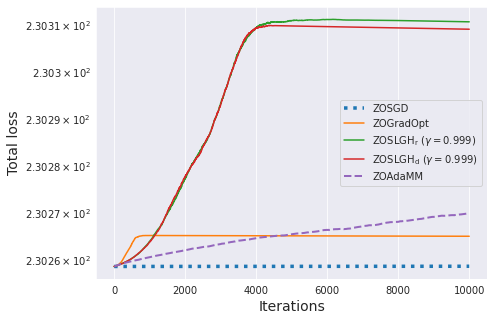}
        \end{center}
     \end{minipage}
    }
    \vspace{-1mm}
    \subfigure[CIFAR-10, Image ID = 41]{
    \label{fig:veryhard2_cifar10}
    \begin{minipage}[t]{0.51\linewidth}
        \begin{center}
        \includegraphics[width=1.9in]{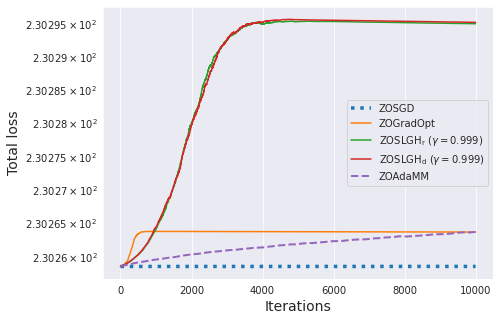}
        \end{center}
     \end{minipage}
    }
\centering
\caption{Additional plots of total loss versus iterations on CIFAR-10 (log scale). (a)-(c) All algorithms can successfully decrease the loss value when images are easy to attack. In particular, in plot (c), SGD-based algorithms can find better solutions than GH-based algorithms.  (d)-(f) Only GradOpt fails to attack due to its slow convergence. (g) Only ZOSGD is stuck around a local minimum $x=0$. (h) Only our $\text{ZOSLGH}_\text{r}$ algorithm succeeds in escaping the local minimum, and thus it can decrease the loss value more than 200 than other algorithms. (i), (j): These images are so difficult to attack that no algorithms can succeed in attacking.}
\label{fig:total_plot_cifar10}
\end{figure}

\begin{figure}[H]
    \subfigure[CIFAR-10, Image ID = 39]{
    \label{fig:hard39_l2_cifar10}
    \begin{minipage}[t]{0.47\linewidth}
        \centering
        \includegraphics[width=2.5in]{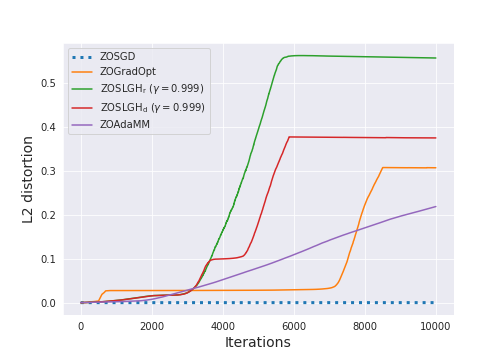}
    \end{minipage}%
    }
    \subfigure[CIFAR-10, Image ID = 104]{
    \label{fig:hard104_l2_cifar10}
    \begin{minipage}[t]{0.47\linewidth}
        \centering
        \includegraphics[width=2.5in]{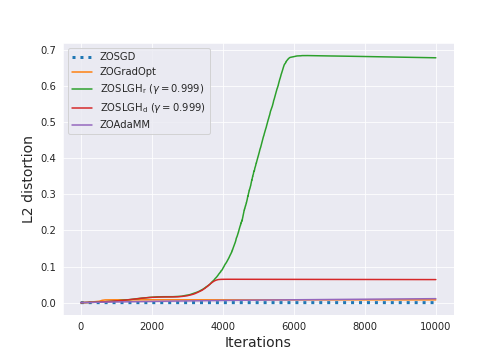}
     \end{minipage}
    }
    
    \subfigure[CIFAR-10, Image ID = 14]{
    \label{fig:veryhard14_l2_cifar10}
    \begin{minipage}[t]{0.47\linewidth}
        \centering
        \includegraphics[width=2.5in]{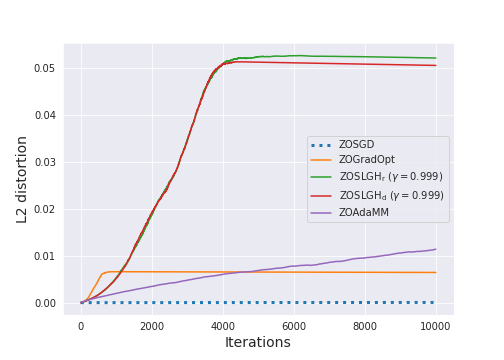}
    \end{minipage}%
    }
    \subfigure[CIFAR-10, Image ID = 41]{
    \label{fig:veryhard41_l2_cifar10}
    \begin{minipage}[t]{0.47\linewidth}
        \centering
        \includegraphics[width=2.5in]{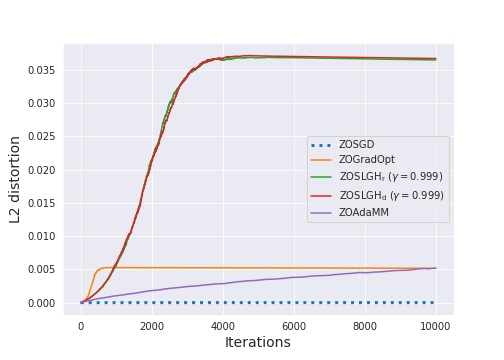}
     \end{minipage}
    }
\centering
\caption{Plots of $L_2$ distortion versus iterations for images that are difficult to attack on CIFAR-10. Each plot of (a)-(d) corresponds to Figure~\ref{fig:hard1_cifar10}-Figure~\ref{fig:veryhard2_cifar10}.}
\label{fig:l2_plot_cifar10}
\end{figure}

\paragraph{Additional plots}
Figures~\ref{fig:total_plot_cifar10} and \ref{fig:l2_plot_cifar10} show additional plots for total loss and $L_2$ distortion, respectively. We can see that our ZOSLGH algorithms successfully decrease the total loss value except in cases where images are so difficult to attack that no algorithms succeed in attacking (Figure~\ref{fig:veryhard1_cifar10}, \ref{fig:veryhard2_cifar10}). Plots in Figure  \ref{fig:l2_plot_cifar10} imply that the algorithms are stuck around a local minimum $x=0$ when they are failed to decrease the loss value.

\begin{figure}[H]
    \subfigure[CIFAR-10, Image ID = 8]{
    \label{fig:comp_cifar10_id8}
    \begin{minipage}[t]{0.47\linewidth}
        \centering
        \includegraphics[width=2.5in]{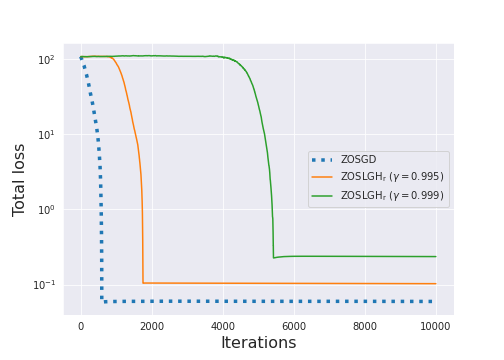}
     \end{minipage}
    }
    \subfigure[CIFAR-10, Image ID = 66]{
    \label{fig:comp_cifar10_id66}
    \begin{minipage}[t]{0.47\linewidth}
        \centering
        \includegraphics[width=2.5in]{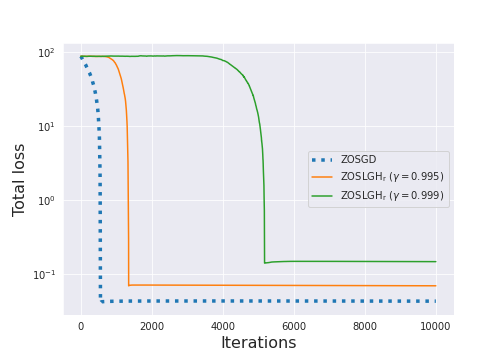}
    \end{minipage}%
    }
    \subfigure[CIFAR-10, Image ID = 105]{
    \label{fig:comp_cifar10_id105}
    \begin{minipage}[t]{0.47\linewidth}
        \centering
        \includegraphics[width=2.5in]{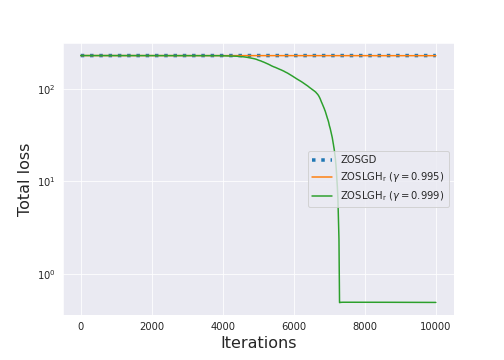}
    \end{minipage}%
    }
    \subfigure[CIFAR-10, Image ID = 89]{
    \label{fig:comp_cifar10_id89}
    \begin{minipage}[t]{0.47\linewidth}
        \centering
        \includegraphics[width=2.5in]{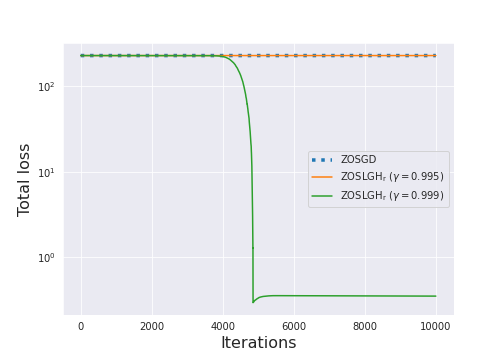}
     \end{minipage}
    }
\centering
\caption{Comparison of total loss transition of ZOSGD, $\text{ZOSLGH}_\text{r}$ with $\gamma=0.995$, and $\text{ZOSLGH}_\text{r}$ with $\gamma=0.999$ (log scale).}
\label{fig:comparison_cifar10}
\end{figure}

%%%%%%%%%%%%%%%%%%%%%%%%%%%%%%%%%%%%%%%%%%%%%%%%%%%%%%%%%%%%%%%%%%%%%
\newpage
\paragraph{Effect of choice of the parameter $\gamma$ in the ZOSLGH algorithm}
We also investigated the effect of choice of the decreasing parameter $\gamma$ in the ZOSLGH algorithm. We compared ZOSGD, $\text{ZOSLGH}_\text{r}$ with $\gamma=0.995$, and $\text{ZOSLGH}_\text{r}$ with $\gamma=0.999$. All other parameters were set to the same values as before.  Figure~\ref{fig:comparison_cifar10} implies that the decreasing speed of $t$ is associated with a trade-off: a rapid decrease of $t$ yields fast convergence, but reduces the possibility to find better solutions.
\newpage
\paragraph{Generated adversarial examples}
Table~\ref{tab:cifar10_attack_examples} shows adversarial images generated by different algorithms and their original images.

\begin{table}[H]
  \caption{Comparison of adversarial images for CIFAR-10 with different algorithms.}
\centering
\begin{tabular}{ccccc}
\toprule
Image ID & 39 & 79 & 89& 115\\ \hline
Original  & \begin{minipage}[b]{0.085\columnwidth}
		\centering
		\raisebox{-.5\height}{\includegraphics[width=\linewidth]{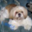}}
	\end{minipage} &
	\begin{minipage}[b]{0.085\columnwidth}
		\centering
		\raisebox{-.5\height}{\includegraphics[width=\linewidth]{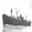}}
	\end{minipage} &
	\begin{minipage}[b]{0.085\columnwidth}
		\centering
		\raisebox{-.5\height}{\includegraphics[width=\linewidth]{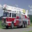}}
	\end{minipage} & \begin{minipage}[b]{0.085\columnwidth}
		\centering
		\raisebox{-.5\height}{\includegraphics[width=\linewidth]{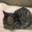}}
	\end{minipage}\\ \ 
Classified as&dog & ship & truck & cat\\ 
$L_2$ distortion: &0 & 0 & 0 & 0\\ \hline

ZOSGD  & \begin{minipage}[b]{0.085\columnwidth}
		\centering
		\raisebox{-.5\height}{\includegraphics[width=\linewidth]{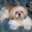}}
	\end{minipage} &
	\begin{minipage}[b]{0.085\columnwidth}
		\centering
		\raisebox{-.5\height}{\includegraphics[width=\linewidth]{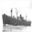}}
	\end{minipage} & \begin{minipage}[b]{0.085\columnwidth}
		\centering
		\raisebox{-.5\height}{\includegraphics[width=\linewidth]{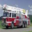}}
	\end{minipage} & \begin{minipage}[b]{0.085\columnwidth}
		\centering
		\raisebox{-.5\height}{\includegraphics[width=\linewidth]{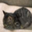}}
	\end{minipage}\\
Classified as&dog (fail.)& airplane & truck (fail.) & horse\\ 
$L_2$ distortion: &$6.7\times 10^{-5}$ & $0.154$ & $5.6\times 10^{-5}$ & $4.5\times 10^{-3}$\\ \hline
ZOAdaMM  & \begin{minipage}[b]{0.085\columnwidth}
		\centering
		\raisebox{-.5\height}{\includegraphics[width=\linewidth]{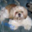}}
	\end{minipage} &
	\begin{minipage}[b]{0.085\columnwidth}
		\centering
		\raisebox{-.5\height}{\includegraphics[width=\linewidth]{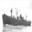}}
	\end{minipage} & \begin{minipage}[b]{0.085\columnwidth}
		\centering
		\raisebox{-.5\height}{\includegraphics[width=\linewidth]{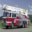}}
	\end{minipage} & \begin{minipage}[b]{0.085\columnwidth}
		\centering
		\raisebox{-.5\height}{\includegraphics[width=\linewidth]{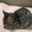}}
	\end{minipage}\\
Classified as&dog (fail.)& airplane & truck (fail.) & horse\\
$L_2$ distortion: &$0.226$ & $0.145$ & $0.131$ & $1.6\times 10^{-3}$\\ \hline

ZOGradOpt  & \begin{minipage}[b]{0.085\columnwidth}
		\centering
		\raisebox{-.5\height}{\includegraphics[width=\linewidth]{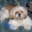}}
	\end{minipage} &
	\begin{minipage}[b]{0.085\columnwidth}
		\centering
		\raisebox{-.5\height}{\includegraphics[width=\linewidth]{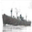}}
	\end{minipage} & \begin{minipage}[b]{0.085\columnwidth}
		\centering
		\raisebox{-.5\height}{\includegraphics[width=\linewidth]{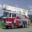}}
	\end{minipage} & \begin{minipage}[b]{0.085\columnwidth}
		\centering
		\raisebox{-.5\height}{\includegraphics[width=\linewidth]{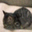}}
	\end{minipage}\\
Classified as&cat & airplane & truck (fail.) & horse\\
$L_2$ distortion: &$0.304$ & $0.254$ & $1.1\times 10^{-30}$ & $0.192$\\ \hline

$\text{ZOSLGH}_{\text{r}}$  & \begin{minipage}[b]{0.085\columnwidth}
		\centering
		\raisebox{-.5\height}{\includegraphics[width=\linewidth]{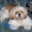}}
	\end{minipage} &
	\begin{minipage}[b]{0.085\columnwidth}
		\centering
		\raisebox{-.5\height}{\includegraphics[width=\linewidth]{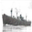}}
	\end{minipage} & \begin{minipage}[b]{0.085\columnwidth}
		\centering
		\raisebox{-.5\height}{\includegraphics[width=\linewidth]{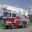}}
	\end{minipage} & \begin{minipage}[b]{0.085\columnwidth}
		\centering
		\raisebox{-.5\height}{\includegraphics[width=\linewidth]{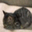}}
	\end{minipage}\\
Classified as&cat & airplane & automobile & horse\\
$L_2$ distortion: &$0.540$ & $0.212$ & $0.282$ & $0.076$\\ \hline

$\text{ZOSLGH}_{\text{d}}$ & \begin{minipage}[b]{0.085\columnwidth}
		\centering
		\raisebox{-.5\height}{\includegraphics[width=\linewidth]{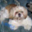}}
	\end{minipage} &
	\begin{minipage}[b]{0.085\columnwidth}
		\centering
		\raisebox{-.5\height}{\includegraphics[width=\linewidth]{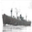}}
	\end{minipage} & \begin{minipage}[b]{0.085\columnwidth}
		\centering
		\raisebox{-.5\height}{\includegraphics[width=\linewidth]{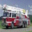}}
	\end{minipage} & \begin{minipage}[b]{0.085\columnwidth}
		\centering
		\raisebox{-.5\height}{\includegraphics[width=\linewidth]{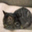}}
	\end{minipage}\\
Classified as&cat & airplane & automobile & horse\\
$L_2$ distortion: &$0.359$ & $0.174$ & $0.241$ & $0.075$\\ 
\bottomrule
\end{tabular}%
\label{tab:cifar10_attack_examples}
\end{table}%

\subsection{MNIST}\label{sec:mnist}
\begin{figure}[H]
    \subfigure[MNIST, Image ID = 7]{
    \label{fig:easy1_mnist}
    \begin{minipage}[t]{0.47\linewidth}
        \centering
        \includegraphics[width=2.5in]{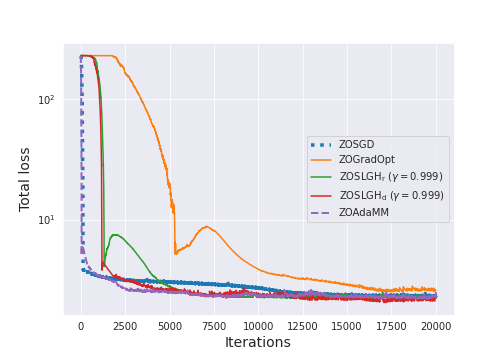}
    \end{minipage}%
    }
    \vspace{-1.2mm}
    \subfigure[MNIST, Image ID = 58]{
    \label{fig:easy2_mnist}
    \begin{minipage}[t]{0.47\linewidth}
        \centering
        \includegraphics[width=2.5in]{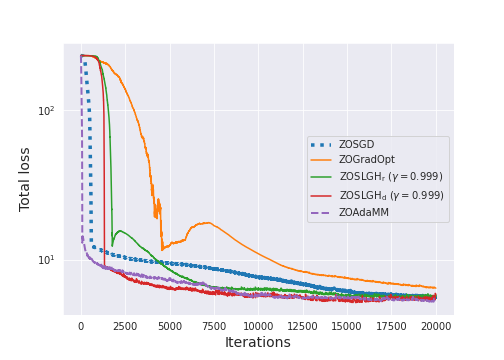}
    \end{minipage}%
    }
    
    \subfigure[MNIST, Image ID = 18]{
    \label{fig:mid1_mnist}
    \begin{minipage}[t]{0.47\linewidth}
        \centering
        \includegraphics[width=2.5in]{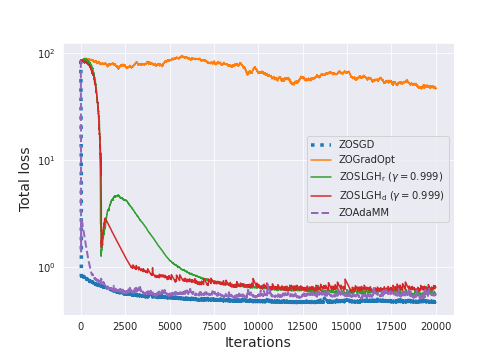}
    \end{minipage}%
    }
    \vspace{-1.2mm}
    \subfigure[MNIST, Image ID = 94]{
    \label{fig:mid2_mnist}
    \begin{minipage}[t]{0.47\linewidth}
        \centering
        \includegraphics[width=2.5in]{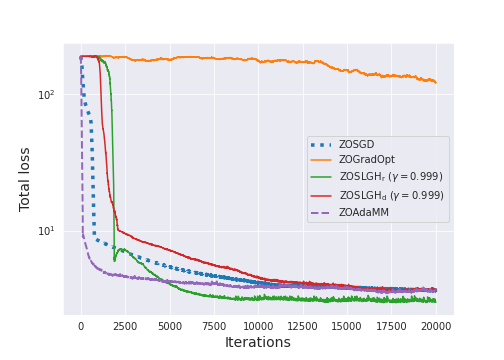}
     \end{minipage}
    }
    \subfigure[MNIST, Image ID = 61]{
    \label{fig:hard1_mnist}
    \begin{minipage}[t]{0.47\linewidth}
        \centering
        \includegraphics[width=2.5in]{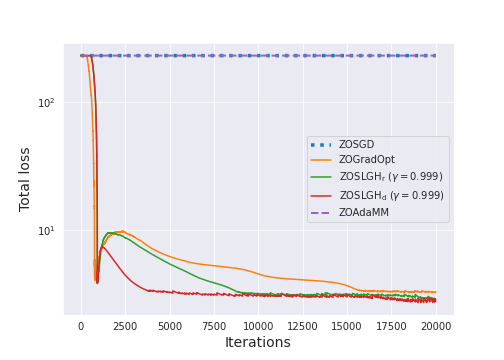}
     \end{minipage}
    }
    \vspace{-1.2mm}
    \subfigure[MNIST, Image ID = 30]{
    \label{fig:hard2_mnist}
    \begin{minipage}[t]{0.47\linewidth}
        \centering
        \includegraphics[width=2.5in]{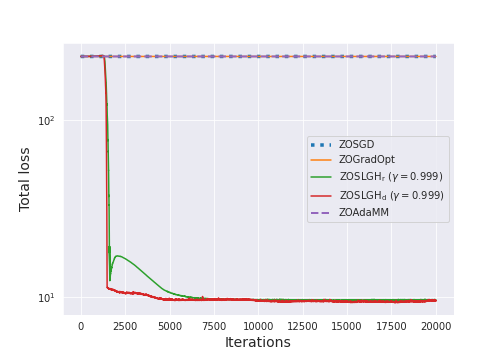}
     \end{minipage}
    }
    \subfigure[MNIST, Image ID = 68]{
    \label{fig:veryhard1_mnist}
    \begin{minipage}[t]{0.47\linewidth}
        \centering
        \includegraphics[width=2.5in]{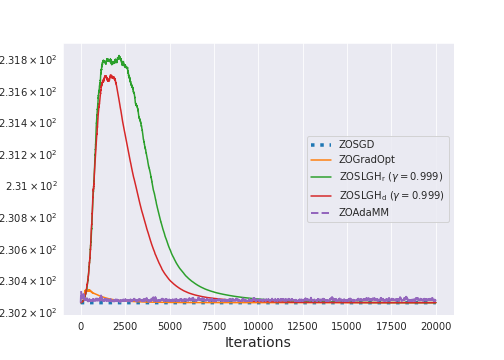}
     \end{minipage}
    }
    \vspace{-1mm}
    \subfigure[MNIST, Image ID = 82]{
    \label{fig:veryhard2_mnist}
    \begin{minipage}[t]{0.47\linewidth}
        \centering
        \includegraphics[width=2.5in]{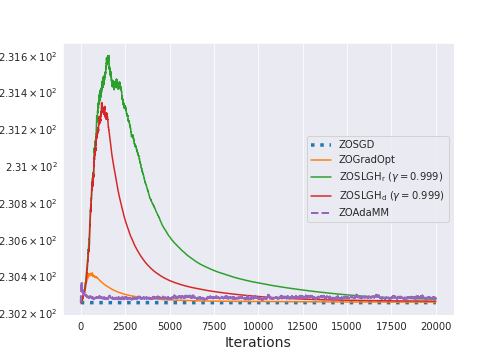}
     \end{minipage}
    }
\centering
\caption{Additional plots of total loss versus iterations on MNIST (log scale). (a)-(b) All algorithms can successfully decrease the loss value when images are easy to attack. (c)-(d) Only GradOpt fails to attack due to its slow convergence. (e) ZOSGD and ZOAdaMM are stuck around a local minimum $x=0$. (f) Only our ZOSLGH algorithms succeed in escaping the local minimum, and thus they can decrease the loss value more than 200 than other algorithms. (g), (h): These images are so difficult to attack that no algorithms can succeed in attacking.}
\label{fig:total_plot_mnist}
\end{figure}

\begin{figure}[H]
    \subfigure[MNIST, Image ID = 61]{
    \label{fig:hard104_l2_mnist}
    \begin{minipage}[t]{0.47\linewidth}
        \centering
        \includegraphics[width=2.5in]{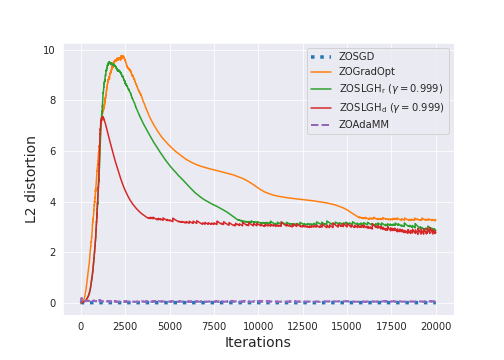}
     \end{minipage}
    }
    \subfigure[MNIST, Image ID = 30]{
    \label{fig:hard39_l2_mnist}
    \begin{minipage}[t]{0.47\linewidth}
        \centering
        \includegraphics[width=2.5in]{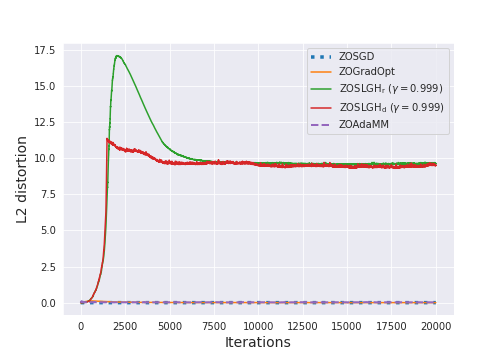}
    \end{minipage}%
    }
    \subfigure[MNIST, Image ID = 68]{
    \label{fig:veryhard14_l2_mnist}
    \begin{minipage}[t]{0.47\linewidth}
        \centering
        \includegraphics[width=2.5in]{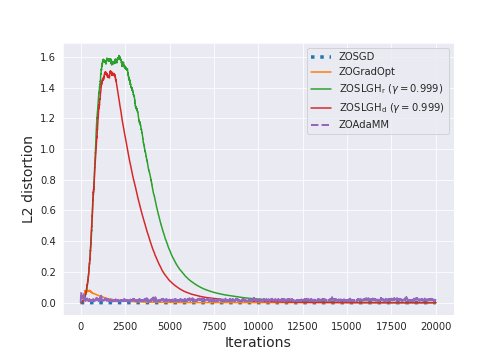}
    \end{minipage}%
    }
    \subfigure[MNIST, Image ID = 82]{
    \label{fig:veryhard41_l2_mnist}
    \begin{minipage}[t]{0.47\linewidth}
        \centering
        \includegraphics[width=2.5in]{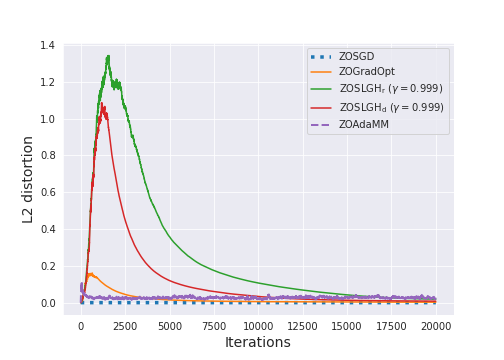}
     \end{minipage}
    }
\centering
\caption{Plots of $L_2$ distortion versus iterations for images that are difficult to attack on MNIST. Each plot of (a)-(d) corresponds to Figure~\ref{fig:hard1_mnist}-Figure~\ref{fig:veryhard2_mnist}.}
\label{fig:l2_plot_mnist}
\end{figure}

Finally, let us show the experimental results on the MNIST dataset. Our ZOSLGH algorithms attain higher success rates than other algorithms on this dataset as well as CIFAR-10 (Table~\ref{table:results_mnist}). Moreover, the average number of iterations to achieve the first successful attack becomes comparable to ZOSGD. The main difference from the results on CIFAR-10 is that the average of $L_2$ distortion at successful time becomes far larger, from $0.050\sim 0.250$ to $4.25\sim 5.20$. This implies that attacks on MNIST are more difficult than those on CIFAR-10. See Figure~\ref{fig:total_plot_mnist} and Figure~\ref{fig:l2_plot_mnist} for additional plots for total loss and $L_2$ distortion. Figure~\ref{tab:mnist_attack_examples} shows adversarial images generated by different algorithms and their original images.

\begin{table}[H]
\centering
  \caption{Performance of a per-image attack over $100$ images of MNIST under $T = 20000$ iterations. ``Succ. rate'' indicates the ratio of success attack, ``Avg. iters to 1st succ.''  is the average number of iterations to reach the first successful attack,  ``Avg. $L_2$ (succ.)'' is the average of $L_2$ distortion taken among successful attacks, and ``Avg. total loss'' is the average of total loss $f(x)$ over 100 samples. Please note that the standard deviations are large since the attack difficulty varies considerably from sample to sample.}
\begin{tabular}{cc|c|c|c|c}
\toprule
\label{table:results_mnist}%
&Methods & \begin{tabular}{c}Succ. rate\end{tabular} & \begin{tabular}{c} Avg. iters\\ to 1st succ.\end{tabular}  & \begin{tabular}{c} Avg. $L_2$\\ (succ.)\end{tabular} & \begin{tabular}{c} Avg. total loss\end{tabular}\\ \hline
SGD algo. &ZOSGD & $67\%$ & $1171 \pm 1954$ & $4.83 \pm 4.13$ & $73.60 \pm 102.70$ \\
&ZOAdaMM& $71\%$ & $\textbf{261} \pm 1068$  & $\textbf{4.25} \pm 3.36$ & $67.49  \pm 100.25$ \\\hline
&ZOGradOpt & $84\%$ & $6166 \pm 4354$ & $5.16 \pm 2.28$ & $28.25 \pm 65.35$ \\ 
GH algo. &$\text{ZOSLGH}_{\text{r}}\ (\gamma=0.999)$ & $\textbf{96\%}$ & $1537 \pm 277$& $\textbf{4.32} \pm 2.44$ & $\textbf{11.83} \pm 37.88$ \\
&$\text{ZOSLGH}_{\text{d}}\ (\gamma=0.999)$ & $\textbf{96\%}$ & $1342 \pm 242$& $\textbf{4.37} \pm 2.58$ & $\textbf{12.09}\pm 38.56$ \\
\bottomrule
\end{tabular}%
\end{table}%

%%%%%%%%%%%%%%%%%%%%%%%%%%%%%%%%%%%%%%%%%%%%%%%%%%%%%%%%%%%%%%%%%%%%%%%%%%%%%%%
%%%%%%%%%%%%%%%%%%%%%%%%%%%%%%%%%%%%%%%%%%%%%%%%%%%%%%%%%%%%%%%%%%%%%%%%%%%%%%%

\begin{table}[H]
  \caption{Comparison of the adversarial images for MNIST with different algorithms.}
\centering
\begin{tabular}{ccccc}
\toprule
Image ID & 10 & 21 & 48& 83\\ \hline
Original  & \begin{minipage}[b]{0.085\columnwidth}
		\centering
		\raisebox{-.5\height}{\includegraphics[width=\linewidth]{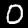}}
	\end{minipage} &
	\begin{minipage}[b]{0.085\columnwidth}
		\centering
		\raisebox{-.5\height}{\includegraphics[width=\linewidth]{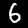}}
	\end{minipage} &
	\begin{minipage}[b]{0.085\columnwidth}
		\centering
		\raisebox{-.5\height}{\includegraphics[width=\linewidth]{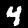}}
	\end{minipage} & \begin{minipage}[b]{0.085\columnwidth}
		\centering
		\raisebox{-.5\height}{\includegraphics[width=\linewidth]{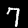}}
	\end{minipage}\\ \ 
Classified as&0 & 6 & 4 & 7\\ 
$L_2$ distortion: &0 & 0 & 0 & 0\\ \hline

ZOSGD  & \begin{minipage}[b]{0.085\columnwidth}
		\centering
		\raisebox{-.5\height}{\includegraphics[width=\linewidth]{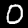}}
	\end{minipage} &
	\begin{minipage}[b]{0.085\columnwidth}
		\centering
		\raisebox{-.5\height}{\includegraphics[width=\linewidth]{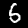}}
	\end{minipage} & \begin{minipage}[b]{0.085\columnwidth}
		\centering
		\raisebox{-.5\height}{\includegraphics[width=\linewidth]{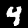}}
	\end{minipage} & \begin{minipage}[b]{0.085\columnwidth}
		\centering
		\raisebox{-.5\height}{\includegraphics[width=\linewidth]{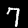}}
	\end{minipage}\\
Classified as&0 (fail.)& 5 & 9 & 7 (fail.)\\ 
$L_2$ distortion: &$4.1\times 10^{-7}$ & $1.194$ & $1.183$ & $1.8\times 10^{-4}$\\ \hline
ZOAdaMM  & \begin{minipage}[b]{0.085\columnwidth}
		\centering
		\raisebox{-.5\height}{\includegraphics[width=\linewidth]{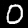}}
	\end{minipage} &
	\begin{minipage}[b]{0.085\columnwidth}
		\centering
		\raisebox{-.5\height}{\includegraphics[width=\linewidth]{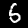}}
	\end{minipage} & \begin{minipage}[b]{0.085\columnwidth}
		\centering
		\raisebox{-.5\height}{\includegraphics[width=\linewidth]{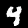}}
	\end{minipage} & \begin{minipage}[b]{0.085\columnwidth}
		\centering
		\raisebox{-.5\height}{\includegraphics[width=\linewidth]{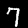}}
	\end{minipage}\\
Classified as&0 (fail.)& 5 & 9 & 7 (fail.)\\ 
$L_2$ distortion: &$4.9\times 10^{-14}$ & $1.334$ & $1.100$ & $4.0\times 10^{-14}$\\ \hline

ZOGradOpt  & \begin{minipage}[b]{0.085\columnwidth}
		\centering
		\raisebox{-.5\height}{\includegraphics[width=\linewidth]{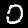}}
	\end{minipage} &
	\begin{minipage}[b]{0.085\columnwidth}
		\centering
		\raisebox{-.5\height}{\includegraphics[width=\linewidth]{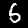}}
	\end{minipage} & \begin{minipage}[b]{0.085\columnwidth}
		\centering
		\raisebox{-.5\height}{\includegraphics[width=\linewidth]{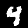}}
	\end{minipage} & \begin{minipage}[b]{0.085\columnwidth}
		\centering
		\raisebox{-.5\height}{\includegraphics[width=\linewidth]{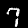}}
	\end{minipage}\\
Classified as& 2&5 & 9 & 9\\
$L_2$ distortion: &$3.898$ & $1.378$ & $1.903$ & $6.379$\\ \hline

$\text{ZOSLGH}_{\text{r}}$  & \begin{minipage}[b]{0.085\columnwidth}
		\centering
		\raisebox{-.5\height}{\includegraphics[width=\linewidth]{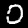}}
	\end{minipage} &
	\begin{minipage}[b]{0.085\columnwidth}
		\centering
		\raisebox{-.5\height}{\includegraphics[width=\linewidth]{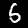}}
	\end{minipage} & \begin{minipage}[b]{0.085\columnwidth}
		\centering
		\raisebox{-.5\height}{\includegraphics[width=\linewidth]{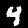}}
	\end{minipage} & \begin{minipage}[b]{0.085\columnwidth}
		\centering
		\raisebox{-.5\height}{\includegraphics[width=\linewidth]{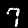}}
	\end{minipage}\\
Classified as& 2&5 & 9 & 9\\
$L_2$ distortion: &$3.867$ & $1.261$ & $1.106$ & $6.075$\\ \hline

$\text{ZOSLGH}_{\text{d}}$ & \begin{minipage}[b]{0.085\columnwidth}
		\centering
		\raisebox{-.5\height}{\includegraphics[width=\linewidth]{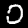}}
	\end{minipage} &
	\begin{minipage}[b]{0.085\columnwidth}
		\centering
		\raisebox{-.5\height}{\includegraphics[width=\linewidth]{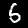}}
	\end{minipage} & \begin{minipage}[b]{0.085\columnwidth}
		\centering
		\raisebox{-.5\height}{\includegraphics[width=\linewidth]{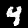}}
	\end{minipage} & \begin{minipage}[b]{0.085\columnwidth}
		\centering
		\raisebox{-.5\height}{\includegraphics[width=\linewidth]{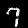}}
	\end{minipage}\\
Classified as& 2&5 & 9 & 9\\
$L_2$ distortion: &$4.048$ & $1.222$ & $1.059$ & $5.722$\\ 
\bottomrule
\end{tabular}%
\label{tab:mnist_attack_examples}
\end{table}%

\end{document}